\documentclass[english]{smfart} \usepackage{amsmath, amsthm, amscd, amsfonts,amssymb} 
\numberwithin{equation}{section} 

\newcommand{\teq}{\arabic{section}.\arabic{equation}}
\newcommand{\teql}{\Alph{section}.\arabic{equation}}

\renewcommand{\theenumi}{\alph{enumi}}

\newcommand{\sqr}[2]{{\vcenter{\vbox{\hrule height.#2pt\hbox{\vrule width.#2pt
height#1pt \kern#1pt\vrule width.#2pt}\hrule height.#2pt}}}}

\newcounter{eqcount}

\renewcommand{\labelenumi}{{{\rm (\teq \alph{enumi})}}} 
\newenvironment{edesc}{\refstepcounter{equation}\begin{enumerate}}%
{\end{enumerate}}
\newenvironment{triv}{\refstepcounter{equation}\begin{list}%
{{\hbox{\rm(\teq)\ }}} \item }{\end{list}}

\newcommand{\ring}[1]{{\mathbb #1}}
\newcommand\bZ{{\ring{Z}}}

\newcommand\bC{{\ring{C}}} \newcommand\bR{{\ring{R}}}
\newcommand\bF{{\ring{F}}} \newcommand\bQ{{\ring{Q}}}
\newcommand\bH{{\ring{H}}}
\newcommand{\csp}[1]{{\mathbb #1}}
\newcommand{\tsp}[1]{{\mathcal #1}}

\newcommand{\prP}{\csp{P}}

\newcommand{\sC}{{\tsp{C}}} 
 
\newcommand{\sQ}{\tsp{Q}}

\newcommand{\sT}{{\tsp {T}}} \newcommand{\sH}{{\tsp {H}}}
 
\newcommand{\sM}{{\tsp {M}}} 
\newcommand{\sG}{{\tsp {G}}}

\newcommand{\eql}[2]{{\rm (\ref{#1}\ref{#2})}} 

\newcommand{\vect}[1]{{\pmb #1}} 
 \newcommand{\bg}{\vect{g}}
 
\newcommand{\bp}{{\vect{p}}} 
\newcommand{\bv}{{\vect{v}}} 
 \newcommand{\bz}{{\vect{z}}}

\newcommand{\row}[2]{{#1_1,\ldots,#1_{#2}}}

\newcommand{\smatrix}[4]{{\big(\begin{array}{cc}
\!\lower2pt\hbox{$\scriptstyle#1$} &\lower2pt\hbox{$\scriptstyle#2$}\!
\\\! \raise2pt\hbox{$\scriptstyle#3$} &\raise2pt\hbox{$\scriptstyle#4$}
\!\end{array}\big)}}

\newcommand{\texto}[1]{{\textr{#1}}}
\newcommand{\GL}{\texto{GL}} \newcommand{\SL}{\texto{SL}}
 \newcommand{\ind}{\texto{ind}}
\newcommand{\PSL}{\texto{PSL}} \newcommand{\PGL}{\texto{PGL}}
 
 \renewcommand{\ni}{\texto{Ni}}
 \newcommand{\Pic}{\texto{Pic}}
\newcommand{\textr}[1]{{\text{\rm #1}}}
\newcommand{\tr}{\textr{tr}} \newcommand{\ord}{\textr{ord}}
\newcommand{\abs}{\textr{abs}}  \newcommand{\cyc}{\textr{cyc}}
 
\newcommand{\pC}{{\textr{C}}} \newcommand{\inn}{\textr{in}}
 \newcommand{\Aut}{\textr{Aut}}

\newcommand{\rd}{\texto{rd}}

\newcommand{\tG}[1]{{}_{#1}\tilde G}

\newcommand{\Det}{\text{\rm Det}}

\newcommand{\GAP}{{\bf GAP}}

\newcommand{\textb}[1]{{\text{\bf #1}}}
\newcommand{\bfC}{{\textb{C}}}
\newcommand{\longmapright}[2]{\smash{\mathop{\hbox to
#2pt{\rightarrowfill}}\limits^{#1}}}
\newcommand{\longmapleft}[2]{\smash{\mathop{\hbox to
#2pt{\leftarrowfill}}\limits^{#1}}}

\newcommand{\mapright}[1]{\smash{\mathop{\longrightarrow}\limits^{#1}}}

   \newcommand{\nm}{{-}}

\newcommand{\lrang}[1]{{\langle #1\rangle}}

\newcommand{\eqdef}{\stackrel{\text{\rm def}}{=}}

\newfont{\sevenrm}{cmr7}
\newfont{\bsevenrm}{cmbx7}
\newfont{\mathseven}{cmsy7}
\newfont{\bigmath}{cmsy10 scaled 1200}
\newfont{\fiverm}{cmr5}
\newfont{\bfiverm}{cmbx5}
\newfont{\hel}{cmbx10 scaled 1200}
\newfont{\eu}{eufb10}
\newfont{\sseu}{eufm5}
\newfont{\seu}{eufm7}
\newfont{\Cal}{cmmib10}
\newfont{\sCal}{cmmib7}

\theoremstyle{plain}
\newtheorem{thm}{Theorem}[section] 
\newtheorem{lem}[thm]{Lemma}
\newtheorem{princ}[thm]{Principle}

\newtheorem{prop}[thm]{Proposition}
\newtheorem{cor}[thm]{Corollary}

\newtheorem{res}[thm]{Result}
\theoremstyle{definition}
\newtheorem{defn}[thm]{Definition}
\newtheorem{exmp}[thm]{Example}
\newtheorem{guess}[thm]{Conjecture}
\newtheorem{quest}[thm]{Question}
\newtheorem{prob}[thm]{Problem}

\theoremstyle{remark}
\newtheorem{rem}[thm]{Remark}

\newcommand{\xs}{\times^s\!}
\newcommand{\wsp}{{$\,$---$\,$}} 

\newcommand{\psigma}{{\pmb \sigma}}

\begin{document}

\newcommand{\sHt}[2]{{\sH(A_{{#1}},\bfC_{3^{#2}})}}
\newcommand{\sHtp}[2]{{\sH_+(A_{{#1}},\bfC_{3^{#2}})}}
\newcommand{\sHtm}[2]{{\sH_-(A_{{#1}},\bfC_{3^{#2}})}}
\newcommand{\sHtpm}[2]{{\sH_{\pm}(A_{{#1}},\bfC_{3^{#2}})}}

\newcommand{\pari}{{\textr{par}}}
\newcommand{\Ct}[1]{{\bfC_{{3^{#1}}}}}
\let\phi=\varphi
\newcommand{\balpha}{{\vect{\alpha}}} \newcommand{\bbeta}{{\vect{\beta}}}
\newcommand{\et}{{\text{\sl et}}}
\newcommand{\ram}{{\text{\rm ram}}}
\newcommand{\G}{{\ring{G}}}
\newcommand{\app}{{App}}
\newcommand{\sV}{{\tsp V}}
\newcommand{\HH}{{{}_+\sH_\infty}}
\newcommand{\NNi}{{{}_+\ni_\infty}} 
\newcommand{\CC}{{{}_+\pC}}
\newcommand{\bfCC}{{{}_+\bfC}}
\newcommand{\mq}{\textr{MQ}}
\newcommand{\one}{{\pmb 1}}
\newcommand{\C}{\text{\rm C}}
\newcommand{\cha}{{\text{\rm char}}}
\newcommand{\lcm}{{\text{lcm}}}
\newcommand{\Out}{{\text{\rm Out}}}
\newcommand{\Sp}{{\text{\rm Sp}}}
\newcommand{\Conj}{{\text{\rm Conj}}}
\newcommand{\ab}{{\text{\rm ab}}}
\newcommand{\hc}{{\text{\rm ${\frac1 2}$-canonical\ }}}
\newcommand{\N}[1]{{\{1,\dots,{#1}\}}} 
\newcommand{\Spin}{{\text{\rm Spin}}}
\setcounter{tocdepth}{3}
\newcommand{\sh}{{\text{\bf sh}}}
\newcommand{\mb}[1]{{\mbox{\boldmath{#1}}}}
\font\eightrm=cmr8  \font\eightit=cmsl8 \let\it=\sl 
\newcommand{\fva}{\cite{FrV2}}
\newcommand{\fval}[1]{\cite[#1]{FrV2}}
\newcommand{\fvb}{\cite{FrV3}}
\newcommand{\fvbl}[1]{\cite[#1]{FrV3}}
\newcommand{\reg}{\text{\rm reg}}
\newcommand{\Gen}[2]{{\text{\rm Gen}(#1,#2)}}
\newcommand{\HMGen}[2]{\text{\rm H-MGen}(#1,#2)}
\newcommand{\HM}{\text{\rm H-M}}
\newcommand{\Pow}{\text{\rm Pow}}
\newcommand{\PowGen}[2]{{\text{\rm PowGen}(#1,#2)}}
\newcommand{\gu}{\text{g-u}}
\newcommand{\mpr}{{\text{\bf mp}}} 
\newcommand{\wid}{{\text{\bf wd}}} 
\renewcommand{\GAP}{{\cite{GAP}}}
\newcommand{\SM}{\text{\rm SM}}
\newcommand{\MT}{{\text{\bf MT}}}
\newcommand{\Fib}{{\text{\bf Fib}}}
\newcommand{\Cu}{{\texto{Cu}}}

\title[Cusps on MTs]{The Main Conjecture of Modular Towers\\ and 
its higher rank generalization} 

\author[M.~D.~Fried]{Michael D.~Fried}

\date{\today}

\address{Math.~Dept., MSU-Billings, Billings MT 59101}
\email{mfried@math.uci.edu, mfri4@aol.com}

\begin{abstract} \begin{small} Publication: In {\sl Groupes de Galois arithmetiques et 
differentiels\/} 
(Luminy 2004; eds. D.~Bertrand and P.~D\`ebes),  Sem.~et Congres, Vol.~{\bf 13} 
(2006), 165--230. \end{small}

The genus of projective curves  discretely separates decidedly different two variable algebraic
relations. So, we can focus on the connected moduli $\sM_g$ of genus $g$ curves. Yet, modern applications      
require a data variable (function) on such curves. The resulting 
spaces are versions, depending on our need from this data variable, of {\sl Hurwitz spaces\/}. A {\sl Nielsen class\/}
(\S\ref{nz.1}) is a set defined by $r\ge 3$ conjugacy classes $\bfC$ in the data variable monodromy $G$. It gives a striking genus analog.  

Using Frattini covers of $G$, every Nielsen class produces a projective system of related Nielsen classes for any prime $p$ dividing
$|G|$. A nonempty (infinite) projective system of  braid orbits in these Nielsen classes is an infinite $(G,\bfC)$ {\sl component (tree) 
branch\/}. These correspond to projective systems of irreducible (dim
$r-3$) components from  
$\{\sH(G_{p,k}(G),\bfC)\}_{k=0}^\infty$, the $(G,\bfC,p)$ Modular Tower (\MT). The classical modular curve towers 
$\{Y_1(p^{k+1})\}_{k=0}^\infty$ (simplest  case: $G$ is dihedral, $r=4$, $\bfC$ are involution classes) are an avatar.

The (weak) Main Conjecture \ref{wconj} says, if $G$ is {\sl $p$-perfect}, there are no rational points at high  
levels of a component branch. When 
$r=4$,   \MT s (minus their cusps) are systems of upper half plane quotients covering the $j$-line.   
Our topics.  
\begin{itemize} \item \S\ref{spinGen} and \S\ref{finerGraphs}: Identifying component branches on a \MT\ from 
{\sl g-$p'$, $p$ and Weigel cusp branches\/} using the \MT\ generalization of {\sl spin structures}.  
\item \S\ref{wcNub}: Listing cusp branch properties that imply the (weak) Main Conjecture and extracting the small list of  
towers that could possibly fail the conjecture.
\item
\S\ref{achSchurDomain}: Formulating a (strong) Main Conjecture for higher rank \MT s (with examples): almost all primes produce a modular
curve-like system.  \end{itemize}
\end{abstract}

\alttitle{La conjecture principale sur les tours modulaires et sa
g\'en\'eralisation en rang sup\'erieur}

\begin{altabstract}
Le genre des courbes projectives est un invariant discret qui permet    
une premi\`ere classification des relations alg\'ebriques en deux
variables. On peut ainsi se concentrer sur les espaces de modules  
connexes ${\sM}_g$ des courbes de genre $g$ donn\'e. Pourtant de
nombreux probl\`emes n\'ecessitent la donn\'ee suppl\'ementaire d'une
fonction sur la courbe. Les espaces de modules correspondants sont les
espaces de Hurwitz, dont il existe plusieurs variantes, r\'epondant \`a 
des besoins divers. Une classe de Nielsen (\S 1) est un ensemble,
constitu\'e \`a partir d'un groupe $G$ et d'un ensemble ${\bf C}$ de
$r\geq 3$ classes de conjugaison de $G$, qui d\'ecrit la monodromie
de la fonction. C'est un analogue frappant du genre.

En utilisant les rev\^etements de Frattini de $G$, chaque classe de
Nielsen fournit un syst\`eme projectif de classes de Nielsen
d\'eriv\'ees, pour tout premier $p$ divisant $|G|$. Un syst\`eme
projectif non vide (infini) d'orbites d'actions de tresses dans ces
classes de Nielsen est une branche infinie d'un arbre de composantes.
Cela correspond \`a un syst\`eme projectif de composantes
irr\'eductibles (de dimension $r-3$) de $\{\sH(G_{p,k}(G),{\bfC})\}_{k=0}^\infty$, la tour modulaire. La tour classique des courbes
modulaires $\{Y_1(p^{k+1})\}_{k=0}^\infty$ (le cas le plus simple o\`u
$G$ est le groupe di\'edral $D_{2p}$, $r=4$ et ${\bf C}$ la classe
d'involution r\'ep\'et\'ee 4 fois) en est un avatar.

La conjecture principale (faible) dit que, si $G$ est $p$-parfait,  
il n'y a pas de points rationnels au del\`a d'un niveau suffisamment
\'elev\'e d'une branche de composantes. Quand $r=4$, les tours
modulaires (priv\'ees des pointes) sont des syst\`emes de quotients 
du demi-plan sup\'erieur au-dessus de la droite projective de
param\`etre $j$. Nous th\`emes.
\begin{itemize} \item \S 3 et \S 4: Identification des branches de composantes sur  
une tour modulaire \`a partir des branches de pointes $g-p^\prime$, $p$ et
Weigel, gr\^ace \`a la g\'en\'eralisation des structures de spin.
\item \S 5: \'Enonc\'e d'un ensemble de propri\'et\'es des branches de pointes
impliquant la conjecture principale (faible) et r\'eduction \`a un
nombre limit\'e de cas de tours pouvant encore \'eventuellement
la mettre en d\'efaut.
\item \S 6: Formulation d'une conjecture principale forte pour des
tours modulaires de rang sup\'erieur (avec des exemples): presque tous  
les premiers conduisent \`a un syst\`eme semblable \`a celui des
courbes modulaires. \end{itemize} 
\end{altabstract}


\thanks{Support from 
NSF \#DMS-99305590, \#DMS-0202259 and \#DMS-0455266. This contains many advances on my March 12, 2004, Luminy talk (subsumed by overheads in
\cite{Ont05}). One of those centers on  Weigel cusps and whether they exist. This interchange with Thomas Weigel 
occurred in Jerusalem and Milan during the long trip including Luminy. Prop.~\ref{op'overgp'} is due to Darren Semmen, a constant
modular representation consultant. Conversations with Anna Cadoret, Pierre Debes and Kinya Kimura influenced me to be more complete than
otherwise I would have been.}  

\vskip -10in \maketitle
\tableofcontents

Luminy in March 2004 gave me a chance to show the growing maturity of Modular Towers (\MT s). Documenting its advances, however, 
uses  two other sources: Papers from this conference; and a small selection from the author's work. 
\S\ref{Lumtalks} lists the former. While the first two papers in that list have their own agendas, they show the influence of \MT s. The last
two papers aim, respectively, at the arithmetic and group theory of \MT s. This paper concentrates on (cusp) geometry.   As \cite{FrBook} is not
yet complete, I've listed typos corrected from the print version of \cite{BFr02}\wsp our basic reference \wsp in the on-line 
version  (\S\ref{BFrtypos}). From it came the serious examples (see partial list of \S\ref{BrFrattProp}) that graphically demonstrate 
the theory.  

A glance at the Table of Contents shows \S\ref{finerGraphs} is the longest and most theoretical in the paper. It will figure in 
planned later papers. We have done our best in \S\ref{achSchurDomain} to get serious examples to illustrate everything in \S\ref{finerGraphs}.
(Constraints include assuring we had in print enough on the examples to have them work as we wanted.) So, we suggest
referring to
\S\ref{finerGraphs} after finding motivation from other sections. 

Many items in 
this paper would seem to 
complicate looking at levels of a \MT: types of cusps, Schur multipliers of 
varying groups, component orbits. It behooves us to have an 
organizing tool to focus, label and display crucial and difficult computations. Further, we
find that arithmetic geometers with little group theory background just don't know 
where to start. What surely helps handle some of these problems is the
$\sh$-incidence matrix. I suggested to Kay Maagard that the braid package (for 
computing Nielsen class orbits) would gain
greatly if it had a sub-routine for this. He said he would soon put such in 
\cite{MShStV}.

We use the $\sh$-incidence Matrix on $\ni(A_4,\bfC_{\pm 3^2})^{\inn,\rd}$ in \S\ref{graph-comp} to show 
what we
mean. More elaborate examples for level 1 of this \MT\ and also for 
$\ni(A_5,\bfC_{3^4})^{\inn,\rd}$ are in
\cite[Chaps.~8 and 9]{BFr02}. All these are done without \cite{GAP} or other computer calculation, and they figure in many places in this paper
as nontrivial examples of the mathematical arguments that describe the structure of \MT\ levels. Still,  
\cite[\S9.2.1 and 9.2.2]{BFr02} list what \cite{GAP} produced for all branch cycles (see \S\ref{learngpcusp} and \S\ref{BrFrattProp}) for
both ($j$-line covering) components at level 1 in the
$(A_5,\bfC_{3^4},p=2)$ \MT. 

\section{Questions and topics} \label{nz.1} In this paper the branch point parameter $r\ge 3$ is usually 4 (or 3). Results (based on
\S\ref{spinGen} and
\S\ref{finerGraphs}) on
\MT s with $r$ arbitrary are in a companion paper
\cite{altGps} that contains proofs of several results from the author's long-ago preprints. For example: It describes all 
components  of Hurwitz spaces attached to
$(A_n,\bfC_{3^r})$, alternating groups with 3-cycle branch cycles running over all $n\ge 3$, $r\ge n-1$. 

\subsection{The case for investigating \MT s} \label{caseMT} A group $G$ and $r$
conjugacy classes $\bfC=\row \C r$ from $G$ define a {\sl Nielsen class\/}  (\S\ref{NielClDict}). The
Hurwitz monodromy group $H_r$  acts on (we say {\sl braids}) elements  in representing Nielsen classes.  Components of
\MT\ levels correspond to $H_r$ orbits. Other geometry, especially related to cusps, corresponds to statements about  subgroups of $H_r$ on
Nielsen classes.  

Sometimes we use the notation $r_\bfC$ for the number $r$ of
conjugacy classes. Mostly, however, we concentrate on \MT s defined by reduced (inner) Nielsen classes
$\ni(G,\bfC)^{\inn,\rd}$ where $r_\bfC=4$  (sometimes one conjugacy class, repeated four times). Then, the sequence of reduced inner Hurwitz
spaces ($\{\sH(G_{p,k}(G),\bfC)^{\inn,\rd}\}_{k=0}^\infty$ below) defining their levels are curves. Here $H_4$, acting on a corresponding
projective sequence of Nielsen classes, factors through a mapping class group we denote as $\bar M_4$. It is naturally isomorphic to
$\PSL_2(\bZ)$. 

In this case, a projective sequence of finite index subgroups of $\PSL_2(\bZ)$  acting on the upper half-plane, indexed by powers of a prime
$p$, do correspond to these levels. Yet,  this sequence appears indirectly in MTs, unlike the classical approach to  the special case of
modular curve sequences.  The closure
$\bar
\sH(G_{p,k}(G),\bfC)^{\inn,\rd}$ is a ramified cover of the $j$-line (\S \ref{j-line}) that includes {\sl cusps\/} (lying over
$j=\infty$). Each cusp identifies with a {\sl Nielsen class cusp set\/} (as in \eql{cusp-orbDef}{cusp-orbDefa}). 

Like modular curves towers, the usual cusp type is a $p$ cusp. Also, like modular curve towers, 
special cusp sets correspond to actual cusps with special geometric properties. The technical theme of this
paper:   \MT s with {\sl g-$p'$ cusps\/}   (\S\ref{cuspTypes})
have a special kinship to modular curves (a subcase). That is because g-$p'$ cusps potentially generalize a classical meaning for those modular
curve cusps akin to representing degenerating Tate elliptic curves. This relates to the topic of  {\sl tangential base points}
(Princ.~\ref{tangBPPrinc} and
\S\ref{secondMC}). The other kind of cusp type called o-$p'$ has no modular curve analog. We give many examples of these occurring on \MT s
where $p=2$ and $G_0$ is an alternating group. 

Direct interpretation of cusps and other geometric properties of
\MT\ levels compensates for how they appear indirectly as upper half-plane quotients. This allows defining \MT s for $r>4$. These
have many applications, and  an indirect relation with Siegel upper half-spaces, though no direct analog with modular curves. 

\subsubsection{Why investigate \MT s?}  We express \MT s as a   
response to these topics.  

\begin{itemize} \item[$T_1$.] They  answer to commonly arising questions: 
\begin{itemize} \item[$T_1$.a.] Why has it taken so long to solve the Inverse Galois Problem?  
\item[$T_1$.b.] How does the Inverse Galois Problem relate to other deep or important problems? \end{itemize} 
\item[$T_2$.] Progress on \MT s generates new applications:  
\begin{itemize} \item[$T_2$.a.] Proving the Main Conjecture shows \MT s have some properties analogous to those for  modular curves.
\item[$T_2$.b.] Specific \MT\ levels have many recognizable applications. \end{itemize} 
\end{itemize}

Here is the answer to $T_1$.a.~in a nutshell. \MT s shows a significant part of the Inverse Galois Problem includes precise
generalizations of many  renown statements from modular curves. Like those statements, \MT\ results say you can't find very many
of certain specific structures over $\bQ$. 

For example, \S\ref{compSTC} cites \cite{cadoret05b} to say  the weak (but not the strong)
Main Conjecture of \MT s follows from the Strong Torsion Conjecture (STC) on abelian varieties.  Still, there is more to say: Progress on
our Main Conjecture implies specific insight and results on the STC (subtle distinctions on
the type of torsion points in question), and relations of it to the Inverse Galois Problem. 

\subsubsection{Frattini
extensions of a finite group $G$ lie behind \MT s} \label{gpStatements} Use the notation $\bZ/n$ for congruences $\!\!\mod n$ and $\bZ_p$
for the $p$-adic integers. Denote the profinite completion of $\bZ$ by $\tilde \bZ$ and its automorphisms (invertible profinite integers) by
$\tilde \bZ^*$.   

Suppose
$p$ is a prime dividing
$|G|$. Group theorists interpret  $p'$ as an {\sl adjective\/} applying to sets related to $G$: A set is $p'$ if $p$ does {\sl not\/} divide
orders of its elements. 

We say $G$ is {\sl $p$-perfect\/} if it has no
$\bZ/p$ quotient. For $H\le G$, denote the subgroup of $G$ generated by commutators ($hgh^{-1}g^{-1}$, $h\in H, g\in G$) by
$(H,G)$. Then, $G$ is {\sl perfect\/} if and only it is $p$-perfect for each $p$ dividing $|G|$ (equivalent to $(G,G)=G$). 
\S\ref{pperfect} explains the point of the $p$-perfect condition. 

A covering homomorphism $\phi: H\to G$ of pro-finite groups is {\sl Frattini\/} if for any proper subgroup $H^*< H$, the image
$\phi(H^*)$ is a proper subgroup of $G$. Alternatively, the kernel $\ker(\phi)$ of $\phi$ lies in the Frattini subgroup
(intersection of all proper maximal subgroups of $G$) of $G$. For
$P$ a pro-$p$ group, the closure of the group containing $p$th powers and commutators is its {\sl Frattini\/} subgroup
$\Phi(P)$. Iterate this $k$ times for $$\Phi^k(P)< \Phi^{k-1}(P)<\cdots < P.$$ 

Consider  a {\sl reduced\/} Nielsen class (\S\ref{redNCcusps}) defined by  $r$ ($p'$) conjugacy classes $$\bfC=(\row {\text C} r)\text{ in
a finite group }G=G_0.$$ Defining  the characteristic (projective) series of Nielsen classes from this requires the  characteristic
(projective) sequence $\{G_k\}_{k=0}^\infty$ of $p$-Frattini covers of $G_0$.  Each
$G_k$ covers
$G$ and is a factor of the universal
$p$-Frattini cover
$\psi: \tG p\to G$, versal for all extensions of $G$ by $p$-groups (\cite[\S1.2]{pierre}, \cite[Chap.~20]{FrJ86}): 
$$\{G_k=G_{p,k}(G)\eqdef
\tG p/\Phi^k(\tilde P_p)\}_{k=0}^\infty\text{ with }
\tilde P_p=\ker(\psi: \tG p\to G).$$  Then, $G_{k+1}\to G_k$ is the maximal Frattini cover of $G_k$
with elementary abelian $p$-group as kernel. Further, $\ker(G_{k+1}\to G_k)$ is a $G_k$ module whose composition factors
consist of irreducible $G_0$ modules. The most important of these is $\one_{G_k}=\one_{G_0}$, the trivial 1-dimensional $G_k$ module.

\cite[\S2.2]{Fr02} shows how to find the rank of the pro-$p$, pro-free  group
$\tilde P_p$. Its subquotients figure in the geometry of the attached \MT\ levels. 

Consider any cover $H\to G$ of profinite groups with kernel $(\ker(H\to G)$ a (pro-)$p$ group. If $\C$ is a $p'$ conjugacy class in $G$, then
above it in $H$ there is a unique $p'$ conjugacy class. This is the most elementary case of the Schur-Zassenhaus Lemma. When we have this
situation it is natural to retain the notation $\C$ for the conjugacy class in $H$, so long as we are clear on which group contains the
class. Conversely, if $\C$ is a $p'$ conjugacy class of $H$ it has a unique image $p'$  conjugacy class in $G$. 

This setup applies whenever we refer to \MT s, as in this. The \MT\ attached to
$(G,\bfC,p)$ is a projective sequence of spaces
$\{\sH(G_k,\bfC)^{\inn,\rd}\eqdef \sH_k\}_{k=0}^\infty$. We also use this  lifting principle even when  $H\to G$ is not a
Frattini cover (as in
\S\ref{gp'Weig}).  

\subsubsection{\MT s and the Regular Inverse Galois Problem} \label{MTview} Use the acronym  RIGP for the Regular Inverse Galois Problem. For
any field
$K$, 
$K^\cyc$ is
$K$ with all roots of 1 adjoined. Let $F$ ($\le \bC$ for simplicity) be a field and 
$G=G_0$ any finite
$p$-perfect group.  An $F$ regular realization of  $G^*$ is a Galois cover $\phi^*:X^*\to \prP^1_z$ over $F$ with group $G^*$ (with 
automorphisms also defined over $F$). Then, the branch point set $\bz$ of $\phi^*$  is an $F$ set, with corresponding
conjugacy classes
$\bfC^*$ in
$G^*$. 

We use the {\sl Branch Cycle Lemma\/} (BCL, \S\ref{cuspBranches};  \cite[Thm.~1.5]{pierre} has example uses when $\bQ=F$).
It says the branch points and respective conjugacy classes satisfy a compatibility condition: For each $\tau\in \Aut(\bC/F)$, 
$z_i^{\tau}=z_j$ implies \begin{equation} \label{bclF} (\C_i^*)^{n_\tau}=\C_j^*\text{ with }\tau\mapsto n_\tau \in G(\bQ^\cyc/F\cap
\bQ^\cyc)\le
\tilde \bZ^* (\S\ref{gpStatements}).\end{equation} We say the conjugacy classes are {\sl $F$-rational\/} if \eqref{bclF} holds without our
having to know anything more about the branch points than they are an $F$ set. That is, if (as a set with multiplicity) $(\bfC^*)^n=\bfC^*$ for
each
$n\in G(\bQ^\cyc/F\cap
\bQ^\cyc)$. 

A significant conclusion is that if $G^*$ is centerless, and $\bfC^*$ is $F$-rational, then such $\phi^*\,$s correspond one-one with 
$F$ points on the space
$\sH(G^*,\bfC^*)^\inn$ (\cite[Thm.~1]{FrV2}; each then gives an $F$ point in $\sH(G^*,\bfC)^{\inn,\rd}$). The quotients of $\tG p$ 
differ in a style akin to the difference between
$D_p$ and
$D_{p^{k+1}}$; in some ways not a big difference at all. So, we ask if they are all regular realizations  
from one rubric? 

\begin{edesc} \label{min-max}\item \label{min-maxa} Minimum: Can all be realized with  some  bound on the number of
branch points (dependent on
$G_0$ and
$p$)?  \item \label{min-maxb}  Maximum: Can all be realized with the same branch point set $\bz$? 
\end{edesc}   For many fields $F$, including number fields  (Rem.~\ref{cycEf}), the hypothesis of
Prop.~\ref{boundBranch} implies its conclusion (\cite[Thm.~2.6]{pierre} outlines the proof). That is, if \eql{min-max}{min-maxa}, then there
is a specific   \MT\  with $F$ points  at each level. 

\begin{prop} \label{boundBranch} Assume there is $r_0$ so each $G_k$ has an $F$ regular realization, with $\le r_0$ branch points. 
Then,  there is a \MT\ from $(G,\bfC)$ with  $r_{\bfC}\le r_0$ and each $\sH(G_k,\bfC)^{\inn}$ (and therefore $\sH(G_k,\bfC)^{\inn,\rd}$), $k\ge
0$, has an
$F$ point.\end{prop}  

The last half answer to Quest.~$T_1$.a is the conjecture that the {\sl conclusion\/} (and therefore the hypothesis) of Prop.~\ref{boundBranch}
doesn't hold for number fields.   

\begin{guess}[Weak Main Conjecture] \label{wconj} Suppose $G_0$
is $p$-perfect and $K$ is a number field. Then, there cannot be $K$ points at every level of a \MT. So, regular
realizations of all the
$G_k\,$s over $K$ requires an unbounded number of branch points. \end{guess} A modular curve case of this is that $Y_1(p^{k+1})$
(modular curve
$X_1(p^{k+1})$ minus its cusps) has no $K$ points for $k>>0$. Thm.~\ref{wmcneedspcusps} says the Main Conj.~holds for $(G_0,\bfC,p)$ unless
there is a $K$ projective sequence of components $\{\sH_k'\subset \sH(G_k,\bfC)^{\inn,\rd}\}_{k=0}^\infty$ and either none of the
$\sH_k'$ has a $p$ cusp; or $\sH_{k+1}'/\sH_k'$ is equivalent  to a degree $p$ rational function   
$f_k:\prP^1_z\to
\prP^1_z$ with 
$f_k$ either a polynomial, or totally ramified over two places.  

\begin{rem}[$F$ for which Prop.~\ref{boundBranch} holds] \label{cycEf} Recall, compatible with
\eqref{bclF}, an element
$g$ in a profinite group is $F$-rational if $g^n$ is conjugate to
$g$ for all $n\in G(\bQ^\cyc/F\cap \bQ^\cyc)\le \tilde \bZ^*$. Denote the  
the field generated by roots of 1 of $p'$ order by   $\bQ^{\cyc,p'}$ and let $F_{p'}=F\cap \bQ^{\cyc,p'}$. 
\cite[Thm.~4.4]{FrK97} shows that if {\sl no\/} $p$-power element $g\in \tG p$ is $F$-rational, then $F$ satisfies
Prop.~\ref{boundBranch}. Further, this holds if
$[F_{p'}:\bQ]< \infty$.  
\end{rem}

\subsubsection{Limit groups} \label{wMC} Finding $F$ regular realizations, and their relation to Conj.~\ref{wconj}, breaks
into three  considerations for the collection of $p$-Frattini covers $G^*\to G$. 

\begin{edesc} \label{bclRub} \item  \label{bclRuba} For what $G^*\,$s is $\sH(G^*,\bfC)^{\inn,\rd}$ nonempty (so it can have
$F$ points)? \item  \label{bclRubb} Which of those nonempty $\sH(G^*,\bfC)^{\inn,\rd}\,$s have some absolutely irreducible $F$ component 
$\sH'(G^*,\bfC)^{\inn,\rd}$? 
\item  \label{bclRubc} Which of the $\sH'(G^*,\bfC)^{\inn,\rd}\,$s have $F$ points. 
\end{edesc} 

Limit groups (a braid orbit invariant) are a profinite summary of what \eql{bclRub}{bclRuba} is
about (\S\ref{limNielsen}): A positive answer for $G^*$ holds in \eql{bclRub}{bclRuba} if and only if $G^*$ is a quotient of a limit group
 for some braid orbit on  
$\ni(G_0,\bfC)$. Note: There may be several limit groups for a given level 0 braid orbit (as in App.~\ref{Z3limgps}). Braid orbits in
$\ni(G_0,\bfC)$ containing g-$p'$ cusps have the whole of
$\tG p$ as one limit group (Princ.~\ref{FP2}). 
\S\ref{gp'WeigC} documents much evidence this is also necessary.  

Fields $F$ that are $\ell$-adic completions of a number field are examples for which the maximum condition \eql{min-max}{min-maxb} holds (see
\cite[\S2.4]{pierre}; though $[F_{p'}:\bQ]= \infty$ in Rem.~\ref{cycEf}). That means there is an 
$F$ component branch (\S \ref{cuspResults} \wsp all levels defined over $F$) on some 
\MT\ with a projective system of $F$ points $\{\bp_k\in \sH(G_k,\bfC)^\inn\}_{k=0}^\infty$. By contrast, though 
\eql{min-max}{min-maxa} (with Prop.~\ref{boundBranch}) postulates $F$ points at all levels of some \MT, over a number field we know they
cannot form a projective system \cite[Thm.~6.1]{BFr02}. 

Denote the completion of a field
$K$ at a valuation $\nu$ of $K$ by $K_\nu$. Evidence from the case of shifts of Harbater-Mumford representatives (H-M reps.) suggests an
affirmative answer for the following. \S\ref{cuspResults} explains the hypotheses opening Quest.~\ref{completionPts}.    

\begin{quest} \label{completionPts} Let $K$ be a number field  with  $\{\sH'(G_k,\bfC)^{\inn}\}_{k=0}^\infty$  a $K$ component branch 
 defined by a g-$p'$ cusp branch. Does it have 
a projective system of $K_\nu$ points  for each $\nu$ over any prime $\ell$ not dividing $|G_0|$? \end{quest}

App.~\ref{F2Z2} and App.~\ref{F2Z3} give cases of Nielsen classes with limit groups other than
$\tG p$. App.~\ref{F2Z2} is a different angle on modular curves, where a universal Heisenberg group obstruction explains the
unique limit group.  

App.~\ref{F2Z3} includes applying Thm.~\ref{highObst} (and Ex.~\ref{highObstExs}). Here, each layer of an H-M cusp branch has above it at least
two components, one not an H-M component. Something similar happens for the main example \MT\ of \cite{BFr02} (for $G=A_5$;
Ex.~\ref{A5limitgps}). So, each level of these examples has at least two components, one with $\tG p$ in its limit group, and the other with
$\tG p$ not in its limit groups.  

A rephrase of \eql{bclRub}{bclRubb} would be to decide which limit groups  produce a $\bQ$ component branch. When the limit group is $\tG
p$ and the component branch is from an H-M cusp branch it is sufficient that all H-M reps.~fall in one braid orbit (see
\S\ref{startOIT}). We expect this to  generalize to g-$p'$ reps. The
criterion of \cite[Thm.~3.21]{Fr4} for H-M reps.~to fall in one braid orbit holds at all levels of a \MT, if it holds at level 0. Still, that
criterion never holds when 
$r=4$, the main case of this paper.  

Finally, given that we know the answers for a particular Nielsen class to \eql{bclRub}{bclRuba} and \eql{bclRub}{bclRubb},
\eql{bclRub}{bclRubc} gets to the nub of our Main Conjecture: High \MT\ levels should have no
rational points over a number field $K$ (at least when the limit group is $\tG p$).    

\S\ref{3rdMC} gives a solid example of how to use the cusp rubric to compute. It shows the nature of the two components, $\sH_0^+\cup
\sH_0^-$  at level 0 of a significant \MT.  Both have genus 0, and
$\sH_0^+$ is an H-M component. Indeed,  it contains  {\sl all\/} H-M cusps (Ex.~\ref{HMbranch}, shifts of special
reps.~in g-$p'$  cusps).  The other has nontrivial lifting invariant (\S\ref{smallInv}) and  so nothing above it at level 1. Both are parameter
spaces of genus 1 curves, and both are upper half plane quotients. Yet,  neither is a modular curve.

\subsection{Five parts on a \MT\ structure} \label{expParts} 
From this point  $r=4$. So,  \MT\ levels are $j$-line covers \cite[Prop.~2.3 and \S2.3.1]{BFr02}.  We list this paper's six main topics.  

\begin{edesc} \label{thisPap} \item  \label{thisPapa} \S\ref{redNCcusps}: Tools for computing cusp widths (ramification orders) and
elliptic ramification of levels. \item \label{thisPapbb} \S \ref{cuspBranches}, \S\ref{cuspTypes} and \S\ref{limNielsen}: Relating infinite 
branches on the cusp and component trees, a classification of cusp types and limit Nielsen
classes.  
\item  \label{thisPapb} \S\ref{gp'Weig}  and \S\ref{weigCor}: Describing infinite component branches.  
\item \label{thisPapbbb} \S\ref{wcNub}: \!Outlining 
for $r=4$ how to prove the \!(weak) Main Conjecture.
\item  \label{thisPapc} \S\ref{initMC}: Formulating the Strong Main Conjecture and comparing its expectations with
that for modular curve towers.  
\item  \label{thisPapd} \S\ref{limNielsen}, \S\ref{secondMC} and \S\ref{3rdMC}: Showing specific \MT\ components apply to
significant Inverse Galois and modular curve topics.  
\end{edesc} 
These contribute to $T_1$.b (\eql{thisPap}{thisPapa}, 
\!\eql{thisPap}{thisPapb} and \eql{thisPap}{thisPapc}) and $T_2$ (\eql{thisPap}{thisPapbb}, \!\eql{thisPap}{thisPapbbb} and
\eql{thisPap}{thisPapd}).  

\subsubsection{Results on cusps} \label{cuspResults} 
Conj.~\ref{genusState} interprets the Main Conjecture as a statement
on computing genera of components. That starts the proof outline that  \eql{thisPap}{thisPapb}  alludes
to. \S\ref{NielsenClasses} turns that computation into group theory and combinatorics. 

Our main results 
relate cusps at a \MT\ level to the components on which they lie. The language uses a cusp (resp.~component) tree
$\sC_{G,\bfC,p}$ (resp.~$\sT_{G,\bfC,p}$) on a \MT\ (\S\ref{CTstart}). The natural map $\sC_{G,\bfC,p}\to \sT_{G,\bfC,p}$ is from containment of
cusps in components. This interprets from a cusp set being in a braid orbit  \eqref{cusp-orbDef}.  

An infinite (geometric) component branch (\S\ref{CTstart}) is a maximal projective sequence 
$$B'=\{\bar \sH'_k\subset  \bar \sH(G_k,\bfC)^{\inn,\rd}\}_{k=0}^\infty\text{ of (geometric) Hurwitz space components}.$$ 

With $F$ a field, call $B'$ an {\sl $F$ component branch\/} if all levels have definition field $F$. An
infinite cusp branch is a maximal projective sequence $$B=\{\bar \bp_k\in \bar \sH(G_k,\bfC)^{\inn,\rd}\}_{k=0}^\infty\text{ of (geometric)
points over }j=\infty.$$  There also exist finite branches, where the last component $\sH_k'$ has nothing above it on
$\sH(G_{k+1},\bfC)^{\inn,\rd}$. Our Main Conjecture only applies to infinite $K$ component branches where $K$ is a number field. Still,
describing the infinite component branches forces  dealing   with the finite branches.  From
\S\ref{redNCcusps}, $B$ corresponds to a sequence of cusp sets defined by a projective
system
$\{{}_k\bg\in  \ni(G_k,\bfC)^{\inn}\}_{k=0}^\infty$ of Nielsen class elements.  Characterizations of such a 
$B$ come from definitions of
$p$, g-$p'$ and o-$p'$  cusps (\S\ref{cuspTypes}).  Three {\sl Frattini Principles\/} \ref{FP1}, \ref{FP2} and \ref{FP3} imply one of these
three happens. 
\begin{edesc} \label{branchtype} \item \label{branchtypea} For $k$ large, $\bar \bp_k$ is a $p$ cusp ($p$ branch).  
\item \label{branchtypeb} For all $k$, $\bar \bp_k$ is a g-$p'$ cusp (g-$p'$ branch).  
\item \label{branchtypec}  For $k$ large, $\bar \bp_k$ is an o-$p'$ cusp (Weigel branch).  
\end{edesc} In  case
\eql{branchtype}{branchtypea} there could be a string consisting of g-$p'$  and/or o-$p'$ cusps before  the $p$ cusp part of the
sequence. For many  g-$p'$ cusps there are no o-$p'$ cusps above them (cusps of shifts of H-M reps., for example as in
Prop.~\ref{op'overgp'}). So, if at level 0 you only have such g-$p'$ cusps, no projective sequence will include both g-$p'$ and o-$p'$ cusps.  

Still, Prop.~\ref{op'overgp'} produces \MT s\
where an o-$p'$ cusp lies over some g-$p'$ cusps at each {\sl high\/} level. When finite exceptional strings don't occur at the start of cusp
branches, we call them {\sl pure}.   Any
\MT\ level can be the start of the tower by applying a fixed shift of the indices. Then these names would apply to cusps at that level. 

\subsubsection{g-$p'$ (cusp) versus Weigel cusp branches} \label{doWeigExist} Any cusp branch $B$ 
determines a component branch $B'$. This allows naming an infinite component branch $B'$ of $\sT_{G,\bfC,p}$ by the name of the cusp branch.  

For example, a g-$p'$ branch (as in Princ.~\ref{FP2}) on the cusp tree produces a g-$p'$ branch on the component tree. A succinct
phrasing of Princ.~\ref{FP2}: 
\begin{triv} Any g-$p'$ cusp starts at least one (infinite) g-$p'$ branch.
\end{triv}  
\noindent A succinct converse of this would help so much to decide which \MT s most resemble modular curve towers. Here is our best
guess for such a converse. 

\begin{guess}[g-$p'$ Conjecture] \label{gp'givesPSC} Show for $K$ a number field, each 
$K$ component branch (\S\ref{cuspResults}) on a \MT\ is defined by
some g-$p'$ cusp branch.  \end{guess} 

Many papers  
consider  H(arbater)-M(umford) cusp (Ex.~\ref{HMbranch}) and component branches (\cite{cadoret05b}, \cite{DDe04}, 
\cite{DEm04}, \cite{WeFM}; not using the term branch). 

By contrast Weigel cusp branches are an enigma.  Identifying g-$p'$ cusps and a corresponding branch of
$\sC_{G,\bfC,p}$ has given the successes for finding infinite branches of
$\sT_{G,\bfC,p}$.  The gist of Conj.~\ref{noWeigB} is they are necessary for a component branch. 
\S\ref{SectnoWeigB} lists evidence for it.  
Examples in  \S\ref{nogp'} show the main issues. 

\begin{guess} \label{noWeigB} With $K$ a number field, there are no Weigel cusp branches on any infinite $K$ component branch of a \MT.
\end{guess} 

If Conj.~\ref{noWeigB} is true, then for any (infinite) $K$ component branch either a g-$p'$ branch defines it or it has only $p$ cusp
branches (see \S\ref{schurQuot}). We also suspect the latter cannot hold,
for such component branches would lack classical aspects.

\subsubsection{Setup for proving the (weak) Main Conjecture} 
The group $H_4$ acts (through $\bar M_4$) compatibly on all Nielsen class levels of a \MT. So any $q\in H_4$ acts on a projective system
$\{{}_k\bg\}_{k=0}^\infty$ defining a cusp branch $B$, with  $\{({}_k\bg)q\}_{k=0}^\infty$ defining a new 
sequence of cusps.  (A different projective system of representatives for $B$ likely gives a different projective system of cusps from
the $q$ action.) 

From
this, many cusp branches may define the same component branch. So any component branch could simultaneously be a g-$p'$,
$p$ and Weigel component branch.      

Thm.~\ref{wmcneedspcusps} says, the (weak) Main
Conjecture \ref{wconj} essentially follows if there must be more than one $p$ cusp branch on a component branch. Since modular
curve towers, and all presently analyzed \MT s have $\infty$-ly many $p$-cusp branches, this seems a sure bet. An affirmative result like 
\cite{BFr02} paved the way if
$\bar\bp_k$ is a
$p$ cusp or the cusp of a shifted H-M rep.  So, here is the hardest remainder (modulo Conj.~\ref{noWeigB}) for \cite{FrWProof}:   
\begin{triv} \label{sumOutlineProof} For $k$ large, a g-$p'$ cusp braids to a $p$ cusp. \end{triv}  \noindent We abstract the framework from 
\cite[\S8]{BFr02} for H-M cusps and $p=2$ in \S\ref{braid-to-p} to show both its likelihood and nontriviality. 

\subsection{\MT s of arbitrary rank and full component branches} \label{appMTRank} For both applications and technical analysis we expand
in two ways on what spaces come attached to a definition of a \MT. 

\subsubsection{Intermediate spaces and groups acting on free groups} \label{actFree} 
Our applications use spaces intermediate to $\sH_k\to U_\infty$ (notation of \S\ref{j-line}), just as modular curves use $Y_0(p^{k+1})$ as a
space intermediate to $Y_1(p^{k+1})\to U_\infty$. This gives the notions of  {\sl full\/} cusp and component graphs
(\S\ref{interSpaces}; these  are rarely  trees).  

Also, starting with a finite group $H$ acting faithfully on a free group $F_u$ (or
a lattice
$\bZ^u$) replacing a finite group
$G$, gives the concept of a \MT\ of rank $u$. This allows running over all primes, not explicitly excluded by our usual assumptions: $G$ is
$p$-perfect and 
 $\bfC$ consists of $p'$ classes. 

We have two immediate reasons for doing this. 
\begin{edesc} \item \S\ref{startOIT}: \!For a version of Serre's {\sl \!O(pen)I(mage)T(heorem)\/} (OIT) 
\cite{SeAbell-adic}.
\item Res.~\ref{mazur-merel}: To compare  \MT s with the most compelling arithmetic
statement we know on modular curve towers. 

\end{edesc}  

\begin{res}[Mazur-Merel] \label{mazur-merel} For each number field $K$, there is a constant $A_K$ (dependent only on $K$) so there are no
rational points on $Y_1(p^{k+1})$ (modular curve $X_1(p^{k+1})$ minus its cusps) if $p^{k+1}> A_K$. \end{res} Our (strong) Main Conjecture
(Conj.~\ref{MMMainConj}) formulates this to \MT s of arbitrary rank.  
\cite{altGps} has applications  to  statements independent of \MT s.  Though \MT\ levels are rarely modular curves (quotients of {\sl
congruence\/} subgroups of
$\PSL_2(\bZ)$ acting on the upper half plane), modular
curve  thinking guides their use. 

\subsubsection{Expanding on cusp and component branches} \label{interSpaces} Using groups intermediate to the $G_k\,$s 
 produces {\sl ($p$)-limit\/} Nielsen classes $\ni(G^*,\bfC)$ with
$G^*$ a maximal quotient (limit group) of $\tG p$ having $\ni(G^*,\bfC)$ nonempty. Limit groups are braid invariants on
projective systems of Nielsen class elements.  Unless $G^*$ is $\tG
p$, these give {\sl full\/} \MT s whose infinite branches don't have components of 
$\{\sH_k\}_{k=0}^\infty$ cofinal among them.  This generalization has three motivations. 

\begin{edesc} \label{hrMot} \item \label{hrMota} To include {\sl all\/} modular curves for $p$ odd in this rubric (not just those closely
related to
$Y_1(p^{k+1})\,$s)  requires a rank 2 \MT. 
\item  \label{hrMotb}  Higher rank
\MT s for 
$p$ can have special {\sl F-quotients\/} (\S\ref{F-q})\wsp still based on the universal
$p$-Frattini cover \wsp from low-level quotients. 
\item \label{hrMotc} Using  \cite{Weigel2} gives some precise limit group properties.    
\end{edesc} 

App.~\ref{F2Z2} gives a full comparison of \MT s with all modular curves. 
It shows the unique limit group for \eql{hrMot}{hrMota} is $(\bZ_p)^2\xs \{\pm 1\}$ (for $p\ne 2$). 
Cor.~\ref{weigCor2} explains how each limit group is defined by a unique obstruction. Here that
obstruction is  universal across all primes, coming from a Heisenberg group. \S\ref{schurQuot} shows how \eql{hrMot}{hrMotc} helps decide
when the limits groups are  
 $\tG p$, the case of our Main Conjecture.  
\S\ref{startOIT} is on  how  F-quotients in \eql{hrMot}{hrMotb} point to generalizations of Serre's OIT.   

\subsubsection{Component branches and Schur multipliers} \label{schurQuot} \cite[\S8]{BFr02} gave a procedure for figuring 
components on a \MT\ level. Making the computations at level 0 requires detailed handling of 
conjugacy classes $\bfC$ for the group $G_0$. Level 0 components in the case of simple groups have contributed much of the success of
the braid approach to the Inverse Galois Problem. Though predicting how components and cusps work at level 0 is still an art,  
various families of groups (simple and otherwise) do exhibit similar patterns when using related conjugacy classes (witness 
$A_n$ and 3-cycles
\cite{altGps}). Given  the level 0 work, we organize for higher levels in three steps.  
\begin{edesc} \label{cuspP} \item \label{cuspPa} Inductive setup from level $k$ to  $k+1$: List cusps at level $k$ within each
braid orbit, and choose one representative
${}_k\bg$ for each braid orbit ${}_kO_b$.
\item \label{cuspPb} List all preimages in $\ni(G_{k+1},\bfC)$ lying over ${}_k\bg$ and use this to list all cusps ${}_{k+1}O_c$ at
level
$k+1$ lying over cusps  ${}_kO_c$ in ${}_kO_b$.  
\item \label{cuspPc} Then, partition cusps lying over ${}_k\bg$ according to their braid orbits.  
\end{edesc} 

The $G_k$ module $M_k=\ker(G_{k+1}\to
G_k)$ controls going from $G_k$ to $G_{k+1}$.  A characteristic sequence of $M_k$ subquotients (called Loewy layers; example
\S\ref{loewyLayers} will help the reader)  are semi-simple
$G_0$ modules. Since \cite{Fr4} we've known it is the  $\one_{G_0}\,$s in the Loewy layers that are critical to  
properties of higher \MT\ levels. 

 The cardinality of the fiber in 
\eql{cuspP}{cuspPb} is a braid invariant. The first business is a version of \eql{bclRub}{bclRuba}: Decide effectively when the
fiber is nonempty. Cor.~\ref{weigCor1} shows it is the $\one_{G_0}\,$s  in
the first Loewy layer of $M_k$ \wsp the maximal elementary $p$ quotient of the Schur multiplier of $G_k$ (\S\ref{pperfect}) \wsp that controls
this.    

Suppose   $O$ is a braid orbit in
$\ni(G=G_0,\bfC)$. Then, $O$ defines a profinite cover
$\psi_O: M_O\to G$ with this versal property (Lem.~\ref{extLem}). For any quotient $G'$ of $\tG p$, each braid orbit $O'\le \ni(G',\bfC)$ over
$O$ corresponds  to   
$\psi': M_O\to G'$ factoring through $\psi_O$. Weigel's Th.~\ref{WeigelThm}  says  $M_O$ is an {\sl oriented $p$-Poincar\'e duality
group\/}. 

One consequence: Cor.~\ref{weigCor1} says that if the fiber over the orbit ${}_kO_b$ is empty (as in \eql{cuspP}{cuspPb}), then 
some $\bZ/p$ quotient in the {\sl first\/} Loewy layer of
$\ker(G_{k+1}\to G_k)$ obstructs it. To wit, if $R\to G$ is the central extension with $\ker(R\to G_k)$ giving this $\bZ/p$
quotient, then  $M_O\to G_k$ for  ${}_kO_b$ does not extend to  $M_O\to R$. 

Further, Cor.~\ref{weigCor2} says that if
$G^*$ is a limit group in a Nielsen class and it is different from
$\tG p$, then the following hold.
\begin{edesc} \label{obstProp} \item 
$G^*$ has exactly one nonsplit extension by a $\bZ/p[G^*]$ module $M'$. 
\item $M'$ is the trivial (one-dimensional) $\bZ/p[G^*]$ module. \end{edesc} App.~\ref{F2Z2} and \ref{F2Z3} 
give explicit examples identifying $M'$. 

The example of \S\ref{3rdMC} combines 
the {\sl  $\sh$-incidence matrix\/} with the natural division into cusp types from
\S\ref{cuspTypes} to show how we often manage figuring \eql{cuspP}{cuspPc}. Princ.~\ref{FP3} frames in pure group theory how to
deal with o-$p'$ cusps. So, it sets a module approach  for, say, Conj.~\ref{noWeigB}. Here's how this refined tool  
relates cusps with their components.  

Suppose $\bg\in O\le \ni(G_k,\bfC)$ defines an o-$p'$ cusp. Then,   
having an o-$p'$ cusp $\bg'\in O'\in \ni(G_{k+1},\bfC)$ over $\bg$ restates as a versal property for two profinite groups extensions that induce
$\psi_{O'}:M_{O'}\to G_{k+1}$. This characterizes with group theory whether there are Weigel cusp branches through $O$. These formulas will
generalize to
\MT s of arbitrary rank and any value of
$r$.  

\subsection{Generalizing complex multiplication and Serre's OIT} \label{startOIT} App.~\ref{F2Z3} gives a significant example when
there are several limit groups $G^*$ (one, at least, $\not =\tG p$) and \wsp as we show \wsp the  spaces are not modular curves. So, it is
nontrivial that we can here be explicit in formulating an OIT and a \MT\ version of complex multiplication.  

\subsubsection{Decomposition groups} Suppose
$j'\in U_\infty(F)$ (\S\ref{j-line}; with
$F$ a number field) is a $j$ value. Then, there is a decomposition group $D_{j'}$  from $G_F$ acting on projective systems of points
$\Fib_{j'}(G^*,\bfC)$ on the full \MT\ over $j'$ defined by $(G^*,\bfC)$. \cite[Thm.~6.1]{BFr02}  (when $G^*=\tG p$) 
says no orbit has length one. 
It is far stronger than the Main Conjecture to have 
$D_{j'}$  with {\sl large\/} orbits on $\Fib_{j'}(G^*,\bfC)$, for all $j'$. 

To go, however, beyond naivet\'e requires estimating how large $D_{j'}$ is. Lem.~\ref{GFaction}  
explains how to use cusp branch types: Practical knowledge of how $G_F$ acts on systems of components comes from knowing how $G_F$
acts on specific types of cusps. 

The historical example is where we know all H-M cusps fall in one braid orbit. Then, 
\cite[Thm.~3.21]{Fr4} says a component containing the H-M cusps has definition field given by the BCL (\S\ref{MTview}; this is $\bQ$ if 
$\bfC$ is $\bQ$-rational). 
\cite{cadoret05} exploits this for arbitrary $r$ to produce many Nielsen classes where the corresponding reduced Hurwitz space contains
absolutely irreducible curves over $\bQ$ (the first result of its kind).  

We expect  
g-$p'$ cusp {\sl types\/} to be the main tool for many results. For example,   
\cite[Thm.~3.21]{Fr4} should generalize to describe component branches with all levels defined over some fixed number field. We guess this is  
 exactly when all g-$p'$ cusps of a fixed {\sl type\/} fall in a bounded (independent of the level)  number of orbits. 

Here is another example. \S\ref{comp-tang} notes that 
a  g-$p'$ cusp branch 
$B$ provides a {\sl tangential base-point} in the sense of Nakamura. Related cusps  would allow
following the proof of Serre's OIT for
\lq\lq large\rq\rq\
$j$-invariant, by considering the arithmetic of these cusps over all rank 1 complete fields.

\subsubsection{Seeking OIT examples} 
App.~\ref{F2Z2} has a (rank 2, as in \S\ref{actFree}) \MT\ attached to $F_2\xs\,\bZ/2$. It describes the full \MT\ whose levels identify with
standard modular curves. Here, for all (odd) $p$, there is a unique limit \MT, and a unique
(proper)  F-quotient of it. For each there is a (full) component graph, which we respectively denote by $\sT_{\GL_2}$ and $\sT_{\text{CM}}$. 

So,
in this language, we expect $j'$ values that produce decomposition groups that correspond to $\sT_{\GL_2}$ (or to
$\GL_2$)  and to $\sT_{\text{CM}}$ (or to $\text{CM}$). That this is so is  Serre's OIT, in our
language. Our next example shows how to extend this to general higher rank \MT s. Seeking an OIT type
result uses analog  properties from Serre's example. It is crucial that we expect there to be Frattini properties for monodromy
groups of \MT\ component branches, as in
\eqref{serreProp}.  

App.~\ref{F2Z3} has a rank 2 \MT\ attached to $G=F_2\xs\, \bZ/3$ that 
shows possibilities for general results like Serre's OIT. We see the g-$p'$ cusp
criterion (Princ.~\ref{FP2}) for identifying  infinite component branches in a \MT. For both $p=2$ and
$p\equiv -1 \mod 3$, one limit group is $\tG p=\tilde F_{2,p}\xs \bZ/3$, and its \MT\ has no F-quotient. At least for $p=2$, there are other
limit groups, explicitly showing Cor.~\ref{weigCor1}. We conjecture $D_{j'}$ in these cases always has
a  type we call $F_2$.

For $p\equiv +1
\mod 3$, 
$\tilde F_{2,p}\xs \bZ/3$ is also a limit group, but its \MT\ has a unique F-quotient.  In this case we expect $D_{j'}$  has either type $F_2$
or  a type we call
$\text{CM}$ (and both types occur).  

\subsubsection{Low \MT\ levels apply to the RIGP and to Andre's Theorem} 
\eql{thisPap}{thisPapd}
alluded to the specific applications of its level 0 and 1 components for $p=2$. None of its levels are modular curves. Also, unlike
modular curve levels, these levels have several components. \S\ref{3rdMC} labels the two level 0 components as $\sH_0^+$ and $\sH_0^-$. Level 1
has six, labeled $\sH_1^x$ with the $x$  decoration  signifying some special property. Here appear generalizations of spin invariants
(as in
\S\ref{schurQuot}) that produce varying types of component branches. 

For
$p=2$, and level 0, $\sH_0^{\pm}$ (parametrizing families of genus 3 curves) map to their absolute (reduced) Hurwitz space
versions
$\sH_0^{\pm,\abs}$. Each, like a modular curve, parametrizes genus 1 curves with extra structure and
embeds naturally in
$\prP^1_j\times
\prP^1_j$. 

Suppose in this embedding the
components have infinitely many coordinates in complex quadratic extensions of $\bQ$. Then, we might be suspicious when $p=2$ that this \MT\
would have some complex multiplication property.  A theorem, however, of Andr\'e's (Prop.~\ref{appAndre}) says they don't. This
further corroborates our guess that for $p=2$, almost all $D_{j'}$ have type $F_2$. 

Two level 1 components (\S\ref{level1p2}) contain H-M reps. We show what a serious challenge is deciding whether their defining field is $\bQ$,
with its  effect on the RIGP  (applied to the exponent 2 Frattini cover of $A_5$).

\section{Ingredients for a \MT\ level} We start with some notation and an explanation of how Schur multipliers appear here. Then we briefly try
to comfort a reader about Hurwitz spaces as families of covers of the Riemann   sphere: $\prP^1_z=\bC_z\cup\{\infty\}$. 

\subsection{$p$-perfectness and Schur multipliers} \label{pperfect} Consider 
$r$ conjugacy classes, $\bfC$, in $G$ and $\bg=(\row g r)\in G^r$. Then, $\bg\in \bfC$  means $g_{(i)\pi}$ is in $\C_i$, for some $\pi$
permuting $\{1,\dots, r\}$. Also, 
$\Pi(\bg)\!\eqdef\! \prod_{i=1}^r g_i$ (order matters). Lem.~\ref{perfnongen} shows how $p$-perfect enters.  

\begin{lem} \label{perfnongen} If  $p$ is a prime with $G$
not $p$-perfect and $\bfC$ are $p'$ classes of $G$, then  elements in $\bfC$ are in the kernel of $G$ to 
the corresponding $\bZ/p$ quotient. So, if $\bg\in \bfC$ then
$\lrang{\bg}=G$ is impossible: $\ni(G,\bfC)$ (and the Hurwitz space) is empty. \end{lem}

Here is another technical plus from the $p$-perfect condition. There is a Frattini cover $R_p\to G$ with $\ker(R_p\to G)$
in the center of $R_p$ and equal to  the $p$ part of the Schur multiplier of $G$. Further,  $R_p\to G$ is {\sl
universal\/} for central
$p$ extensions of $G$  (for example,
\cite[\S3.6.1]{BFr02}; call it the representation cover for $(G,p)$). We use the notation $\SM_G$ (resp.~$\SM_{G,p}$) for the Schur
multiplier (resp.~$p$-part of the Schur multiplier) of
$G$.  If $G$ is $p$-perfect for all $p||\SM_G|$, then the fiber product over
$G$ of all such $R_p$ is truly a universal Frattini central extension of $G$. \S\ref{schurDisc} lists properties we use of Schur multipliers.

Identifying components of \MT\ levels is a recurring theme. Whether a component at level $k$ has some component above it at level
$k+1$ \wsp the level $k$ component is {\sl unobstructed\/} \wsp is controlled by Schur multipliers. Lem.~\ref{liftLem}
and Cor.~\ref{weigCor1} are our main tools. Applying them is the heart of describing the  type of infinite branches in a
\MT. We conclude with comments on the literature. 

The definition of homology groups of $G$ (with coefficients in $\bZ$)  came from topology. These were the homology groups of a space with
fundamental group $G$ whose simply connected cover is contractible. \cite[p.~2]{brown} discusses  how Hopf used it to
describe 
$H_2(G,\bZ)$. Write $G=F/R$ with $F$ free. Then, $H_1(G,\bZ)=G/(G,G)$ and $H_2(G,\bZ)=R\cap (F,F)/(F,R)$ (the Schur multiplier of
$G$). 

The expression for $H_1$ is from general principles. For $H_2$ it is  not obvious.  It is usual to compute $H_2$ using tricks to identify 
$E$  that suits \eqref{schurProp}. If $G$ is perfect, then there is a universal (short exact) sequence \begin{equation}
\label{schurProp} 0\to H_2(G,\bZ)\to E\to G\to 1.\end{equation} The group $E$ factors through all central extensions of $G$ \cite[p.~97,
Ex.~7]{brown} (by a unique map through $p$ group extensions if $G$ is $p$-perfect). By contrast, the universal Frattini cover $\tG {}\to G$ of
$G$ is versal: It factors through all extensions of
$G$ including
$E$, but the factoring map isn't unique.  Then,  $R_p\to G$ is the extension
of $G$ from modding
$E$ out by the $p'$  part of $\ker(E\to G)$. It is easy that $p$-perfectness is the same as  $R_p\to G$ being a universal 
$p$-central extension of $G$. 

Also, if $G$ is $p$-perfect and centerless, then all the characteristic Frattini quotients (\S\ref{gpStatements})  $G_{p,k}$ are too. That
implies 
$\sH(G_{p,k},\bfC)^\inn$ (see below) has {\sl fine\/} moduli \cite[Prop.~3.21]{BFr02}. Take $R_{p,k}$ as the representation cover of
$(G_{p,k},p)$. Then,    
$\sH(R_{p,k},\bfC)^\inn$ does not have fine moduli. Both statements produce many Hurwitz space applications.  

\subsection{One cover defines a family of covers} \label{familyDes} An analytic cover, $\phi:
X\to
\prP^1_z$ of compact Riemann surfaces, ramifies over a finite set of points 
$$\bz=\row z r\subset 
\prP^1_z:\ \prP^1_z\setminus \{\bz\}=U_\bz.$$   
Such a  $\phi$ defines a system of
covers by applying Riemann's existence theorem and deforming the branch points (keeping them distinct).  We explain. 

Represent projective $r$ space $\prP^r$ as
nonzero polynomials of degree at most $r$ modulo scalar multiples. Then, polynomials ($r$ unordered points) with at
least two equal zeros form its {\sl discriminant\/} locus $D_r$. Denote $\prP^r\setminus D_r$ by 
$U_r$. By moving  branch points $\bz$, you can form along any path in $U_r$ a unique continuation
of the cover $\phi$. 

Given $\bz$ and classical generators at $\bz^0$ (\cite[\S2.1-2.2]{BFr02} or
\S\ref{gp'Weig}), this interprets  homotopy classes of paths in
$\pi_1(U_r,\bz)$ as Hurwitz monodromy 
$H_r$ (\S\ref{NielClDict}). Its action on Nielsen classes then  reproduces this deformation of covers.  

Suppose given $(G_0,\bfC,p)$ with $p'$ classes $\bfC=(\row {\C} r)$. \cite[\S1.2]{pierre} reminds how this produces a projective
sequence 
$\{\sH_k^\inn\}_{k=0}^\infty$,
 of inner Hurwitz spaces. Assuming it is nonempty, the level $k$ space has dimension $r$ and is an affine variety \'etale
over
$U_r$. These levels correspond to inner Nielsen classes as in
\S\ref{NielsenClasses}. 

Any $\bp\in \sH_k^\inn$ corresponds to an equivalence class of Galois covers $\phi_\bp:X_\bp\to \prP^1_z$, with group
denoted $\Aut(X_\bp/\prP^1_z)$.  The representative includes a specific isomorphism $\mu: \Aut(X_\bp/\prP^1_z)\to G_k(G)$. Another cover
$\phi':X'\to
\prP^1_z$ is in the same {\sl inner\/} class if the following holds. There is a continuous $\psi: X'\to X_\bp$, commuting with the maps to
$\prP^1_z$, inducing conjugation  by some  
$g\in G_k(G)$ between identifications of $\Aut(X_\bp/\prP^1_z)$ and $\Aut(X'/\prP^1_z)$ with $G_k(G)$. We say the cover is in the Nielsen
class
$\ni(G_k(G)=G_{p,k}(G),\bfC)^\inn$. 

More detail is in \cite[\S2]{BFr02}, \cite[Chap.~4]{FrBook}, \cite[Chap.~10]{VB}. 
The first two especially discuss the motivation and basic definitions for \MT s. 

\subsection{Reduced inner spaces} \label{j-line} We use {\sl reduced\/} inner Nielsen
classes. This references triples $(\psi,\mu, \beta)$ (not just $(\psi,\mu)$ as in \S\ref{familyDes}): $\beta\in \PGL_2(\bC)$, and
$\phi_\bp\circ\psi=\beta\circ
\phi'$. 

\subsubsection{The $j$-invariant} \label{jinv} To an unordered 4-tuple $\bz\in U_4$ we associate the {\sl $j$-invariant\/} $j_\bz$ of $\bz$, a
point of
$U_\infty\eqdef \prP^1_j\setminus
\{\infty\}$. To simplify, normalize  so $j=0$ and 1 are the usual elliptic
points corresponding to $j_\bz$ having non-trivial (more than a Klein 4-group; \S\ref{NielsenClasses} and  \cite[Chap.~4,
\S4.2]{FrBook}) stabilizer in
$\PGL_2(\bC)$.

Given $j'\in U_\infty \setminus \{0,1\}$, there is an uncanonical one-one association: covers with $j$-invariant $j'$ in the reduced Nielsen
class $\Leftrightarrow$ elements of the reduced Nielsen classes (\S\ref{redNCcusps}).  So, reduced Nielsen classes produce   
$\{\sH_k=\sH(G_k(G),\bfC)^{\inn,\rd}\}_{k=0}^\infty$: a projective sequence 
 of inner reduced Hurwitz spaces. 

The map $\sH_{k+1}\to \sH_k$ is a cover over every unobstructed component
(\S\ref{pperfect}) of $\sH_k$. By cover we include that it is  possibly ramified for $k$ at points  over $j=0$ or 1. Each {\sl
nonempty\/} component of  
$\sH_k$ is    an upper half-plane quotient and $U_\infty$  
cover (ramified only over $j=0$ and 1) \cite[\S 2]{BFr02}. 

Since the components of $\{\sH_k\}_{k=0}^\infty$ are curves, they have natural
nonsingular projective closures $\{\bar \sH_k\}_{k=0}^\infty$, with each $\bar \sH_k$ extending to give a finite map to $\prP^1_j$. As
expected, we call the (geometric) points of $\bar \sH_k\setminus
\sH_k$ the level $k$ {\sl cusps\/}. 
 
To see why we use reduced spaces consider the following statement (encapsulating \eql{apps}{appsb}) where {\sl
$\infty$-ly many\/} means no two are  reduced equivalent. 
\begin{triv} \label{inequivEx} For there to be  
$\infty$-ly many 4 branch point, reduced inequivalent $\bQ$  regular realizations of $G_1(A_5)$, the H-M components of
$\sH(G_1(A_5),\bfC_{\pm 5^2})^{\inn,\rd}$ must have infinitely many
$\bQ$ points. \end{triv} The (two) H-M components in question have genus 1. We ask if they have infinitely many $\bQ$ points. 
Even one $\bQ$ point $\bp$ (not a cusp) on one of these components would give a geometric cover $\phi_\bp:X_\bp\to
\prP^1_z$ over $\bQ$ with group $G_1(A_k)$. Further, running over $\beta\in \PGL_2(\bQ)$ the covers $\beta\circ\phi_\bp:X_\bp\to
\prP^1_z$ give $\infty$-ly many inner inequivalent covers with the same group also over $\bQ$. These, however, are all reduced equivalent. It
is more significant to consider the outcome of  \eqref{inequivEx}.

The following statement implies Conj.~\ref{wconj}   (special case of \cite[Thm.~6.1]{BFr02}; outline in
\cite[Thm.~2.6]{pierre}).

\begin{guess} \label{genusState}  For large
$k$, all components of $\bar \sH_k$ have genus exceeding 1. \end{guess} 

\subsubsection{Definition fields} \label{defFields} 
All \MT\ levels, with their moduli space structure,  have minimal definition field the same common cyclotomic field (\S\ref{MTview}). If $\bfC$
is $\bQ$-rational, then this definition field is $\bQ$. Still, it is the absolutely irreducible components of 
levels that require attention. For example, if our base field is $\bQ$, and some \MT\ level has no $\bQ$ components, then
this (or any higher) level can have no $\bQ$ points.  This case of the weak Main Conjecture is then trivial (for $\bQ$). 

\S\ref{comp-tang} reminds of methods to find \MT s with component branches over $\bQ$.  They don't,
however, apply when $r_\bfC=4$. So, some  component branch of a MT might have no number field 
definition: No matter what is $K$ with $[K:\bQ]<\infty$, there may be a value of $k$ so the level $k$ component  has
definition field outside $K$.  Lem.~\ref{GFaction}  uses cusp branches to limit, though not yet eliminate, this possibility. Thus, our approach
to the Main Conjecture  aims at deciding it  based only on the \MT\ (cusp) geometry. 

\subsection{Nielsen classes, Hurwitz monodromy and computing genera}  \label{NielsenClasses}

We can compute the genera of
the components of $\bar \sH_k$ using the Riemann-Hurwitz formula by answering the following questions.  
\begin{edesc} \label{cuspData} \item \label{cuspDataa}  What are the $\bar \sH_k$ components.
\item \label{cuspDatab} What \!are the cusp\! widths  \!(ramification \!orders \!over \!$\infty$) in \!each
component.
\item \label{cuspDatac} What points ramify in each component over elliptic
points ($j=0$ or 1). \end{edesc} 

\subsubsection{A Nielsen class dictionary} \label{NielClDict} Use notation of
\S\ref{pperfect}.  {\sl Reduced\/} Nielsen classes let us calculate
components, cusp and  elliptic ramification. We'll see how the Frattini property
controls growth of cusp widths (ramification) with $k$. 

Here are  definitions of Nielsen classes, and their absolute (requires adding a
transitive permutation representation $T: G\to S_n$) and inner quotients. In the absolute case we equivalence Nielsen class elements $\bg$ and
$h\bg h^{-1}$ with $h$ in the normalizer $N_{S_n}(G)$ of $G$ in $S_n$.

$$\begin{array}{rl} \text{Nielsen classes:}&\ni(G,\bfC) = \{\bg\in \bfC\mid \lrang{\bg}=G;
\Pi(\bg)=1\} \\
\text{Absolute classes:}&   \ni(G,\bfC)/N_{S_n}(G,\bfC) \eqdef\ni(G,\bfC,T)^\abs; \text{ and}  \\
\text{Inner  classes:}&
\ni(G,\bfC)/G  \eqdef\ni(G,\bfC)^\inn.\end{array}$$ 

Elements
$q_i$, $i=1,2,3$ (braids), generate the degree 4 {\sl Hurwitz monodromy group\/} $H_4$. Each acts on any Nielsen classes by a twisting on
its 4-tuples. Example:  
$$ \label{twisting} q_2: \bg\mapsto 
(\bg)q_2=(g_1,g_2g_3g_2^{-1}, g_2,g_4). 
$$  

For $\beta\in \PGL_2(\bC)$, reduced equivalence of covers (as in \S\ref{j-line}) works as follows: 
$$\phi: X\to \prP^1_z \Longleftrightarrow 
\beta\circ \phi: X\to \prP^1_z.$$ This equivalence preserves the $j=j_\bz$-invariant of the  branch point 
set 
$\bz=\bz_\phi$.

Reduced equivalence on Nielsen
classes results from each set $\bz_\phi$ having some Klein 4-group  subgroup of 
$\PGL_2(\bC)$ fixing it.  This corresponds to 
modding out the Nielsen class by $\sQ''=\lrang{(q_1q_2q_3)^2, q_1q_3^{-1}}\le H_4$ \cite[Prop.~4.4]{BFr02}. 

So, the action of $H_4$ on reduced Nielsen classes factors through the {\sl mapping class group}: $\bar M_4\eqdef H_4/\sQ''\equiv
\PSL_2(\bZ)$. \cite[\S2.7]{BFr02} has normalized this identification with $\PSL_2(\bZ)$ (see \S\ref{redNCcusps}). It uses generators  
\begin{equation} \label{gammas} \begin{array}{rl} & \lrang{\gamma_0,\gamma_1,\gamma_\infty},  
\gamma_0=
q_1q_2,\gamma_1=\sh= q_1q_2q_3,\gamma_\infty = q_2, \\ &\text{satisfying the product-one relation: }\gamma_0\gamma_1\gamma_\infty=1.
\end{array}\end{equation}  

\subsubsection{Reduced Nielsen classes and cusps} \label{redNCcusps} Regard the words $\gamma_0,\gamma_1,\gamma_\infty$ in the
$q_i\,$s of \eqref{gammas}  as in $H_4$. Usually the
$\gamma$ notation expresses them as  acting in the quotient group
$\bar M_4$, on reduced Nielsen classes. 

Here is the notation for absolute (resp.~inner) reduced representatives:    
$$\begin{array}{rl} \ni(G,\bfC)/\lrang{N_{S_n}(G,\bfC),\sQ''}\eqdef &\ni^{\abs,\rd}
\text{ and}\\
\ni(G,\bfC)/\lrang{G,\sQ''}\eqdef &\ni^{\inn,\rd}.\end{array} $$ 

The element $\sh$ acts like the {\sl shift\/}. It sends a reduced rep.~$\bg=(\row g 4)$ to the reduced class of $(g_2,g_3,g_4,g_1)$.
On reduced Nielsen classes, $\sh$ has order 2 (not 4 as it does on Nielsen classes). Similarly, $\gamma_0$ has order 3 on
reduced Nielsen classes (absolute or inner). Yes, these identify with the generating elements in $\PSL_2(\bZ)$ having orders 2 and 3
corresponding respectively to
$j=1$ and
$j=0$! The action of $\gamma_\infty=q_2$ then gives a combinatorial interpretation of cusps. 

\begin{defn} \label{cuspGp} The  {\sl cusp group\/} (a subgroup of $H_4$) is  
$\text{Cu}_4=\lrang{q_2,\sQ''}$. \end{defn} 

Orbits of $\text{Cu}_4$ (resp.~$\bar M_4$) on  Nielsen classes correspond to cusps (resp.~components) of the corresponding Hurwitz spaces
\cite[Prop.~2.3]{BFr02}. In computational notation, running over
$\bg\in
\ni(G_k,\bfC)^{\inn,
\rd}$: 
\begin{edesc} \label{cusp-orbDef} \item \label{cusp-orbDefa} Cusps on $\bar \sH_k$ $\Leftrightarrow$
$(\bg)\text{Cu}_4$, a {\sl cusp set\/} in the Nielsen classes. 
\item \label{cusp-orbDefb} Components on $\bar \sH_k$ $\Leftrightarrow$ $(\bg)\bar M_4$, a braid orbit on Nielsen classes.
\end{edesc} 
We often refer to $\bg\in \ni(G,\bfC)$ as a cusp, shortening reference to its cusp set. 

\subsubsection{Riemann-Hurwitz on components} 
Now we interpret Riemann-Hurwitz: $(\gamma_0,\gamma_1,\gamma_\infty)$ act on a $\bar M_4$ orbit
$\Leftrightarrow$  branch cycles for a component of $\bar
\sH(G,\bfC)^\rd
\to
\prP^1_j$. 

\begin{edesc} \label{cuspct} \item  \label{cuspcta}  Ramified points over 0 $\Leftrightarrow$ orbits of
$\gamma_0$.
 \item  \label{cuspctb}  Ramified points over 1 $\Leftrightarrow$ orbits of $\gamma_1$.
\item  \label{cuspctc}  The index contribution $\ind(\gamma_\infty)$ from a cusp with rep.~$\bg\in \ni(G,\bfC)^{\inn,\rd}$ is 
$|(\bg)\text{Cu}_4/\sQ''|-1$. 
\end{edesc}

Reminder: The {\sl index\/} of $g\in S_n$ with $t$ orbits is $\ind(g)\eqdef n-t$. App.~\ref{F2Z3} does one example computation of
\eqref{cuspct}.
\cite[\S2.8]{BFr02} computes modular curve genera from this viewpoint, while 
\cite[\S2.10]{BFr02}  and 
\cite[Cor.~8.3]{BFr02} show how the $\sh$-incidence matrix works effectively to do much harder genus computations where the group is
respectively
$A_5$ and $G_1(A_5)$. 

\subsection{More on Schur multipliers and Frattini covers of a subgroup} \label{schurDisc} 
We list results on Schur multipliers and Frattini covers used, say, in  examples like Ex.~\ref{A5limitgps} and Ex.~\ref{schurComps}. One thing
they say is that a
$\bZ/p$ quotient at the head of $M_k=\ker(G_{k+1}\to G_k)$ makes a special contribution to the $\bZ/p$ quotients at the head of all
$M_t\,$s,
$t\ge k$. So, the appearance of a Schur multiplier of a simple group at level 0 affects all levels of a \MT. 

\subsubsection{Two Schur multiplier topics}  Use notation of \S\ref{pperfect}. A 
$\bZ/p$ quotient of $\SM_G$ has  {\sl height\/} the largest $u$
with $\SM_{G,p}\to \bZ/p$ factoring through $\bZ/p^u$.  

\begin{edesc}  \item Given a $\bZ/p$ quotient of $\SM_G$, what is its  height? 
\item When do  $\bZ/p$ quotients of $\SM_{G,p}$  arise from pullback of Schur multipliers of classical groups?
\end{edesc} 

\cite{Fr02} and
\cite{FrS05} have a general classification of Schur
multipliers by how they append to
$M_k=\ker(G_{p,k+1}\to G_{p,k})$. Also, a Schur multiplier appearing at level $k$ {\sl replicates\/} to higher levels in a form called {\sl
antecedent} (\S\ref{repComps}). 

The archetype is the sequence of groups
$\{G_{2,k}(A_n)\}_{k=0}^\infty$,
$n\ge 4$. For each $k$, there is a $\bZ/2$ quotient of the Schur multiplier of $G_{2,k}(A_n)$ that is 
the antecedent of  the 2-Frattini central Spin cover $\Spin_n\to A_n$. Often antecedents inherit
properties from the original Schur multiplier. Here are two examples. 

\begin{edesc} \label{height} \item \label{heighta}  If $u$ is the height of a $\bZ/p$ quotient of $\SM_G$, then 
it is also the height of its antecedent in $\SM_{G_{p,k}}$ \cite[\S4.4]{FrS05}. 
\item \label{heightb}  For $p=2$, if a $\bZ/2$ quotient of $\SM_G$ is the pullback to $\Spin_N$ of an embedding  $G\le A_N$, 
some $N$, then an effective test decides if the antecedent of  $\SM_{G_{p,k}}$ is from an embedding 
$G_{p,k} \le A_{N'}$, some $N'$.
\end{edesc}  \cite[\S9.4]{BFr02} shows by example how  \eql{height}{heightb} contributes. It separates the two braid orbits of
$\ni(G_1(A_5), \bfC_{3^4})$ (as at the top of \S\ref{nz.1}) by the lifting invariant (\S\ref{smallInv}) from the pullback of $G_1(A_5)\le
A_{N'}$ with various values of
$N'$ (40, 60 and 120). This isn't so effective as to decide in one fell swoop the story of braid orbits for $\{\ni(G_k(A_5),
\bfC_{3^4})\}_{k=0}^\infty$. Still, that is our heading. 

Finally, Prop.~\ref{densityOnes} shows, even for $p=2$, Schur multipliers relating to spin covers of groups don't exhaust 
all  Schur multipliers that conceivably affect computations on \MT\ levels. 
\cite[\S5.7]{BFr02} explains its dependence on \cite{GS78}: That the condition that $M_0$ (and so $M_k$) being 1-dimensional is
equivalent to $G_0$ being a slight generalization (supersolvable) of dihedral groups. As a special case,  if $M_0$ is not 1-dimensional, 
then $\one_{G_k}$ (see \S\ref{gpStatements}) appears with an explicit positive  density in $M_k$ for $k$ large. Though 
effective, for small  $k$ it is subtle to predict the 
appearance of $\one_{G_k}$ and, for all $k$, where in the Loewy display the  $\one_{G_k}\,$s appear. 

Recall: Over an
algebraically closed field the set of simple $G_0$ modules has the same 
cardinality as the set of $p'$
conjugacy classes. Let $S$ be any simple $G_0$ module. Let $K$ be  
algebraically closed and
retain the notation $M_k$ after tensoring with $K$. We use $\lrang{S,M_k}$, and 
related compatible notation, 
for the total multiplicity of $S$ in all  Loewy layers of the $G_k$ module 
$M_k$. Let $O_{p'}(G)$ be
the maximal  normal $p'$subgroup of  finite group $G$ (it is the same for each $G_k$).

\begin{prop}[\cite{darren3}] \label{densityOnes}  If $\dim_K(M_0) \ne 1$,  then 
$$\lim_{n\mapsto
\infty}
\frac{\lrang{S,
M_k} }{\dim_K( M_k)} = \frac{ \lrang{S, K[G/ O_{p'}(G)]}}{\dim_K(K[G/
O_{p'}(G)])}.$$
\end{prop} 

\subsubsection{Frattini covers of a subgroup of $G$} I can't find the following useful lemma (applied in Rem.~\ref{remembedHG},
Lem.~\ref{lemSchur-Frat} and Princ.~\ref{FP3}) in my previous publications. 

\begin{lem} \label{embedHG} Let $H\le G$. Then, 
for each $k$ there is an embedding (not unique) $\beta_k: G_{p,k}(H)\to G_{p,k}(G)$ lying over the embedding of $H$ in $G$. \end{lem} 

\begin{proof} The lemma follows from Schur-Zassenhaus if $H$ is a $p'$ group where we use $G_{p,k}(H)$ to be $H$ itself. Now assume $H$ is
not $p'$. The pullback
$\text{inj}_k^{-1}(H)$ of
$H$ in
$G_{p,k}(G)$ is an extension with
$p$ group kernel having exponent
$p^k$.   From the versal property of $G_{p,k}(H)$ that produces $\beta_k: G_{p,k}(H)\to
\text{inj}_k^{-1}(H)\le G_{p,k}(G)$ lying over the embedding of
$H$ in
$G$.

Denote pullback of $H$ in $\tG p$ by $\text{inj}^{-1}(H)$. Since   ${}_p\tilde H \to H$ is the {\sl minimal\/}
cover of $H$ with kernel pro-free $p$-Sylow \cite[Prop.~20.33]{FrJ86}, there is homomorphism $\text{inj}^{-1}(H) \to {}_p\tilde H$. This 
induces  $\psi_k:
\text{inj}_k^{-1}(H)\to G_{p,k}(H)$ in the other direction. The compositions $\psi_k\circ\beta_k: G_{p,k}(H)\to G_{p,k}(H)$ are onto: 
They lie over the identity on $H$ and $G_{p,k}(H)\to H$ is a Frattini cover. So, acting on a finite group, they must be one-one. In
particular,
$\beta_k$ is one-one.    
\end{proof}

\begin{rem} \label{remembedHG} The proof that gives  $\beta_k$ in Lem.~\ref{embedHG}  extends it inductively to some $\beta_{k+1}$. So,
we may choose  $\{\beta_k\}_{k=0}^\infty$ compatibly, coming from an injection $\beta: {}_p \tilde H\to \tG p$. Also, if $G_k\to G$
factors through any $\mu: G'\to G$, then we may compose $\beta_k$ with $\mu$. When notation allows, continue to denote the resulting map
$G_{p,k}(H)\to G'$ by $\beta_k$.\end{rem}

\section{Projective systems of braid orbits} \label{spinGen} We consider two natural trees attached to the levels of a \MT. 

\subsection{Projective systems of components} \label{CTstart} Restrict the maps  $\bar \sH_{k+1}\to \bar \sH_k$ to cusps and components to 
respectively define a  {\sl cusp--tree\/} $\sC_{G,\bfC,p}$ and a {\sl component-tree\/} $\sT_{G,\bfC,p}$ directed by increasing levels. A {\sl
branch\/} on one of  these trees is a maximal (directed upward) path; so it starts at level 0. Containment of cusps in their components induces 
a map from $\sC_{G,\bfC,p}$ to  $\sT_{G,\bfC,p}$. 

\subsubsection{Cusp branches} \label{cuspBranches} The Nielsen class view of this regards the vertices of $\sC_{G,\bfC,p}$
(resp.~$\sT_{G,\bfC,p}$) as $\text{Cu}_4$ (resp.~$\bar M_4$) orbits on the collections $\{\ni(G_k,\bfC)^{\inn,\rd}\}_{k=0}^\infty$. Yet, 
we need the spaces to consider absolute Galois groups acting on these trees. 

Let $F_\bfC$ be the subfield in the
cyclotomic numbers fixed by  $\{n\in \tilde \bZ^*\mid  \bfC^n=\bfC\}$, where  equality is of sets  with multiplicities.   \cite[Prop.~1]{FrV2} 
says (in general) the spaces $\sH(G,\bfC)^\inn$ (with their maps to $U_r$ interpreted as moduli spaces) have minimal definition field $F_\bfC$.
This implies $F_\bfC$ is a definition field for $\sH(G,\bfC)^{\inn,\rd}$ (with its similar moduli properties), and so for the system of
spaces
$\{\bar \sH(G_k,\bfC)\}^{\inn,\rd}\}_{k=0}^\infty$ (with their compatible maps to $\prP^1_j$). 

\begin{lem} \label{GFaction} The absolute Galois group  $G_{F_\bfC}$ acts compatibly on the vertices of $\sC_{G,\bfC,p}$ and $\sT_{G,\bfC,p}$.
So,
$G_{F_\bfC}$ acts compatibly as permutations on (finite or infinite) branches of $\sC_{G,\bfC,p}$ and $\sT_{G,\bfC,p}$. 

Assume a cusp branch $\!B$ 
defines component branch $\!B'$. \!If, modulo braiding, \!$G_{F_\bfC}$ \!has a finite orbit on (resp.~\!fixes)  $B$, then it has a finite orbit
on (resp.~\!fixes) 
$B'$. 
\end{lem} 

\S\ref{comp-tang} notes we know many places where the \lq\lq finite orbit on $B$\rq\rq\ hypothesis of Lem.~\ref{GFaction} holds, with $B$ an  
H-M cusp branch (Ex.~\ref{HMbranch}).    The  modular curve tower $\{X_1(p^{k+1})\}_{k=0}^\infty$ has just one component-branch. We 
understand its cusp-branches well. Manin-Demjanenko (\cite[Chap.~5]{SeMW} or \cite[\S5.3]{Fr02}) gave this case of 
Conj.~\ref{wconj} long before Faltings' Theorem. (We apply Faltings to treat general \MT s.) It is typical to define a branch of
$\sT_{G,\bfC,p}$ by labeling it from the image of a branch of $\sC_{G,\bfC,p}$. See Princ.~\ref{FP2} and Ex.~\ref{HMbranch}. 

There is nothing to prove in Conj.~\ref{wconj} if the  $\sH_k$ (or  $\ni(G_k,\bfC)$) are empty for large $k$. This happens in
one of the two components of the \MT\ for $(A_n,\bfC_{3^r},p=2)$ with $r\ge n\ge 4$ (or if $r=n-1$ and $n$ is even) 
\cite[Main Result]{altGps}. For  $n=4=r$ see App.~\ref{Z3limgps}. This gives a necessary situation for a number field $K$ for 
 considering Conj.~\ref{wconj}:  There is an infinite component branch   
\begin{equation} \label{lowgenus} B'\eqdef \{\bar \sH_k'\Leftrightarrow \bar
M_4 \text{ orbit }\ni_k'\}_{k=0}^\infty \text{ fixed by $G_K$ (as in Lem.~\ref{GFaction})}.\end{equation} 

\S\ref{cuspTypes} divides cusps into three types. It is easier to describe the cusps than to place them in components. 
\S\ref{cuspGrowth}  describes how projective systems of $p$ cusps contribute to indices in the Riemann-Hurwitz formula. 

\newcommand{\Ge}{\text{\rm Ge}}

\subsubsection{Sequences of component genera} \label{compGenus}
Restrict the  
$\gamma\,$s of \eqref{gammas} to $ \ni_k'$ in \eqref{lowgenus}. This gives $(\gamma_{0,k}',\gamma_{1,k}',\gamma_{\infty,k}')$ defining the genus
$g_{\bar
\sH'_k}$ of 
$\bar
\sH'_k$:     
\begin{equation} \label{genusseq} 2(\deg(\bar
\sH_k'/\prP^1_j)+g_{\bar \sH'_k}-1)= 
\ind(\gamma_{0,k}')+\ind(\gamma_{1,k}')+\ind(\gamma_{\infty,k}').\end{equation} 

Below we denote the genera sequence for the branch $B'$ by $\Ge_{B'}\eqdef\{g_{\bar
\sH_k'}\}_{k=0}^\infty$.  The strongest results toward the Main
Conjectures require two contributions: 
\begin{edesc} \label{infoBranches} \item \label{infoBranchesa} Deciphering the infinite branches from the finite branches. 
\item \label{infoBranchesb} Separating cusp branches into types that indicate their contributions to Riemann-Hurwitz. 
\end{edesc} \cite[Lect.~1]{Ont05}  
starts by computing  modular curve genera from a
\MT\ viewpoint. \S\ref{cuspTypes} describes those cusp types, including the significant special cusps called g-$p'$, and the  corresponding
g-$p'$ cusp branches. The following is a prototype modular curve property, and \cite{FrS05} uses it as an explicit target. 

\begin{quest} \label{wildGuess} Suppose $K$ is a number field and $B'$ is an (infinite) K component branch with $B'$ the image of a g-$p'$
cusp branch $B\in \sC_{G,\bfC,p}$. Is it possible to give a closed  expression for the elements of $\Ge_{B'}$?
\end{quest} 

\subsubsection{Reduction to  $G_0=G$ has no $p$-part to its center} \label{redAbel} One part of Princ.~\ref{FP1} says that $p$ cusps contribute
highly to cusp ramification. That result is a subtle use of Prop.~\ref{nopcenter}. This reduces considering \MT s (or at least the Main
Conjecture) to the case where for all $k$, the $p$-part of the center is trivial. Denote the center of a group
$G$ by $Z(G)$, and the $p$-part of the center by $Z_p(G)$.

\begin{prop}[$p$-Center Reduction] \label{nopcenter} Suppose $G=G_0$ is a $p$-perfect group  with $Z_p(G)\ne \{1\}$. Then,
there is a $p$-Frattini cover  $\psi^c: G \to G^c$ with $Z_p(G^c)$ trivial (and $G^c$ is $p$-perfect). Any $p'$ conjugacy class $\C$ of $G$ has
a unique image class in $G^c$  which we also donate by $\C$ (\S\ref{gpStatements}). In particular, Main Conj.~\ref{wconj} holds for
$(G^c,\bfC^c,p)$, if and only if it holds for $(G,\bfC,p)$.
\end{prop} 

\begin{proof} Let $U_p$ be the maximal normal $p$-Sylow of $G$, and let $\Phi(U_p)$ be the Frattini subgroup of $U_p$.  Then, $G \to
G/\Phi(U_p)$ is a $p$-Frattini cover.

First consider the case $G$ is $p$-split: $G=U_p\xs G/U_p$. From $G$ being $p$-perfect,  $G/U_p$ has no fixed points on
$U_p/\Phi(U_p)$. So $Z_p(G/\Phi(U_p))=\{1\}$. General case: Form $G/\Phi(U_p)$. We're done if $Z_p(G/\Phi(U_p))$ is trivial. Otherwise
iterate this to achieve $G^c$. 

Now consider the last sentence of the proposition. Since $G\to G^c$ is a $p$-Frattini cover, the universal $p$-Frattini cover of $G^c$ is the
same as that of $G$. Denote the
$k$th characteristic Frattini extension of $G^c$ by $G_k^c$. From the construction, there is a $k_0$ so that $G^c_{k_0}\to G^c$ factors through
$G\to G^c$. Conclude easily for each $k$ there is a corresponding $k'$ so that $G^c_{k'}\to G^c$ factors through $G_k\to G$. Also, 
the map $\psi_k: G_k\to G$ composed with $\psi^c$ factors through $G_k^c\to G$. 

In particular, this means for $k>>0$ there is a $k'$ so that $\sH(G_{k'}^c,\bfC)^{\inn,\rd}$ naturally
maps (surjectively, over any field containing their simultaneous definition fields) to 
$\sH(G_{k},\bfC)^{\inn,\rd}$. So: if $\sH(G_{k},\bfC)^{\inn,\rd}(K)=\emptyset$, then $\sH(G_{k'}^c,\bfC)^{\inn,\rd}(K)=\emptyset$; if
$\sH(G_{k'}^c,\bfC)^{\inn,\rd}(K)\ne \emptyset$, then $\sH(G_{k},\bfC)^{\inn,\rd}(K)\ne \emptyset$; etc. Conclude $(G^c,\bfC,p)$ and
$(G,\bfC,p)$ simultaneously pass or fail the conclusion of the Main Conjecture. 
\end{proof}

\begin{rem}[Center considerations] \label{consCenter} Do not conclude from Prop.~\ref{nopcenter} that 
\MT s can't handle groups with centers. All our sections and also of \cite{BFr02}  must consider that $\tG p$ is full of subquotient
sequences of the form $\psi': R'\to G'$, a central extension of $G'$, with $\ker(\psi')$ a quotient of $G'\,$s Schur multiplier. As in
\S\ref{weigCor}, it is the maximal elementary $p$-quotient of 
$G_k\,$'s Schur multiplier that controls major properties of \MT\ levels.  

Use the notation of Prop.~\ref{nopcenter}.   Denote the $p'$ part of
$Z(G)$ by $Z_{p'}(G)$. Then, for all $k$, 
$Z_{p'}(G_k)=Z_{p'}(G_k^c)=Z_{p'}(G)$. See this by identifying $\tG p$ with the universal $p$-Frattini of $G/Z_{p'}$ fiber product with $G$
over $G/Z_{p'}$. We could have continued the map
$\psi^c: G\to G^c$ through
$G_c\to G_c/Z_{p'}(G)$. That would, however, complicate the final conclusion of Prop.~\ref{nopcenter}. No longer
could we canonically  identify the image conjugacy classes with $\bfC$. So, while \MT s already deals seriously with the $p$-part of centers, 
K.~Kimura's master's thesis \cite{kimura} has a point in considering phenomena that arise from the $p'$-part. \end{rem}

\subsection{g-$p'$ and o-$p'$ cusps, and Frattini Principles 1 and 2} \S\ref{cuspTypes} defines the three cusp types 
using a representative $\bg=(\row g 4)$ of the cusp orbit. We expect g-$p'$ cusp branches to give outcomes like that of Quest.~\ref{wildGuess}. 
Modulo Conj.~\ref{noWeigB}, we expect
some g-$p'$ cusp branch defines any component branch {\sl with all levels having a fixed number field as definition field\/}. 
\S\ref{cuspExamps} considers cases when we can use g-$p'$ cusps to get a handle on o-$p'$ cusps. 

\subsubsection{The cusp types} \label{cuspTypes} Use these notations: $$H_{2,3}(\bg)\eqdef \lrang{g_2,g_3}\text{
and }H_{1,4}(\bg)=\lrang{g_1,g_4};$$ and $(\bg)\mpr \eqdef
\ord(g_2g_3)$, the order of the {\sl middle product}. Primary contributions after level 0 to \eqref{genusseq} come from
{\sl $p$ cusps}: $p| (\bg)\mpr$. Here are the other types. 

\begin{edesc} \label{gp'op'} \item \label{gp'op'a} {\sl g(roup)-$p'$}: 
$H_{2,3}(\bg)$ and $H_{1,4}(\bg)$ are
$p'$ groups. 
\item \label{gp'op'b} {\sl o(nly)-$p'$}: $p\not | (\bg)\mpr$, but the cusp is not g-$p'$. 
\end{edesc} 

 Let $\{{}_k\bg=({{}_kg}_1, {{}_kg}_2,{{}_kg}_3,{{}_kg}_4)\in
\ni_k'\}_{k=0}^\infty$ be a projective system of cusp representatives. Then ${}_k\bg$ corresponds to a braid orbit $\ni_k'\subset
\ni(G_k,\bfC)$, and therefore to a component $ \sH_k'\subset \sH(G_k,\bfC)^{\inn,\rd}$.  Denote the corresponding projective system of cusps by
$\{\bp_k\in \sH_k'\}_{k=0}^\infty$.
When  a point $\bp'$ on  some space lies over another point $\bp$, denote the ramification order ({\sl index, or width\/}) of $\bp'/\bp$ by
$e(\bp'/\bp)$. Crucial to our Main Conjecture is the phenomenon that $p$ cusp widths grow automatically as we go up \MT\ levels. The formal
statement, coming mostly from 
\cite[\S8.1]{BFr02}, is our first use of the Frattini property. Recall: $Z_p(G)$ is the $p$-part of the center of $G$
(\S\ref{redAbel}). 

\begin{princ}[Frattini Princ.~1]  \label{FP1} If $p^u|
({}_k\bg)\mpr$, $u\ge 1$, then
 $p^{u+1}|({}_{k+1}\bg)\mpr$. 

Assume $Z_p(G)$ is trivial. Then,  for  $p$ odd (resp.~$p=2$) and  $k\ge 0$ (resp.~$k>>0$) $e(\bp_{k+1}/\bp_k)$ is $p$.
\end{princ} 

\begin{proof}[Comments on explicitness]  The first part is a consequence of \cite[Lift Lem.~4.1]{FrK97} (for example). It comes from this 
simple statement: All lifts to $G_{k+1}$ of an element of order $p$ in $G_k$ have order $p^2$.  That concludes the first part. 

Denote the operator that takes
any 
$(a,b)\in G^2$ to
$(aba^{-1},a)$ by $\gamma$. Then,
\cite[Prop.~2.17]{BFr02} \wsp \S\ref{BFrtypos} has a typo free statement with $(g_1,g_2)$ replacing $(a,b)$ \wsp tells  how to
compute the length of the orbit (using no equivalence between pairs) of
$\gamma$ generated by
$(a,b)$. The length of the $\gamma^2$ orbit is $$o(a,b)\eqdef o=\ord(a\cdot b)/|\lrang{a\cdot b}\cap Z(a,b)|.$$  Then, one
of the following holds for the length $o'(a,b)=o'$ of the
$\gamma$ orbit on $(a,b)$. Either: $a=b$ and $o'=1$, or;  
\begin{triv} \label{abram} if $o$ is odd and $b(a\cdot b)^{\frac{o-1}2}$ has order 2, then $o'=o$;  or else $o'=2\cdot o$. \end{triv}

\cite[Lem.~8.2]{BFr02} of necessity was intricate, for it's goal  was to nail $e(\bp_{k+1}/\bp_k)$ from data
on the group theoretic cusp from ${}_k\bg$ and ${}_{k+1}\bg$. This was to  precisely list genera of examples. We now say
this result in a more relaxed way. 

Assume $Z_p(G_0)$ is trivial. From \cite[Prop.~3.21]{BFr02} the same therefore holds for $Z_p(G_k)$ for all
$k\ge 0$. All we care about in our conclusion is the $p$ part of $e(\bp_{k+1}/\bp_k)$.   We divide the contribution to the $p$ part
$e_p(\bp_k/\infty)$ into two cases:
$p$ odd, and
$p=2$. 

When $p$ is odd, \cite[Lem.~8.2]{BFr02} gives  $e_p(\bp_k/\infty)$ as the $p$-part of $o({}_kg_2,{}_kg_3)$. If the $p$-part
of $|\lrang{{}_kg_2\cdot {}_kg_3}\cap Z({}_kg_2,{}_kg_3)|\eqdef k_p(2,3)$ is trivial, then the result is the $p$-part of $\ord({}_kg_2\cdot
{}_kg_3)$. Since the $p'$-part of $\ord({}_kg_2\cdot {}_kg_3)$ is unchanging with $k$, the
first statement in the proposition gives $e(\bp_{k+1}/\bp_k)=p$. 

To see why $k_p(2,3)=1$, use that the action of $\sQ''$ expresses the cusp width also
from $({}_k g_4,{}_k g_1)$ (\S\ref{NielClDict}). The result must be the same, using an analogous expression $k_p(1,4)$. Since $({}_k
g_4\cdot {}_k g_1)^{-1}={}_kg_2{}_kg_3$, then 
$k_p(2,3)=k_p(1,4)$. Now if both are nontrivial, it means $$Z_p({}_kg_2,{}_kg_3)\ge ({}_kg_2\cdot {}_kg_3)^{\ord({}_kg_2\cdot {}_kg_3)/p}\le  
Z_p({}_kg_4,{}_kg_1).$$ Since ${}_kg_1, {}_kg_2,{}_kg_3, {}_kg_4$ generate $G_k$, this implies $Z_p(G_k)$ is
nontrivial. 

For $p=2$, the computation works similarly, except for factors of 2-power order (bounded by 4) in $e_2(\bp_k/\infty)$  from the action of
$\sQ''$ and the distinction between
$o=o'$ and
$o=2\cdot o'$ given in \eqref{abram}. These are, however, regular behaviors.  Observations like those
about $\sQ''$ in \S\ref{fineredMod}, allow replacing $k>>0$ by a more precise statement.    
\end{proof} 

\begin{princ}[Frattini Princ.~2] \label{FP2} The definition of     
$p'$ and g-$p'$ cusp doesn't depend on its 
rep.~in  $(\bg)\text{\rm Cu}_4$ \cite[Prop.~5.1]{FrS05}. If ${}_0\bg\in \ni(G_0,\bfC)$ represents a g-$p'$ cusp,  then above it
there is a g-$p'$ cusp branch $\{{}_k\bg\in \ni(G_k,\bfC)\}$. 
\end{princ} 

\begin{proof} Use $(g_1,g_2,g_3,g_4)$ for ${}_0\bg$. Let $H\le {}_0 G$ be a $p'$ group. Then, consider the pullback $\psi^{-1}(H)$ in $\tG
p$. The  profinite version of Schur-Zassenhaus says  the extension
$\psi^{-1}(H)\to H$ splits \cite[20.45]{FrJ86}. Apply this to each $p'$ group $H_{1,4}({}_0\bg)$ and
$H_{2,3}({}_0\bg)$. This gives  $H_{1,4}', H_{2,3}'\le \tG p$, defined up to conjugation by $\tilde P_p$, mapping one-one to their
counterparts modulo reduction by
${}_p\tilde P$.  

Let $g_1',g_4'\in H_{1,4}'$ (resp.~$g_2',g_3'\in H_{2,3}'$) be the elements over
$g_1,g_4\in H_{1,4}$ (resp.~$g_2,g_3\in H_{2,3}$). 
Then, $g_2'g_3'$ is conjugate to $(g_1'g_4')^{-1}$ by some $h\in {}_p\tilde P$.
Replace
$H_{1,4}'$ by its conjugate by $h$ to find $\bg'=(\row {g'} 4)\in \ni(\tG p,\bfC)$ lying over ${}_0\bg$. The images of $\bg'$ in each 
$\ni(G_k,\bfC)$ give the desired g-$p'$ cusp branch.   
\end{proof}

\begin{exmp}[$\sh$  of an H-M rep] \label{HMbranch} \S\ref{NielClDict} has the definition of the shift $\sh$. A H(arbater)-M(umford) rep.~in the
reduced Nielsen class
$\ni(G,\bfC)^\rd$ (applies to inner or absolute equivalence) has the shape $\bg=(g_1,g_1^{-1},g_2,g_2^{-1})$. Then, 
$(\bg)\sh$ is clearly a g-$p'$ cusp. It has width 1 or 2. A formula distinguishes between the cases (proof of Prop.~\ref{FP1}).
Typically our examples have 
$H_{2,3}(\bg)\cap H_{1,4}(\bg)=\lrang{1}$, or else $G=\lrang{g_1,g_2}$ has a nontrivial cyclic $p'$ kernel dividing the orders of
$\lrang{g_i}$, $i=1,2$. 
\end{exmp} 

\subsubsection{Consequences of fine reduced moduli} \label{fineredMod} The reduced spaces of the levels of a component branch are moduli
spaces. Using them as moduli spaces behooves us to know when they have (reduced) fine moduli: objects that represent points do so in a
unique way. There isn't a prayer they have fine moduli unless the corresponding unreduced spaces $\sH(G_k,\bfC)^{\inn}$ have fine
moduli. For that, the if and only if criterion, given that $G_0$ is $p$-perfect, is that $G_0$ has no center \cite[Prop.~3.21]{BFr02}. 

Given this, 
\cite[Prop.~4.7]{BFr02} gives the if and only if criterion for level $k_0$ of a branch to have (reduced) fine moduli. This says: 
Two computational conclusions hold from the action of $H_4$ and $\bar M_4$ on the corresponding level $k_0$ braid orbit $\ni_{k_0}'$:  
\begin{edesc} \label{fmst} \item \label{fmsta} $\sQ''$ has all its orbits on $\ni_{k_0}'$ of length 4; and 
\item \label{fmstb} both $\gamma_{0,k_0}'$ and $\gamma_{1,k_0}'$ act without fixed point. 
\end{edesc} 

Both Thm.~\ref{wmcneedspcusps} and  \S\ref{BrFrattProp}, on the {\sl Branch Frattini Propery},  use
Lem.~\ref{fineredBr}. 

\begin{lem} \label{fineredBr} For any $k$, $\bar \sH_{k+1}'/\bar \sH_k'$ 
ramifies only over cusps (points  over $j=\infty$) if and only if \eql{fmst}{fmstb} holds. If \eql{fmst}{fmstb} holds for $k=k_0$, then it
holds also for  $k\ge k_0$, and for each such $k$, $p$ is the ramification index for each prime ramified in the cover $\bar\sH_{k+1}'/\bar
\sH_k'$.  So, this holds if the component branch $B'$ has
fine moduli (for $k= k_0$).  
\end{lem}

\begin{proof} The cover $\bar \sH_k'\to \prP^1_j$ ramifies only over $j=0,1,\infty$. The lengths of the disjoint cycles for  $\gamma_{0,k}'$
(resp.~$\gamma_{1,k}'$)  on $\ni_k'$  correspond to the orders of ramification of the points of  $\bar \sH_k'$ lying over 0 (resp.~1). 

\!Apply  multiplicativeness of ramification  to  $\bar\sH_{k+1}'\longmapright{\psi_{k+1,k}\ }{20} \bar
\sH_k'\mapright{\psi_k}
\prP^1_j$. If $\bp_{k+1}\in \bar\sH_{k+1}'$, denote $\psi_{k+1,k}(\bp_{k+1})$ by $\bp_k$. Then, $\bp_{k+1}/\psi_k\circ
\psi_{k+1,k}(\bp_{k+1})$ has ramification index  
\begin{equation} \label{abhyLem} e(\bp_{k+1}/\psi_k\circ \psi_{k+1,k}(\bp_{k+1}))=
e(\bp_{k+1}/\bp_k)e(\bp_{k}/
\psi_k\circ \psi_{k+1,k}(\bp_{k+1})).\end{equation} 
If $\psi_k\circ \psi_{k+1,k}(\bp_{k+1})=0$, then $e(\bp_{k+1}/0)=1$ and $e(\bp_{k}/0)=1$ are each either 1 or 3 (\S\ref{redNCcusps}).
Conclude from \eqref{abhyLem}:  $$e(\bp_{k+1}/0)=3\text{ and }e(\bp_{k+1}/ \psi_k\circ \psi_{k+1,k}(\bp_{k+1}))=1$$ both hold if and only if
$e(\bp_{k}/0)=3$. 

Statement \eql{fmst}{fmstb} for $\gamma_{0,k_0}'$ says $e(\psi_{k_0}/0)$ is 3 for each $\psi_{k_0}$ lying over 0. This inductively implies
no point of $\bar\sH_{k+1}$ lying over $0\in \prP^1_j$ ramifies over $\bar \sH_k$ if $k\ge k_0$. The same argument works for $\gamma_{1,k_0}$
and concludes the lemma.  
\end{proof}  

\begin{exmp}[When reduced fine moduli holds] \label{nononcuspBr} 
For all the examples of \cite[Chap.~9]{BFr02}, reduced fine moduli holds with $k_0=1$ in Lem.~\ref{fineredBr}. \cite{FrWProof} shows for $p=2$
any H-M component branch has fine moduli. We hope to expand that considerably before publishing a final version. If Conj.~\ref{gp'givesPSC} is
true, then that implies any (infinite) component branch of any of the many $A_4$ and $A_5$ ($p=2$ and any type of $2'$ conjugacy classes) \MT s
have reduced fine moduli.
\end{exmp}

\begin{exmp}[\eql{fmst}{fmstb} can hold without fine moduli]  
Here again, we have a modular curve comparison with a highlight from \cite[\S 4.3.2]{BFr02}. While there is
a one-one map (onto) map $\sH(D_{p^{k+1}},\bfC_{2^4})^{\inn}\to Y_1(p^{k+1})$ (\S\ref{twoExs}), the spaces, as moduli spaces, are not exactly
the same. The latter has fine moduli, but the former does not. The distinction is that the moduli problem for
$\sH(D_{p^{k+1}},\bfC_{2^4})^{\inn}$ is {\sl finer\/} than that for
$Y_1(p^{k+1})$: There are \lq\lq more\rq\rq\ genus 1 Galois covers of $\prP^1_z$ with $D_{p^{k+1}}$ monodromy  than there are corresponding
elliptic curve isogenies. Still, \eql{fmst}{fmstb} holds.  
\end{exmp}

\subsubsection{Relations between g-$p'$ and o-$p'$ cusps} \label{cuspExamps}  For our arithmetic conjectures we only care about infinite $K$
component branches (\S\ref{cuspResults}) where $K$ is some number field. For this discussion we accept Conj.~\ref{gp'givesPSC}. That means in
dealing with the possibility of o-$p'$ cusp branches, we only need to consider those that appear on a g-$p'$ component branch.   Since o-$p'$
cusp branches are so important, we hope thereby to be as explicit with them as with g-$p'$ cusps.  

This occurs, for example, if an o-$p'$ cusp
is over a g-$p'$ cusp. To simplify, start with a  g-$p'$ cusp
${}_0\bg$ at level 0 with
$({}_0\bg)\mpr\eqdef v$ of order $c$. Prop.~\ref{op'overgp'} shows the conditions of
\eqref{norm} sometimes hold (though not for shifts of H-M reps., Ex.~\ref{HMbranch}). 

Expressions in \eqref{norm}  are in
additive notation in  $M_0=\ker(G_1\to G_0)$; the group ring $\bZ/p[G_0]$ acts on the right. For $g\in G_0$ and $m\in  M_0$, denote the
subspace of
$M_0$ that commutes with $g$ (on which $g$ acts trivially) by $\text{\rm Cen}_g$, and its 
translate by $m$ by $\text{\rm Cen}_g-m$. Denote 
$1+v+\cdots + v^{c-1}: M_0\to M_0$ by 
$L(v)$. 

\begin{prop} \label{gp'-op'imp} Suppose  
$\bg'\in\ni(G_{1},\bfC)$ lying over ${}_0\bg$ is neither a g-$p'$, nor a $p$ (so is an o-$p'$),  cusp. Let  $\bg\in
\ni(G_1,\bfC)$ be a g-$p'$ cusp over ${}_0\bg$ as in the conclusion of Princ.~\ref{FP2}.  Then, with no loss we may
assume $$\bg'=((m^*)^{-1}g_1m^*, g_2,m_3 g_3m_3^{-1}, (m_4m^*)^{-1}g_4 (m_4m^*))$$ with
$m^*,m_3,m_4\in M_k$ and $(g_2,g_3)$ is not conjugate to $(g_2,g_3')$. 

Then, the order of $(\bg')\mpr$  is  $c$ and the following are equivalent to $\bg'$ being o-$p'$.
\begin{edesc} \label{norm} \item \label{norma} Product-one: $m_3({{}_0g_3}-1) + m_4({{}_0g_4}-1) + m^*(v-1)=0$. 
\item \label{normb} $p'$ middle-product:  $m_3({}_0g_3 - 1)$ is an element of $M_0(v-1)$.  
\item \label{normc} Not g-$p'$: It does not hold that $m_3({}_0g_2-1)\in \text{\rm Cen}_{g_3}({}_0g_2-1)$.    
\end{edesc} 
\end{prop} 

\begin{proof} Since $\bg'$ is an o-$p'$ cusp, we may assume $H_{2,3}(\bg')$ is not a $p'$ group. 
Characterize this by saying $H_{2,3}(\bg)$ is not conjugate to $H_{2,3}(\bg')$. By conjugating, we may
assume $g_2=g_2'$ and $g_3'=m_3g_3m_3^{-1}$ for some $m_3\in M_k\setminus\{0\}$. For $(g_2',g_3')$ to be conjugate 
to $(g_2,g_3)$ is equivalent to some $m\in M_0\setminus \{0\}$ commutes with $g_2$ while $m-m_3$  commutes with $g_3$.   The
other normalization conditions are similar. Then, 
\eql{norm}{norma} is  $\Pi(\bg')=1$ in additive notation. 

Compute $(g_2m_3 g_3m_3^{-1})^c=(g_2g_3m_3^{g_2}m_3^{-1})^c$ to get $(g_2g_3)^c=1$ times an element $u\in M_0$. That
$u$, in additive notation, is just $$(m_3)({}_0g_3-1)(1+v+v^2+\cdots v^{c-1})=(m_3)({}_0g_3-1)L(v).$$ Since $g_2g_3$ is assumed $p'$, that gives
$u=0$, or $(m_3)({}_0g_3-1)$ is in the kernel of $L(v)$.  As, however, $v$ has $p'$ order, the characteristic polynomial $x^c-1$ of $v$ has
no repeated roots. So, $M_0$ decomposes as a direct sum $\bZ/p[x]/(x-1)\oplus \bZ/p[x]/L(x)$ with $v$ acting in each factor as multiplication
by $x$.  Thus, the kernel of 
$L(v): M_0\to M_0$ is exactly the image of $(1-v)$.  That is, $v$ having $p'$ order is equivalent to $m_3({}_0g_3 - 1)$ is an element of
$M_0(v-1)$. That completes showing
\eql{norm}{normb}. 

Finally, suppose there is $m\in M_0$ that conjugates $(g_2',g_3')$ to $(g_2,g_3)$. Compute to see this is equivalent to $\text{\rm
Cen}_{g_2}-m_3\cap \text{\rm Cen}_{g_3}=\emptyset$. An $m$ in this overlap would satisfy  $m_3({}_0g_2-1)=
m({}_0g_2-1)$. Statement \eql{norm}{normc} says there is no such $m$. 
\end{proof} 

Apply \eqref{norm} to the shift of an H-M rep (Ex.~\ref{HMbranch}). Then, $c=1$ and $m_3$ commutes with $g_3$, contrary to assumption. So
the cusp defined by the shift of an H-M rep.~cannot have an o-$p'$ cusp over it.   Still, Prop.~\ref{op'overgp'} shows some g-$p'$ cusp
branches produce a profusion of  o-$p'$ cusps over g-$p'$ cusps. 

\begin{prop} \label{op'overgp'}  Let $\{{}_k\bg\in \ni(G_k,\bfC)\}_{k=0}^\infty$
represent a g-$p'$ cusp branch from Princ.~\ref{FP2}. Let $c_i$ be the order of ${}_0g_i$, $i=1,2,3,4$. Assume
\begin{triv} \label{op'Hyp} $O_p'(G_0)$ is trivial and
$1/c_2+1/c_3 + 1/c_4 + 1/c <1$. \end{triv} Then,  for $k$ large, an o-$p'$ cusp ${}_{k+1}\bg'\in \ni(G_{k+1},\bfC)$  lies over $\bg_k$. 
\end{prop} 

\begin{proof} Use notation of Prop.~\ref{gp'-op'imp}, starting with a g-$p'$ cusp ${}_0\bg$ at level 0. Take 
$m^*=1$. Consider what  \eqref{norm} forces on   
$\bg'=(g_1,g_2, m_3 g_3m_3^{-1}, m_4^{-1} g_4m_4)$ to force it to be an o-$p'$ cusp  in $\ni(G_1,\bfC)$. 
Condition \eql{norm}{normb} says: $$(\bg')\mpr\text{ is  }p' \Leftrightarrow
m_3({}_0g_3 - 1)\in M_0(v-1).$$ Also we must assure $m_3({}_0g_2-1)$ is not in $\text{\rm Cen}_{g_3}({}_0g_2-1)$. 

Combine all conditions of \eqref{norm}. Then,  there is an o-$p'$ cusp if and only if 
\begin{equation} \label{fundInter}  M_0({}_0g_3 - 1) \cap M_0({}_0g_4 - 1) \cap M_0(v-1) \setminus
 \text{\rm Cen}_{g_3}({}_0g_2-1)\neq \emptyset.
\end{equation} 

By the {\sl relative\/} codimension or dimension of a subspace of $M_k$, we 
mean the codimension or dimension divided by the dimension of $M_k$.  
While we can't expect \eqref{fundInter} to hold at level 0, we show it holds with conditions \eqref{op'Hyp}  if we substitute ${}_k\bg$
for ${}_0\bg$ (and $M_k$ for $M_0$) for $k$ large. 

If the relative codimension
of  $$M_0({}_0g_3 - 1) \cap M_0({}_0g_4 - 1) \cap M_0(v-1)$$ plus the relative dimension of $\text{Cen}_{g_3}$ is asymptotically less than $1$,
then  \eqref{fundInter}  holds for  $k>>0$.  Prop.~\ref{densityOnes} (using $O_{p'}(G)=\{1\}$) gives this for $k>>0$ if 
\eqref{op'Hyp} holds. So, these conditions imply an o-$p'$ cusp over ${}_k\bg$ for $k$ large.  \end{proof}

\begin{exmp}[Case satisfying \eqref{op'Hyp}] \label{manyop'} Let $G_0$ be the alternating group $A_7$ and let $p=7$. Define the Nielsen
class selecting ${}_0\bg$ with $g_2,g_3\in A_5$ both 5-cycles generating $A_5$ and having $v=g_2g_3$ a 3-cycle. From
\cite[Princ.~5.13]{BFr02} there is just one choice (up to conjugation) if  $g_1$ and $g_2$ are in the two different
conjugacy classes of order 5: $g_2  =  (5\, 4\, 3\, 2\, 1)$ and $g_3= (2\, 4\, 3\, 5\, 1)$, and $g_2g_3=(5\,3\,4)$. Now choose $g_1$ and $g_4$
analogously as 5-cycles acting on
$\{3,4,5,6,7\}$ so 
$g_4g_1$ is $(4\,3\,5)$.  Here, $H_{2,3}({}_0\bg)$ and $H_{1,4}({}_0\bg)$ are both copies of  $A_5$.  
All the $c_i\,$s are 5,  while $c = 3$.  The inequality \eqref{op'Hyp} holds: 
$
1/5 + 1/5 + 1/3 + 1/5 = 14/15 < 1$. \end{exmp}

\section{Finer graphs and infinite branches in $\sC_{G,\bfC,p}$ and $\sT_{G,\bfC,p}$} \label{finerGraphs}   We don't know what contribution 
o-$p'$  cusps in
\S\ref{cuspTypes} make to the genera of components at level $k$ on a \MT. Are they like g-$p'$ cusps in defining
projective systems of components through o-$p'$ cusps. Or, if you go to a suitably high level are all the cusps above them $p$ cusps? 
Conj.~\ref{noWeigB} says the latter holds.  
\S\ref{SectnoWeigB} consists of support for and implications of this.

Schur multipliers of quotients of the universal $p$-Frattini cover $\tG p$ of $G$ are at the center of these conclusions in the form of
lifting invariants (\S\ref{smallInv}). We must deal with these many Schur multipliers when considering graphs finer than
$\sC_{G,\bfC,p}$ and $\sT_{G,\bfC,p}$.  

\subsection{Limit Nielsen classes}  \label{limNielsen} For a full analysis of higher rank \MT\ examples such as in (\S\ref{twoExs}),
\S\ref{extGraphs} extends the previous component and cusp branch notions. This extension uses all quotients of the universal $p$-Frattini
cover (not just characteristic quotients). Given the definition of cusps from \cite[Lect.~4]{Ont05} for arbitrary values of $r$, the
concepts of this section work there, too. 
 
\subsubsection{Extending graphs to include any quotients of ${}_p\tilde G$} \label{extGraphs} 
 Let
$\sG_{G,p}$ be all finite covers
$G'\to G$ through which $\tG p\to G$ factors. Given $(G,\bfC,p)$, consider components and cusps of 
$\{\ni(G',\bfC)^{\inn}\}_{G'\in \sG_{G,p}}$. As in previous cases, they form directed graphs $\sC^f_{G,\bfC,p}$ and $\sT^f_{G,\bfC,p}$
(the $f$ superscript for {\sl full}) with maps between them. 

Now, however, there may be many kinds of maximal directed paths ({\sl branches\/}) not just distinguishing finite from infinite). Also, among
undirected paths there could be loops  because there may be several chief series for the Krull-Schmidt decomposition of
$\ker(G_{p,k+1}\to G_{p,k})$ into irreducible $G_{p,k}$ modules. This doesn't happen for $G_{2,1}(A_n)\to A_n$ for $n=4,5$, but does
for
$G_{2,2}(A_4)\to G_{2,1}(A_4)$ \cite[Cor.~5.7]{BFr02}. 

A directed path on $\sC^f_{G,\bfC,p}$ is defined by
$\{(\bg_{H_i})\Cu_4\}_{i\in I}$ with $I$ a directed set, $H_i$ a quotient of $\tG p$ and $\bg_{H_i}\in \ni(H_i,\bfC)$. If $i'>i$, then
$\tG p\to H_i$ factors through
$H_{i'}$ sending $(\bg_{H_{i'}})\Cu_4$ to  $(\bg_{H_{i}})\Cu_4$. This path defines a unique braid orbit in
$\ni(H_i,\bfC)$ for  $i\in I$: A cusp path (resp.~branch) defines a component path (resp.~branch). 

\begin{lem} \label{directedPath} A  directed path on $\sT^f_{G,\bfC,p}$ defines a set of directed paths on $\sC^f_{G,\bfC,p}$: 
Each node from any of the latter sits on a corresponding node of the former (with the obvious converse). If  $\{(\bg_{H_i})\Cu_4\}_{i\in
I}$ is a directed path, then we can choose its cusp representatives  $\bg_{H_i}$ to also be a projective system. \end{lem} 

\begin{proof} A  directed path on $\sT^f_{G,\bfC,p}$ is defined by a directed system $\{H_i\}_{i\in I}$. For each $i$ there is a node
consisting of $\sH_{H_i}$, a component of $\sH(H_i,\bfC)^{\inn,\rd}$. The points $R_{H_i}$ of the nonsingular $\sH_{H_i}$ over $j=\infty$ have
ramification degrees adding to the degree  map  $\bar \sH_{H_i}\to \prP^1_j$. For $i'\ge i$, the natural map $R_{H_{i'}}\to R_{H_i}$ 
defines a projective system of finite nonempty sets. So, the set of limits is nonempty, and each defines a directed cusp path. Let
$\{(\bg_{H_i})\Cu_4\}_{i\in I}$ be one of these (as in the correspondence of \S\ref{redNCcusps}). 

The collections $(\bg_{H_i})\Cu_4$, $i\in I$ also form a projective system of finite
nonempty sets in the set of subsets of Nielsen classes. So, they too have projective limits. Each is a projective system of the form
$\{\bg_{H_i}\}_{i\in I}$. That gives the final statement. 
\end{proof} 

\begin{defn}[$F$ paths, branches, \dots] \label{Fpaths} For $F$ a field the notion of $F$  cusp path (resp.~cusp branch), component path
(resp.~component branch) on $\sC^f_{G,\bfC,p}$ or $\sT^f_{G,\bfC,p}$ extends naturally that for $\sC_{G,\bfC,p}$ or $\sT_{G,\bfC,p}$ (as in
\S\ref{cuspResults}).
\end{defn} 

\subsubsection{Limit Groups} \label{limGroups}   
Our next definitions use notation from Lem.~\ref{directedPath}. 
\begin{defn} \label{limitGps} A directed path from a projective system $\{\bg_{H_i}\}_{i\in
I}$  has an attached group $\lim_{\infty \leftarrow i\in I} H_i=G^*$. Call this a
{\sl limit group\/} (of $(G,\bfC,p)$)  if
the directed path is maximal. Then,
$\ni(G^*,\bfC)$ is the {\sl limit Nielsen class\/} attached to the maximal path, and $\lim_{\infty\leftarrow i}\bg_{H_i}\in \ni(G^*,\bfC)$
represents the {\sl limit braid orbit\/} of the path. 
\end{defn}
We might also call $G^*$ the limit group of the braid orbit of $\bg_{G}$, or of the component of $\sH(G,\bfC)$ attached to this orbit, 
etc. 

\begin{defn} \label{obstMT} Suppose $\{\bg_{H_i}\}_{i\in I}$ defines a maximal path. Then, for each $k\ge 0$ we can ask if   
$H_i=G_k$, for some $i$. If so, we say the path goes through level $k$ of the \MT (and through braid orbit $O_{\bg_{H_i}}$). If $k_0$ is
the biggest integer with
$\{\bg_{H_i}\}_{i\in I}$ going through level $k$, then call the \MT\  {\sl obstructed\/} along the path at level $k_0$.  
\end{defn} Obvious variants on Def.~\ref{obstMT} refer to a braid orbit $O_\bg$ at level $k$ being obstructed: Every path through
$O_\bg$ is obstructed at level $k$, etc. 

If $O_*$ is the limit braid orbit in $\ni(G^*,\bfC)$ defined by a maximal path, then we say the path is obstructed at $O_*$. We also use
variations on this.  Any quotient $G'$ of $\tG p$ (possibly a limit group)  has attached
component and cusp graphs,
$\sC^f_{G,\bfC,p}(G')$ and $\sT^f_{G,\bfC,p}(G')$, by running over Nielsen classes corresponding to quotients of $G'$. 

\subsubsection{Setup for the (strong) Main Conjecture} \label{setUpSM}  Suppose  $F_u$ is free of rank $u$ and $J$ is finite acting
faithfully on $F_u$. Consider $F_u\xs J$, and let  $\bfC=(\row \C r)$ be conjugacy classes in $J$. (Our examples use
$r=4$.)  

Form $\tilde F_{u,p}$, the pro-$p$, pro-free  completion of $F_u$. Then  
$\Phi^t=\Phi^t_p$  is the $t$th Frattini subgroup of $\tilde F_{u,p}$ (\S\ref{gpStatements}). Consider two sets $P_\bfC$ and $P_\bfC'$ of
primes, with each consisting of those $p$ with  
$\tilde F_{u,p}/\Phi^1\xs J$ not $p$-perfect, or $p$ has this (respective) property:  

\begin{itemize} \item  $P_\bfC$: $p\mid (p,|J|)\not = 1$.   
\item $P_\bfC'$: $p|\ord(g)$  some $g\in \bfC$. \end{itemize}

For
$p\not\in P_\bfC$, denote (finite) $J$ quotients 
of $\tilde F_{u,p}$ covering $(\bZ/p)^u$  by 
$\sV_p(J)$. 
 
\begin{prob} Which $\ni(V\xs J, \bfC)^{\inn}$ are nonempty,  $p\not\in P_{\bfC}$ and  $V\!\in
\sV_p(J)$.  \end{prob} 

For $p\not\in P_\bfC$, form the collection $\sG_{J,p}$ of limit groups over nonempty 
Nielsen
classes (Def.~\ref{limitGps}). The 
$P_\bfC'$ version of this forms characteristic $p$-Frattini quotients of $F_u\xs J$ where $p$ may divide the order
of $J$, but not the orders of elements in $\bfC$.  

By taking $F_u$ of rank 0 ($u=0$), the $P_\bfC'$ version includes the weak Main Conjecture as a special case
of the strong Main Conjecture \ref{MainConj}.  

We also  must consider the finite $J$ quotients $V$ of $\tilde F_{u,p}$ 
where we  ask only that $V$ is nontrivial.  Denote this set by $\sV_p'(J)$.  

\begin{prob} \label{braidOrbitLimit} What are the  $G^*\in \sG_{J,p}$, $p\not\in P_\bfC'$ (or just in $P_\bfC$)?
What are the 
$H_4$ (braid) orbits on $\ni( G^*,\bfC)^{\inn}$? \end{prob}    We say $ G^*\in \sG_{J,p}$ is a $\bfC$ {\sl
$p$-Nielsen  limit}. If
$O$ is a braid orbit in $\ni(G,\bfC)$  we may  consider only maximal paths (branches) over $O$. Then, maximal groups are 
$p$-Nielsen limits  {\sl through\/} 
$O$ ($\bfC$ is now superfluous). So a cusp or component branch through $O$ defines a $p$-Nielsen limit through $O$. Extend this to consider
$p$-Nielsen limits through any nonempty braid orbit on
$\ni(G',\bfC)$, 
$G'$ any $p$-Frattini cover  of $G$. 

\subsubsection{Examples: $u=2$, $|J|$ is 2 or 3} \label{twoExs} Take $F_u=\lrang{x_1,x_2}$, Our two examples in \eqref{F2Zn} illustrate  limit
Nielsen classes, and the questions we pose.  
 
\begin{edesc} \label{F2Zn} \item  \label{F2Zn2}  $\bZ/2$ case: $J=J_2=\bZ/2=\{\pm 1\}$; $-1$ acts on generators of $F_2$ by  $x_i\mapsto
x_i^{-1}$,
$i=1,2$; and $\bfC=\bfC_{2^4}$ is 4 repetitions of -1. 
\item  \label{F2Zn3} $J_3=\bZ/3=\lrang{\alpha}$; $\alpha$  maps $x_1\mapsto
x_2^{-1}$, 
$x_2\mapsto x_1x_2^{-1}$; and $\bfC\!=\bfC_{\pm 3^2}$ is two repetitions each
of $\alpha$, $\alpha^{-1}$. 
\end{edesc}

The apparent simplicity of \eql{F2Zn}{F2Zn2} is misleading: It is the Nielsen class behind Serre's
Open Image Theorem (\cite[\S 6]{Fr05} explains this). The result (in App.~\ref{F2Z2}) is that $\ni(V\xs J_2,\bfC)$
is nonempty precisely when  
$V\in \sV'_{\bfC_{2^4}}$ is abelian.  

App.~\ref{F2Z3} shows all Nielsen classes in \eql{F2Zn}{F2Zn3} are 
nonempty because they contain H-M reps.~(a special case of Princ.~\ref{FP2}). That is, there are infinite component branches. Yet, it remains a
challenge to Prob.~\ref{braidOrbitLimit}. 

\begin{prob} \label{onlyHM} Let $K$ be any number field. Are all infinite $K$ component branches of $\sT_{(\bZ/p)^2\xs \bZ/3,\bfC_{\pm 3^2},p\ne
3}$, case
\eql{F2Zn}{F2Zn3},  defined by H-M rep.~cusp branches?
\end{prob} Prop.~\ref{H3NC} gives an infinite limit group not equal to $\tilde F_{2,2}\xs J_3$: H-M cusp branches don't give
 all  infinite component branches of $\sT^f_{(\bZ/p)^2\xs
\bZ/3,\bfC_{\pm 3^2},p\ne 3}$. 

\begin{rem} It is essential for the RIGP (\S\ref{MTview}) that we consider questions like 
Prob.~\ref{onlyHM} for all $r$, based on Conj.~\ref{gp'givesPSC}. \end{rem}   

\subsection{The small lifting invariant} \label{smallInv} Let 
\!$G$ be  finite,   \!$\psi: \!R\to G$ a Frattini
central extension, and \!$\bfC$ conjugacy classes of \!$G$ with 
elements of order prime to \!$|\ker(\psi)|$. 

For 
$\bg\in \ni(G,\bfC)$, we have a {\sl small lifting invariant\/} $s_\psi(\bg)=s_{R/G}(\bg)=s_{R}(\bg)$ (notation of \S\ref{pperfect}): 
Lift $\bg$ to $\hat \bg\in \bfC$ regarded as conjugacy classes in $R$ and 
form $\Pi(\hat \bg)\in \ker(R\to G)$. It is an invariant on the braid orbit $O=O_\bg$ of $\bg$  which we call $s_R(O)$ \cite[Part
III]{Fr4}. When
$\ker(R\to G)=\SM_{G,p}$, denote this  
$s_{G,p}(O)$. 

At times we regard $\ker(R\to G)$ as a multiplicative (resp.~additive) group: So,  $s_{G,p}(O)=1$ (resp.~$s_{G,p}(O)=0$) when the
invariant is trivial.   

\subsubsection{Component branch obstructions} 
Consider a nontrivial  Frattini central  cover 
$R'\to G'$ through which $\tG p\to G_0$ factors. Then, $\ker(R'\to G')$ is a quotient of
the Schur multiplier of $G'$ (\S\ref{pperfect}). Denote the  collection of such covers 
${\tsp S} {\tsp M}_{G,p}$, and the subcollection of $R'\to G'$ that are a subfactor of $G_{k+1}\to G_k$ with the notation ${\tsp S} {\tsp
M}_{G,p,k}$.  Suppose $(G,\bfC,p)$ satisfies the usual \MT\ conditions.
  
\begin{lem} \label{liftLem} In the above notation for $R'\to G'\in \!{\tsp S}\! {\tsp M}_{G,p,k}$ these are equivalent: 
\begin{itemize} \item the injection from braid orbits in 
$\ni(R',\bfC)$ to braid orbits in  $\ni(G',\bfC)$ has  
$\bg\in \ni(G',\bfC)$ in its image; \item and $s_{R'}(\bg)=1$. \end{itemize}  

For each $k\ge 0$, braid orbits  in $\ni(G_{k+1},\bfC)$ map  onto 
compatible systems of braid orbits $O$ on $\ni(G',\bfC)$ with $R'\to G'\in {\tsp S} {\tsp M}_{G,p,k}$ and
$s_{R'}(O)=1$.  

Similarly, infinite branches of 
$\sT_{G,\bfC,p}$ map onto compatible systems of braid orbits $O$ in $\ni(G',\bfC)$ with $R'\to G'\in {\tsp S} {\tsp M}_{G,p}$  and
$s_{R'}(O)=1$; and this is 
one-one.  
\end{lem} 

\begin{proof}[Comments] Given $\bg\in  \ni(G',\bfC)$ there is a unique lift to $\hat \bg\in (R')^r\cap \bfC$, and $\hat \bg\in \ni(R',\bfC)$ if
and only if $s_{R'}(\bg)=1$. This shows the first paragraph statement. 

Consider any cover $H''\to H'$ through which $\tG p\to G$
factors. We can always refine it into a series of covers to assume $\ker(H''\to H')=M'$ is irreducible (as an
$H'$ module). For asking when braid orbits on $\ni(H'',\bfC)$ map surjectively to braid orbits on $\ni(H',\bfC)$ it suffices to assume
$M'$ is irreducible. \cite[Obst.~Lem.~3.2]{FrK97} says  the map $\ni(H'',\bfC)\to
\ni(H',\bfC)$ is surjective unless $M'$ is the trivial $H'$ module. So, we have only to check surjectivity in those cases, using the lifting
invariant. That establishes the second paragraph statement. 

\S\ref{level1p2} uses $k=1$ for $(A_4,\bfC_{\pm 3^2}, p=2)$ to show the braid orbit map of the second paragraph is not necessarily
one-one. This is from their being two orbits of H-M reps.~in $\ni(G_1(A_4),\bfC_{\pm 3^2})$.   

The only point needing further comment is why the onto map of the last paragraph statement is one-one. That is because the collection of $G'$
with $R'\to G'\in {\tsp S} {\tsp M}_{G,p}$ is cofinal in all quotients of $\tG p$:   Prop.~\ref{densityOnes}.  
\end{proof}  

Frattini Princ.~\ref{FP3} relates cusp branches (on $\sC_{G,\bfC,p}$) to component branches. 
This is a tool for considering if there is an o-$p'$ cusp branch lying over a given o-$p'$ cusp. Resolving
Conj.~\ref{noWeigB} is crucial to deciding what are the infinite \MT\  component branches.  Though elementary, Lem.~\ref{GFaction} is a
powerful principle.  

\begin{princ} \label{tangBPPrinc} Suppose $B'$ is a component branch on $\sT_{G,\bfC,p}$. The only way we can now prove
$G_F$ has a finite orbit on $B'$ (the hypothesis of \eqref{lowgenus}) is to find a cusp branch $B$ that defines $B'$ for which, modulo braiding,
$G_F$ has a finite orbit on
$B$. Further, all successes here are with g-$p'$ branches. \end{princ} 

\S\ref{nogp'} has \MT s with no  g-$p'$ cusps where we don't yet know if they have infinite component branches. 
Conj.~\ref{noWeigB} says they should not.  Prob.~\ref{onlyHM} ($\bZ/3$ rank 2 \MT) has 
similar challenges for Conj.~\ref{gp'givesPSC}: Do g-$p'$ cusps define all infinite component branches. 

\subsubsection{Replicating obstructed components} \label{repComps} Thm.~\ref{highObst} gives \MT s with at least two components at every
level. One is an H-M component with an obstructed component (Def.~\ref{obstMT}) lying above it at the next level ($k\ge 0$). 

Suppose
$\psi_0: R_0\to G_0$ is a Frattini central extension of
$G_0$ with
$\ker(\psi_0)=\bZ/p$: a $\bZ/p$ quotient as in \S\ref{schurQuot}. Further, suppose
$\psi_1:R_1\to G_1$ is a Frattini central extension of $G_1$ with $\ker(\psi_1)$ also a
$\bZ/p$ quotient, but  {\sl antecedent\/} to $\ker(R_0\to G_0)$. This means: $\ker(\psi_1)=\lrang{\tilde a^p}$ with 
$\tilde a\in
\ker(\tG p\to G_0)$  a lift  of a generator of $\ker(\psi_0)$. 

The idea of antecedents generalizes in the following technical lemma. It will seem less technical from the proof by recognizing $M_k'$
interprets as
$M_0$ multiplied by $p^k$. 

\begin{lem} \label{ppowermap} Then, $\ker(R_1\to G_0)$ is a $\bZ/p^2[G_0]$ module. For all $k\ge 1$, there is a
Frattini cover $\psi_k^*:R_k^*\to G_k$ with $\bZ/p^2[G_0]$ acting on $\ker(\psi_k^*)$ isomorphic to its action on $\ker(R_1\to G_0)$. Also, 
$\psi_k^*$ factors through a cover
$G_k^*\to G_k$ with $G_0$ acting on $\ker(G_k^*\to G_k)=M_k^*$ as it does on $M_0$. Further: 
\begin{edesc} \label{repAct} \item \label{repActa}  $M_k^*$ is a quotient of  $M_k$ (\S\ref{schurQuot}) on which $G_k$  acts
through
$G_0$; and  \item \eql{repAct}{repActa} extends to a $\bZ/p^2[G_k]$ action on $\ker(R_k^*\to G_{k})$ that factors through 
$\bZ/p^2[G_0]$ acting on 
$\ker(R_1\to G_0)$. 
\end{edesc} 
\end{lem} 

\begin{proof} The condition that $\ker(R_1\to G_1)$ is a $\bZ/p^2[G_0]$ module is the main condition for an antecedent Schur multiplier, 
part of the characterization of that condition in \cite[Prop.~4.4]{Fr02}. 

The lemma says the $\bZ/p^2[G_0]$ module $\ker(R_1\to G_0)$ \lq\lq replicates\rq\rq\ at 
all levels. It comes from forming the {\sl abelianization\/} $\tG p/(\ker_0,\ker_0)\eqdef \tG p'$ of $\tG p\to G_0$ (as in \S\ref{abelgen} and used
many times in such places as \cite[\S4.4.3]{BFr02}). 

Denote the characteristic Frattini quotients of $\tG p'$ by $\{G_k'\}_{k=0}^\infty$. Then, $M_0$ still identifies naturally with 
$\ker(G_1'\to G_0)$. Since $\ker_0/(\ker_0,\ker_0)=\ker_0'$ is abelian, taking all $p$th powers (additively: image of multiplication by $p^k$)
in $\ker_0'$ gives the $k$th iterate of its Frattini subgroup $\ker_k'$. Then,
$M_k'$ is the 1st Frattini quotient of $\ker_k'$. Since $G_0$ acts on $\ker_k'$ this induces an action on $M_k'$. As $\ker(R_1\to G_0)$ is
also abelian, this replicates at level $k$ as $R_k'$, also by \lq\lq multiplication by $p^k$.\rq\rq 

The conclusion of the lemma follows from recognizing, inductively from the universal $p$-Frattini property, that $R_{k+1}\to G_{k}$ must factor
through $R_{k+1}'\to G_k'$, giving $R_{k+1}^*$ as the pullback over $G_k$ of $R_{k+1}'\to G_k'$, etc. 
\end{proof} 

Continue the notation of Lem.~\ref{ppowermap}. We use it to replicate the event of having two components, one an H-M component,  at one \MT\
level to higher tower levels. 
\begin{thm} \label{highObst} 
Let   ${}_0\bg=({}_0g_1^{-1},{}_0g_1,{}_0g_2,{}_0g_2^{-1})\in \ni(G_0,\bfC)$ be an H-M rep. As in Princ.~\ref{FP2}, take
$\{{}_k\bg\}_{k=0}^\infty$ to define an H-M cusp branch above ${}_0\bg$.  Assume there is  a level 1 braid orbit  represented by
${}_1\bg'\in
\ni(G_1,\bfC)$ with these properties: $${}_1\bg' \mapsto {}_0\bg\text{ and }s_{R_1}({}_1\bg')\ne 1.$$ 

Then, there is a sequence  $\{{}_k\bg'\in \ni(G_k,\bfC)\}_{k=1}^\infty$ with
${}_k\bg'$ lying over ${}_{k-1}\bg$ and 
$s_{R_k^*}({}_k\bg')\ne 1$. Finally, we don't need to start these statements at level 0; they apply for $k\ge k_0$, if the hypotheses hold 
replacing 
$(G_0,M_0,R_0,R_1))$ with                                            
$(G_{k_0},M_{k_0},R_{k_0},R_{k_0+1})$. 
\end{thm}  

\begin{proof} As in \cite[\S9]{BFr02}, with no loss assume $${}_1\bg'=({}_1g_1^{-1},a_1({}_1g_1)a_1^{-1},a_2^{-1}({}_1g_2)a_2,{}_1g_2^{-1})$$
with 
$a_1,a_2\in M_0$ the images of $\hat a_1,\hat a_2\in R_1$ lying respectively over them. A restatement of $s_{R_1/G_1}({}_1\bg')\ne 1$
(multiplicative notation) is this: 
\begin{triv} \label{repas} $a_1^{{}_0g_1}a_1^{-1}a_2^{-1} a_2^{{}_0g_2^{-1}}=1$, but $\hat a_1^{{}_0g_1}\hat a_1^{-1}\hat a_2^{-1} \hat
a_2^{{}_0g_2^{-1}}\ne 1$. 
\end{triv} Now let $a_1,a_2$ represent their respective images in $M_k^*$ and replace ${}_1g_1$ and ${}_1g_2$ in \eqref{repas} by  ${}_kg_1$ and
${}_kg_2$. This produces the ${}_k\bg'$ in the theorem's statement. The corresponding expressions in \eqref{repas} hold because we have a
$\bZ/p^2[G_k]$ isomorphism of $\ker(\psi_k^*)$ with $\ker(\psi_0^*)$. 

The final statement applies the general principle that we can start a \MT\ at any level we want just by shifting the indices.
\end{proof} 

\begin{exmp}[Several components at high levels] \label{highObstExs} \cite[Prop.~9.8]{BFr02} shows level 1 of the $(A_5,\bfC_{3^4},p=2)$ \MT\
has exactly two components, and these  satisfy the hypotheses of Thm.~\ref{highObst} (more in Ex.~\ref{A5limitgps}). Thus, each level $k\ge 1$
of this \MT\ has at least two components. (Level 0 has just one.)  

Level 1 of the
$(A_4,\bfC_{\pm 3^2},p=2)$ \MT\  has two H-M  and four other components, each over the H-M component (from two at level 0; see
\S\ref{level1p2}). Thm.~\ref{highObst} lets us select whatever H-M cusp representatives we want over ${}_0\bg$. So, suppose there are several
braid orbits of H-M branches, and the hypothesis at one level holds. Then, each braid orbit of an H-M cusp branch through that level gives a
pair of components at  higher levels.    Thus,
Thm.~\ref{highObst} says each level $k\ge 2$ of the $(A_4,\bfC_{\pm 3^2},p=2)$ \MT\ has at least eight components. 
\end{exmp} 

\subsection{Weigel's $p$-Poincar\'e Duality Theorem} \label{gp'Weig}
Let $\phi:X\to\prP^1_z$, with branch points
$\bz$,  be a Galois cover in  $\ni(G,\bfC)^\inn$ representing a braid orbit $O$.

 With
$U_\bz=\prP^1_z\setminus \{\bz\}$, use {\sl classical
generators\/} 
$\row x r$ to describe the fundamental group $\pi_1(U_\bz,z_0)$: $\row x r$ (in order corresponding to branch points
of $\phi$, 
$\row z r$) freely generate it, modulo the product-one relation
$\prod_{i=1}^r x_i$ \cite[\S1.2]{BFr02}. Restrict $\phi$ off $\bz$ to give  $\phi^0:X^0\to U_\bz$. 
Let
$\bg\in
\ni(G,\bfC)$ be the corresponding branch cycles giving a representing homomorphism $\pi_1(U_\bz,z_0)\to G$ by $x_i\mapsto g_i$, $i=1,\dots, 
r$. 

Denote the pro-$p$ completion of the fundamental group of the (compact) Riemann surface $X$ by $\pi_1(X)^{(p)}$. 
\cite[Prop.~4.15]{BFr02} produces a quotient
$M_\phi$ of $\pi_1(U_\bz,z_0)$ with    
$\ker(M_\phi\to G)$  identifying with $\pi_1(X)^{(p)}$  (proof
of Lem.~\ref{extLem}). 

We sometimes denote $M_\phi$ by $M_\bg$ when given $\bg\in \ni(G,\bfC)$ defined by  classical generators. Lem.~\ref{extLem} says \wsp up to
braiding \wsp   
$M_\bg\to G$ is independent of $\bg$. Since
$\ker(M_\bg\to G_k)$ is a pro-$p$ group, the notation
$\ni(M_{\bg},\bfC)$ makes sense (as in
\S\ref{gpStatements}).  

\begin{lem} \label{extLem} The action of $H_r$ on $\bg$ is
compatible with its action on $\row x r$. This gives a braid orbit of homomorphisms starting with $M_\bg\to G$. As abstract group extensions
they are isomorphic.   

Also,
$p$-Nielsen limits through
$O$ are maximal among  quotients of $\tG p$ through which $M_\bg\to G$ factors (up to conjugation by  $\ker(M_\bg\to G)$). So, $O$ starts 
a component branch of $\sT_{G,\bfC,p}$ if and only if, running over $R'\to G'\in {\tsp S} {\tsp M}_{G,p}$ (as in
Lem.~\ref{liftLem}),  each  $\psi_{G'}: M_\bg\to G'$ extending
$M_\bg\to G$ extends to $\psi_{R'}: M_\bg \to R'$. 

The obstruction to extending $\psi_{G'}$ to $\psi_{R'}$ is the image in $H^2(M_\bg,\ker(R'\to G'))$ by
inflation of $\alpha\in H^2(G',\ker(R'\to G'))$ defining the extension $R'\to G'$. \end{lem}
 
\begin{proof}[Comments] Let
$W$ be the normal subgroup of
$\pi_1(U_\bz,z_0)$ generated by $x_i^{\ord(g_i)}$, $i=1,\dots,r$.  Identify $U=\ker(\pi_1(U_\bz,z_0)/W\to G)$ with $\pi_1$ 
of  $X$; what Weigel calls a {\sl  finite index surface group}  \cite[Proof of Prop.~5.1]{Weigel2}. (If $\phi$ is not a Galois cover, then it
is more complicated to describe $\pi_1(X)$ by {\sl branch cycles\/} \cite[p.~75--77]{CombComp}.) 

In Weigel's notation, 
$\Gamma=\pi_1(U_\bz,z_0)/W$. Form $M_\bg$ by completing $\Gamma$ with respect to $\Gamma$ normal subgroups in $U$ of index (in $U$) a power of
$p$. For more details see \S\ref{oneliftinv}. Then $M_\bg$ has a universal property  captured in the second paragraph of the lemma. 

In a characteristic 0 smooth connected family of covers the isomorphism class of the monodromy group does not change. That is, the braiding of
$\bg\in
\ni(G,\bfC)$ to
$\bg'$ from a deformation of the cover with branch point set $\bz_0$ over a path in $\pi_1(U_r,\bz_0)$ produces another copy of $G$.  The
same is true if you apply this to a profinite family of covers defining a cofinal family of quotients of $M_\bg$. This shows that braiding
induces an isomorphism on $M_\bg$ as said in the first paragraph of the lemma. 

This gives the first
paragraph statement. The final paragraph statement is likely well-known. See, for example, 
\cite[Prop.~2.7]{Fr4} or \cite[Prop.~3.2]{Weigel2}.    
\end{proof} 

We continue notation of Lem.~\ref{extLem}. The following translates \cite{Weigel2} for our group $M_\bg$. We 
explain terminology and  module conditions for later use. 

\begin{thm} \label{WeigelThm}  $M_\bg$ is a dimension 2 oriented $p$-Poincar\'e duality group. 
\end{thm} 

\begin{proof}[Comments] 
The meaning of the phrase (dimension 2) $p$-Poincar\'e duality
is in \cite[(5.8)]{Weigel2}. It expresses an exact cohomology pairing  \begin{equation} \label{duality} H^k(M_\bg,U^*)\times H^{2-k}(M_\bg,U)\to
\bQ_p/\bZ_p\eqdef I_{M_{\bg},p}\end{equation} where $U$ is any  abelian $p$-power group  that is also a $\Gamma=M_\bg$ module, $U^*$ is its dual
with respect to
$I_{M_{\bg},p}$ and 
$k$ is any integer.
\cite[I.4.5]{SeGalCoh} has the same definition, though that assumes in place of 
$M_\bg$ a pro-$p$-group. By contrast,   $M_{\bg}$ is $p$-perfect, being generated by $p'$ elements (Lem.~\ref{perfnongen}). 
In the extension problems of \S\ref{weigCor}, the quotients of $M_{\bg}$ that interest us are Frattini covers of $G$, so also 
$p$-perfect. 

\cite[p.~38]{SeGalCoh} points to Lazard's result  that a
$p$-adic analytic group of dimension
$d$ (compact and torsion-free) is a Poincar\'e group of dimension $d$. Since, however, our group is residually pro-free, it isn't
even residually $p$-adic analytic. 

Weigel's result is for general Fuchsian groups $\Gamma$, and the {\sl dualizing
module\/}  $I_{\Gamma,p}$,   may not be the same as in this example. It is classical that 
$\pi_1(X)$ (and $\pi_1(X)^{(p)}$) satisfies Poincar\'e duality.  \cite[Chap.~VIII, \S3, Remark]{brown} interprets this exactly as the discussion
of
\S\ref{pperfect} suggests for group cohomology. 
\cite[Prop.~18, p.~25]{SeGalCoh} applies Shapiro's Lemma to show  a dualizing module that works for $\Gamma$ also works for every open 
subgroup. Most of Weigel's proof establishes the converse: That the $I_{\pi_1(X),p}$ used here does act as a dualizing module for
$M_\bg$. 
\end{proof} 

\begin{rem}[Addendum to Lem.~\ref{extLem}] Suppose two extensions $M_{\bg_i}\to G$,  arise from $\bg_i\in \ni(G,\bfC)$, $i=1,2$. Further,
assume they are isomorphic. Then, it is still possible they are not braid equivalent, though examples aren't easy to come by. We allude to
one in \eql{apps}{appsa}: Two extensions corresponding to the two H-M components called $\sH_1^{+,\beta}$, $\sH_1^{+,\beta^{-1}}$. 
The group $G$ in this case is $G_1(A_4)$. It has an automorphism mapping $\bg_1$ to $\bg_2$, giving elements in different braid orbits. 
Since these are H-M components, Princ.~\ref{FP2} gives isomorphic  extensions
$M_{\bg_i}\to
\tG p$,
$i=1,2$ (Princ.~\ref{FP2}) in distinct braid orbits. 
\end{rem} 
 
\subsection{Criterion for infinite branches on $\sT_{G,\bfC,p}$} \label{weigCor} Cor.~\ref{weigCor1} reduces 
finding infinite component branches on $\sT_{G,\bfC,p}$ through a braid orbit (as in \S\ref{schurQuot}) to a sequence of small lifting
invariant checks from the Schur multiplier of each $G_k$, $k\ge 1$.  Cor.~\ref{weigCor2} is 
our major test for  when we have a limit group. 

\subsubsection{One lifting invariant checks unobstructed braid orbits} \label{oneliftinv} This subsection regards the small lifting invariant in
additive notation. Let $O_k \le \ni(G_k,\bfC)$ be a braid orbit and ${}_k\bg$ a representative of this orbit. The cardinality of the fiber in
\eql{cuspP}{cuspPb} over $O_k$ is the degree of a level $k+1$ \MT\ component over its level $k$ image defined by $O_k$. This is a braid
invariant. Cor.~\ref{weigCor1} is (at present) our best test for 
when it is nonempty, unless $\bg$ braids to a g-$p'$ representative (Princ.~\ref{FP2}). 

We may consider $M_\bg$ as a completion of a group,
$D_{\bar \psigma}$, presented as $\lrang{\row {\bar \sigma} r}$ modulo the normal subgroup generated by $\bar\psigma\eqdef \{\bar
\sigma_i^{\ord(g_i)},i=1,\dots,r,\text{ and } \bar \sigma_1 \cdots \bar \sigma_r\}$. Let $K_{\bar \psigma^*}$ by the group from removing
the quotient  relation $\bar \sigma_1
\cdots \bar \sigma_r=1$. Denote  corresponding generators of it by $\row {\bar \sigma^*} r$. Then,  
the cyclic groups $\lrang{\bar \sigma_i^*}/(\bar \sigma_i^*)^{\ord(g_i)}$, $i=1,\dots r$, freely generate $K_\psigma$. 

Complete $K_{\bar \psigma^*}$ with respect to  $p$-power index  subgroups
of
$\ker(K_{\bar \psigma^*} \to G)$, normal in $K_\psigma$, calling the result $\tilde K_{\bar \psigma^*}$ (forming a
natural surjection $\psi_{\bar
\psigma^*}: \tilde K_{\bar \psigma^*} \to M_\bg$). 

\begin{lem} \label{univpprojcover} Mapping the $\tilde K_{\bar \psigma}$ generators $\row {\bar \sigma^*} r$,
in order, to entries of ${}_k\bg$, gives a  homomorphism $\mu_k: \tilde K_{\bar \psigma^*} \to G_k$. If  $\row {h^*} r\in
\bfC\cap G_{k+1}^r$ lie respectively over entries of ${}_k\bg$, then the surjective homomorphism $\mu_{k+1}: \tilde K_{\bar \psigma^*} \to
G_{k+1}$ mapping $\bar \sigma_i^*\mapsto h_i^*$, $i=1,\dots,r$, extends 
$\mu_k$. \end{lem} 

\begin{proof}[Comments] The construction is geometric: Remove an additional point $z'$ from $U_\bz$ to get $\pi_1(U_{\{\bz,z'\}},z_0)$.
We can identify this with notation coming from $D_\psigma$, as the group freely generated by $\bar\psigma$. This identifies $K_{\bar \psigma^*}$
with  its description above. It also identifies $\ker(K_{\bar \psigma^*} \to G_0)$ with the fundamental group of $X'\eqdef X\setminus
\{\phi^{-1}(z')\}$. As   $X'$ is a projective curve with a nonempty set of punctures, this is a free group.  \end{proof} 

\begin{rem}[Addendum to proof of Lem.~\ref{univpprojcover}] The group $M_\bg$ is not $p$-projective. Yet, here is why its cover $\tilde K_{\bar
\psigma^*}$ is.  For $P$  a $p$-Sylow of $G$, we can identify a $p$-Sylow of
$\tilde K_{\bar
\psigma^*}$ with the pro-$p$ completion of the free group $\pi_1(X'/P)$. A profinite group with pro-$p$ $p$-Sylow is
$p$-projective (\cite[Prop.~22.11.08]{FrJ86}, in new edition). \end{rem}
 
\subsubsection{Two obstruction corollaries} \label{twoCorollaries}  Continue the discussion of \S\ref{oneliftinv}. If $\bg\in O_k$,  
then it defines a cover $\psi_\bg: 
M_\bg \to G_k$.   A paraphrase of 
Cor.~\ref{weigCor1} is that if 
$\psi_\bg$ is {\sl obstructed\/} at level $k$ then it is  by some
$\bZ/p$ 
 quotient of  $\ker(G_{k+1}\to G_k)$. Cor.~\ref{weigCor2} tells us precisely what are the exponent $p$ Frattini extensions of a limit group.  

\begin{cor} \label{weigCor1} The fiber over $O_k$ is empty if and only if there is some central Frattini extension $R\to G_k$ with kernel
isomorphic to $\bZ/p$ for which $\psi_\bg$ does not extend to $M_\bg\to R\to G$. 
  \end{cor}

\begin{proof} In the notation of \S\ref{schurDisc} we only need to show this: If the fiber of
\eql{cuspP}{cuspPb} is empty, then 
$s_{R/G_k}(\bg)\ne 0$ for some $\bZ/p$ quotient $R/G_k$ of the first Loewy layer of $M_k$. \cite[Prop.~2.7]{Fr4} says $H^2(G_k, M_k)=\bZ/p$: It
is 1-dimensional. Lem.~\ref{extLem} says the obstruction to lifting $\psi$ to $G_{k+1}$ is the inflation of some fixed generator 
$\alpha\in H^2(G_k, M_k)$ to
$\tilde \alpha\in H^2(M_{\bg}, M_k)$.  

Though $\tilde \alpha$ may seem abstract, the homomorphism  $\mu_{k+1}$ of Lem.~\ref{univpprojcover} allows  us to form an 
explicit cocycle for the obstruction to lifting $M_\bg\to G$. For each  $\bar g\in M_\bg$ choose $h_{\bar g} \in G_k$ as the image in $G_k$ of
one of the elements of $\tilde K_{\bar \psigma^*}$ over $\bar g$. Now compute from this the 2-cocycle 
$$\tilde \alpha({\bar g_1,\bar g_2})=h_{\bar g_1} h_{\bar g_2} (h_{\overline {g_1g_2}})^{-1}, \bar g_1,\bar g_2\in  M_\bg$$ describing the
obstruction.  Since $\psi_{\bar \psigma^*}$ is a homomorphism, the only discrepancy between  $\alpha(\bar g_1,\bar g_2)$ and the identity is
given by the leeway in  representatives for $h_{\overline {g_1g_2}}$ lying over $\overline{g_1g_2}$. So, the cocycle  $\tilde \alpha(\bar
g_1,\bar g_2)$ consists of words in the kernel of  $K_{\bar \psigma^*}\to M_\bg$, and it vanishes if and only if it is possible 
to choose $(h_1^*,\dots,  h_r^*)$ (as in the statement of Lem.~\ref{univpprojcover}) to satisfy $h_1^*\cdots h_r^*=1$. 

By \eqref{duality} duality, $H^2(M_{\bg}, M_k)$ has a perfect pairing with $H^0(M_{\bg}, M_k^*)$, that initially goes into
$H^2(M_{\bg},I_{M_{\bg},p})$ by applying an element of $H^0(M_{\bg}, M_k^*)$ to the values of a 2-cycle in $H^2(M_{\bg}, M_k)$. Identify
$H^0(M_{\bg}, M_k^*)$ with  $$H_0(M_{\bg}, D \otimes M_k) \simeq D
\otimes_{\bZ/p[M_{\bg}]} M_k,$$ with $D=\bZ/p$ the duality module for $\bZ/p[M_{\bg}]$  (on which it acts trivially).  Hence, the tensor
product $D \otimes_{\bZ/p[M_{\bg}]} M_k$ is canonically isomorphic to the maximal quotient of $M_k$ on which $M_{\bg}$ (and therefore
$G_k$) acts trivially \cite[p.~98]{AtWall}. That is, $D \otimes_{\bZ/p[M_{\bg}]} M_k$ identifies with the kernel of the maximal central
exponent $p$ extension of $G_k$. 

Now we check the value of the pairing of 
$\tilde \alpha(\bullet,\bullet)\in H^2(M_{\bg}, M_k)$ against an element  $\beta\in H^0(M_{\bg}, M_k^*)$. Further, regard $\beta\eqdef \beta_R$
as the linear functional on
$M_k$ from $\ker(G_{k+1}\to R)$, with $R\to G_k$  a central extension defining a $\bZ/p$ quotient, as above. 

Being very explicit, this says the value of $\beta_R$ on $\tilde \alpha$ is  the lifting invariant $s_{R}(\bg)$ for the image $\bg$
of
$(\row {h^*} r)$ in $\ni(G_k,\bfC)$. Since the pairing is perfect, conclude the corollary: The obstruction for extending
$M_\bg\to G_k$ to
$M_\bg
\to G_{k+1}$ is trivial if and only if $s_{R}(\bg)$ is trivial running over all such $R\to G_k$. 
\end{proof}

The proof of the last result also applies to limit groups. 
\begin{cor} \label{weigCor2} If $G^*$ is a limit group in a Nielsen class and a proper quotient of 
$\tG p$, then
$G^*$ has exactly one nonsplit extension by a $\bZ/p[G^*]$ module, and that module must be trivial. \end{cor} 

\begin{proof} Suppose $\bg^*\in \ni(G^*,\bfC)$ represents the braid orbit giving $G^*$ as a limit group (Def.~\ref{limitGps}). From the proof of
Cor.~\ref{weigCor1},  we have only to show there cannot be two
$\bZ/p$ quotients of the exponent
$p$ part of the Schur multiplier of $G^*$. 

Suppose $R_i\to G^*$, $i=1,2$, are two distinct central extensions defining $\bZ/p$ quotients. So, their kernels generate a
2-dimensional quotient of the Schur multiplier of $G^*$. Since $G^*$ is a limit group,  $s_{R_i/G^*}(\bg^*)\ne 0$ generates
$\ker(R_i\to G^*)$, $i=1,2$. 

Apply Thm.~\ref{WeigelThm}: $H^2(M_\bg,\bZ/p)=\bZ/p$.  Let  $\alpha_i\in H^2(M_\bg,\bZ/p)=\bZ/p$ be the inflation of the element of
$H^2(G^*,\bZ/p)$ defining
$R_i$,
$i=1,2$. So there are $p'$ integers $a_i$, $i=1,2$, with $a_1\alpha_1+a_2\alpha_2=0$. Also, $a_1s_{R_1/G^*}(\bg^*)+a_2 s_{R_2/G^*}(\bg^*)\ne
0$ defines a $\bZ/p$ quotient of the Schur multiplier of $G^*$. 

This gives a central extension $R^*\to G^*$, and the inflation of an element
of
$H^2(G^*,\bZ/p)$ to $H^2(M_{\bg^*},\bZ/p)$  defining it is 0. Thus  Lem.~\ref{extLem} contradicts that $G^*$ is a limit group since it says 
$M_{\bg^*} \to G^*$ extends to $M_{\bg^*}\to R^*$. 
\end{proof}

\subsubsection{Why Cor.~\ref{weigCor1} is a global result} \label{exofWeigTest} Consider two (braid inequivalent) extensions of $\psi_i:
M_\bg\to G_{k+1}$,
$i=1,2$, of
$\psi: M_\bg\to G_k$. Assume, hypothetically, the following holds (it does not in general):
\begin{triv} \label{weigErr} There is an
extension of
$\psi_1$ to $\psi_1': M_\bg\to G_{k+2}$ if and only if there is an extension of
$\psi_2$ to $\psi_2': M_\bg\to G_{k+2}$. \end{triv} 

Applying Princ.~\ref{FP2} would then give the following (false) conclusion from \eqref{weigErr}. 
\begin{triv} \label{weigErrConc} If $\bg$ is a g-$p'$ cusp, then any component branch of $\sT_{G,\bfC,p}$ through the braid orbit of $\bg$ is
infinite.
\end{triv}

Cor.~\ref{weigCor1} works with
$G^*$, any group through which
$\tG p\to G_0$ factors,   replacing $G_k$ and  with any $G^*$ quotient $M^*$ of $\ker(G_1(G^*)\to
G^*)$ replacing $M_k$.  (Reminder: $G_1(G^*)$ is the 1st characteristic $p$-Frattini cover of $G^*$.) So,  given an
hypothesis like
\eqref{weigErr},  one might try to reduce the  proof of Cor.~\ref{weigCor1} to where $M^*$ is simple. This would allow stronger conclusions, 
eschewing  considering one integer $k$ at-a-time. 

This, however, is a variant of the false conclusion \eqref{weigErrConc}. 
Examples \ref{level1A5} and \ref{level1A4} show \eqref{weigErrConc} is false. They explain why applying Cor.~\ref{weigCor1} to detect an
infinite branch can't be done by just testing the lifting invariant at one level. These examples \wsp
based on \cite[Chap.~9]{BFr02} \wsp help understand this subtle argument.

Also, for a given \MT\ level
$k$, and
$R'\to G'\in {\tsp S} {\tsp M}_{G,p,k}$ (Lem.~\ref{liftLem}), precise genera formulas for  \MT\ branches require knowing if braid
orbits achieve other lift values than the trivial one. Again, these examples illustrate. They rely on  {\sl centralizer condition\/}
\eqref{centCond}. So, we don't yet know how to generalize them to,  say, replace
$A_n$ by $G_k(A_n)$ for $k$ large, even for the antecedent Schur multiplier  because \eqref{centCond} doesn't hold. 

\begin{exmp}[Level 1 of the $(A_5,\bfC_{3^4}, p=2)$ \MT] \label{level1A5} Here $\bfC_{3^4}$ is four repetitions
of the 3-cycle conjugacy class in $A_5$. \cite[Prop.~9.8]{BFr02} shows there are exactly two braid orbits
$O_1$ and
$O_2$ on
$\ni(G_1(A_5), \bfC_{3^4})$ where $p=2$, both over the unique braid orbit $O$ on $\ni(A_5,\bfC_{3^4})$. The 2-part, $\SM_{G_1(A_5),2}$, 
of the Schur multiplier of $G_1(A_5)$ is $\bZ/2$.  Let $R_1\to G_1$ be the $\bZ/2$ quotient it defines. Then,
$s_{R_1/G_1}(O_1)=0$ and $s_{R_1/G_1}(O_2)\ne 0$.  In fact, $O$ and $O_1$ are orbits of H-M reps. So, at level 1 (but not at level 0) all
possible lift invariants are assumed.  This pure module argument used a strong condition: 
\begin{triv} \label{centCond} The rank of the centralizer in
$M_0=\ker(G_1(A_5)\to A_5)$ of    
$g\in \C_3$ is the same as the rank of $\SM_{G_1,2}$, and $R'\to G_1(A_5)$ is antecedent (\S\ref{schurDisc}). \end{triv}
\end{exmp}

\begin{exmp} \label{level1A4} 
\cite{FrS05}
notes \eqref{centCond} also holds for $(A_4,\bfC_{\pm 3^2},p=2)$  (see \S\ref{3rdMC}; $R'\to G'=G_1(A_4)$
is the antecedent Schur multiplier). The Schur multiplier of $G_1(A_4)$ is $(\bZ/2)^2$. Ad hoc arguments  show we achieve the other two values
of the lifting invariant running over 
$R''\to G'=G_1(A_4)$, with $R''\to G$ the two non-antecedent central Frattini extensions giving $\bZ/p$ quotients. \end{exmp}

\subsection{Weigel branches in $\sC_{G,\bfC,p}$  and Frattini Princ.~3} \label{gp'WeigC}  \cite[Lect.~4]{Ont05}  
generalizes g-$p'$ reps.~to all  $r$. We believe having a g-$p'$ cusp branch $B$ is necessary for an
infinite component branch in $\sT_{G,\bfC,p}$ (Conj.~\ref{gp'givesPSC}). 
Here we approach Conj.~\ref{noWeigB} using multiplicative notation for the small lifting invariant
(\S\ref{smallInv}). 

\subsubsection{Set up for  o-$p'$ cusps} 
We introduce a practicum for deciding if a given o-$p'$ cusp $\bg\in
\ni(G_k,\bfC)$ has an o-$p'$ cusp $\bg'\in \ni(G_{k+1},\bfC)$ over it.  (Compare with the more restrictive search for an
o-$p'$ cusp over a g-$p'$ cusp in \S\ref{cuspExamps}.) From this comes Def.~\ref{defWeigCusp} of a  {\sl Weigel cusp\/}. 
Prop.~\ref{op'overgp'} says there are
\MT s\ where o-$p'$ cusps appear at all high levels. Still, the examples we know do not produce Weigel  branches (projective sequences of such
cusps), so they do not contradict Conj.~\ref{noWeigB}. 

Assume $\bg=(g_1,g_2,g_3,g_4)\in \ni(G,\bfC)$ is an o-$p'$ cusp rep. As $p'$
elements generate $H_{2,3}=\lrang{g_2,g_3}=H$ it is $p$-perfect (Lem.~\ref{perfnongen}). 
Consider diagram \eqref{weigDiag}. The bottom (resp.~top) row has the sequence for the
$p$-representation cover
$R_p'$ of $H$ (resp.~$G$). Pullback of $H$ in $R_p$ is a central extension of  $H$. So, 
a unique  map $\beta_H: R_{p'}\to R_p$ makes \eqref{weigDiag} commutative:  

\begin{equation} \label{weigDiag} \CD 1 \longrightarrow \SM_{G,p}\quad\   @>>>
R_p @>>> G @>>> 1\\
@AAA @AA\beta_H A @AA\text{inj}A\\
1 \longrightarrow \SM_{H,p} @>>> R_p'  @>>> H @>>> 1.\endCD \end{equation} 
Unlike its Lem.~\ref{embedHG} analog, $\beta$ may not be an embedding. Example: Let $H$ be simple, with $\SM_{G,p}\ne\{1\}$  ($p$ odd), and
embed it in an alternating group. The following lemma summarizes this to show compatibility of \eqref{weigDiag} with Lem.~\ref{embedHG}.  

\begin{lem} \label{lemSchur-Frat} Properties of \eqref{weigDiag} apply to any $p$-perfect (or $p'$) subgroup $H\le G$. Further, the map 
$\beta_H$ is compatible with the map $\beta: {}_p\tilde H\to R_p$ defined in Rem.~\ref{remembedHG}. \end{lem}  

\subsubsection{The 3rd Frattini Principle} \label{Weignotgp'} 
Princ.~\ref{FP3}  relates cusp types and lifting invariants for component branches. 
Assume  ${}_0\bg=\bg=(g_1,g_2,g_3,g_4)\in \ni(G,\bfC)$ is an
o-$p'$ cusp rep. Denote a 5th $p'$ conjugacy class containing 
$(g_2g_3)^{-1}$ by $\C_5$. Similarly, its inverse is $\C_5^{-1}$. Denote the collection  $\C_2,\C_3,\C_5$
(resp. $\C_1,\C_4,\C_5^{-1}$) by $\bfC_{2,3}$ (resp.~$\bfC_{1,4}$). Also: $$(g_2,g_3,(g_2g_3)^{-1})={}_0\bg'
\text{ and  }((g_4g_1)^{-1}, g_4,g_1))={}_0\bg'',$$ and let $O_\bg$, $O_{{}_0\bg'}$ and  $O_{{}_0\bg''}$ be the respective braid orbits of
the corresponding Nielsen class representatives.

Assume for some $k\ge 0$, ${}_k\bg\in \ni(G_k,\bfC)$. Let $R_{G_k}\to G_k$ be the central extension of $G_k$
with
$\ker(R_{G_k}\to G_k)$ the maximal quotient of
$M_k$ on which $G_k$  acts trivially. Then, we have similar notation with
$$H_{2,3}({}_k\bg)=H_{2,3}\text{ and }H_{1,4}({}_k\bg)=H_{1,4}$$ replacing
$G_k$. Diagram
\eqref{weigDiag}, with $H=H_{2,3}$, induces maps 
$\beta_{2,3}:R_{H_{2,3}(\bg)}\to R_{G}$ from Lem.~\ref{lemSchur-Frat} as the situation deserves.

\begin{princ}[Frattini Principle 3] \label{FP3} With the previous hypotheses 
\begin{equation} \label{indepCusp} s_{G,p}(\bg)=\beta_{1,4}(s_{R_{H_{1,4},p}}((g_4g_1)^{-1},
g_4,g_1))\beta_{2,3}(s_{R_{H_{2,3},p}}(g_2,g_3,(g_2g_3)^{-1})).\end{equation} 

Suppose ${}_k\bg\in \ni(G_k,\bfC)$ is an o-$p'$ cusp. Consider: 
$$\begin{array}{rl} {}_k\bg'=&({}_kg_2,{}_kg_3,({}_kg_2{}_kg_3)^{-1})\in
\ni(G_k(H_{2,3}(\bg)),\bfC_{2,3}) \text{ and }\\ 
{}_k\bg''=&(({}_kg_4{}_kg_1)^{-1}, {}_kg_4,{}_kg_1))\in
\ni(G_k(H_{1,4}(\bg)),\bfC_{1,4}).\end{array}$$  
Suppose $s_{R_{H_{2,3}({}_k\bg)}}({}_k\bg')=1$ and $s_{R_{H_{1,4}({}_k\bg)}}({}_k\bg'')=1$.   Then, there is an o-$p'$ cusp  ${}_{k+1}\bg\in
\ni(G_{k+1},\bfC)$ over ${}_k\bg$. 

Assume there is an infinite component branch on the $(H_{2,3}(\bg),\bfC_{2,3},p)$ \MT\ over  $O_{{}_0\bg'}$,
and also such a component branch on the $(H_{1,4}(\bg),\bfC_{1,4},p)$ \MT\ over  $O_{{}_0\bg''}$. Then, an o-$p'$ cusp branch gives an infinite 
component branch on the \MT\ over $O_{\bg}$. 
\end{princ}

\begin{proof} Consider the 6-tuple, $\bg^*=((g_4g_1)^{-1},
g_4,g_1, g_2,g_3,(g_2g_3)^{-1})$. This is a juxtaposition of two product-one 3-tuples. Since $(g_4g_1)^{-1}(g_2g_3)^{-1}=1$, we easily see
$s_{G,p}(\bg^*)=s_{G,p}(\bg)$. So,
\eqref{indepCusp} follows from direct computation and the compatibility of the maps $\beta_{2,3}$ and $\beta_{1,4}$ defined in different
places. Lem.~\ref{embedHG} lets us assume $G_{p,k}(H_{2,3})$ and $G_{p,k}(H_{1,4})$ are in $G_{p,k}(G)$. Over
${}_k\bg'$ (resp.~${}_k\bg''$) Lem.~\ref{extLem}  produces  
$${}_{k+1}\bg'\in \ni(G_{k+1}(H_{2,3}(\bg)),\bfC_{2,3}) \text{ (resp.~}{}_{k+1}\bg''\in
\ni(G_{k+1}(H_{1,4}(\bg)),\bfC_{1,4})).$$  

Use Schur-Zassenhaus to produce  
$h\in
\ker(G_{k+1}\to G_k)$ that conjugates   
$$({}_{k+1}g_2'\,{}_{k+1}g_3')^{-1}\in
\C_5 \text{ to }{}_{k+1}g_4'{}_{k+1}g_1'\in \C_5.$$ Replace $({}_{k+1}g_2,{}_{k+1}g_3,({}_{k+1}g_2,{}_{k+1}g_3)^{-1})$ with its conjugate
by the image of 
$h$. So, with no loss,  $({}_{k+1}g_1,{}_{k+1} g_2, {}_{k+1} g_3,{}_{k+1} g_4)$ has product-one,  is in
$\ni(G_{p,{k+1}}(G),\bfC)$ and lies over ${}_0\bg$. This concludes the proof. 

The final paragraph is a simple induction on the previous argument. 
\end{proof}

If Conj.~\ref{noWeigB} holds, then the 3rd paragraph hypotheses of Princ.~\ref{FP3}   can't hold. 

\begin{rem}[Extend Princ.~\ref{FP3}] \label{extendFP3} \cite{FrWProof} has a stronger version of the 2nd paragraph of Princ.~\ref{FP3}: 
If 
$\bg\in \ni(G,\bfC)$ is a rep.~for an o-$p'$ cusp with any two of $s_{R_{H_{2,3},p}}({}_0\bg')$, $s_{R_{H_{1,4},p}}({}_0\bg'')$ and
$s_{R_G,p}(\bg)$ equal 1, then the third is also 1. 
\end{rem} 

\subsection{Evidence for and consequences of no Weigel cusp branches} \label{SectnoWeigB} This subsection considers both evidence for and 
challenges to Conj.~\ref{noWeigB}. 

\begin{defn}[Weigel branch] \label{defWeigCusp} If ${}_k\bg$ satisfies the hypotheses  of Princ.~\ref{FP3}, 2nd paragraph, then we
call
$({}_k\bg)\text{Cu}_4$ a {\sl level $k$ Weigel cusp}. A cusp branch which for large $k$ consists of Weigel cusps   is a {\sl Weigel branch}. 
\end{defn} \noindent We also refer to the component branch in 
$\sT_{G,\bfC,p}$ defined by a Weigel cusp branch as a Weigel component branch. 
 
\subsubsection{Example disappearances of o-$p'$ cusps} \label{disappearop'} For $g\in A_n$ 
of odd order, let $w(g)$ be the sum of $(l^2-1)/8 \mod 2$ over all disjoint cycle lengths $l$ in $g$  ($l\not \equiv \pm 1 
\mod 8$ contribute). \cite[Cor. 2.3]{altGps} has a short proof of Prop.~\ref{cnf.3}  based on  when  $\bfC=\bfC_{3^r}$ is $r$ repetitions of
the 3-cycle class  (guiding the original
statement  in  
\cite{Ser90a}).
 
\begin{prop} \label{cnf.3} Suppose $\bg\in \ni(G,\bfC)$ with $G\le A_n$ transitive, and $\bfC$ consists of conjugacy classes in $G$ with
elements of respective odd orders $\row d r$. Assume also the genus of a degree $n$ cover $\phi: X\to \prP^1_z$ with branch cycles
$\bg$ from this embedding has genus 0. Then, $s_{\Spin_n}(\bg)=(-1)^{\sum_{i=1}^r w(g_i)}$.  
\end{prop} 

At level 0 of the $(A_5,\bfC_{3^4})$ \MT\ 
($p=2$), no cusps are $2$ cusps: Widths are 1,1, 3,3, 5, 5 (\cite[\S2.9.3]{BFr02}; shifts of the cusps of width 1 are H-M reps.). By level 1,
all o-$2'$ cusps  disappear, leaving only g-$2'$ cusps (shifts of H-M reps.) as non-$2$ cusps [BFr02; \S9.1]. Combine this with the comment
before Prop.~\ref{op'overgp'}  for the following. 

\begin{prop} The only infinite cusp branches on the $\sC_{A_5,\bfC_{3^4}, p=2}$ cusp tree are g-$p'$ and $p$ cusp branches.
\end{prop}

\begin{prob} Are there component branches on $\sT_{A_5,\bfC_{3^4}, p=2}$ that contain only $p$ cusp branches? \end{prob}  

\subsubsection{Some Weigel cusps and challenges to Conj.~\ref{noWeigB}} \label{nogp'} 

We give an example Weigel cusp in a Nielsen class containing no  g-$p'$ cusps. Use notation from Ex.~\ref{manyop'} and the representative  
for the Nielsen class $\ni(A_5, \bfC_{\pm53})$ given by $\bg=(g_1,g_2,g_3)$ with $g_1  =  (5\, 4\, 3\, 2\, 1)$ and
$g_2= (2\, 4\, 3\, 5\, 1)$, and
$g_3=(4\,3\,5)$.

There are two
conjugacy classes of 5-cycles in
$A_5$:  $\C_{+5}$ and $\C_{-5}$. Further, if
$g\in
\C_{+5}$, then so is
$g^{-1}$. Let $\bfC_{\pm53}$ denote the collection of conjugacy classes consisting of $\C_{+5}$, $\C_{-5}$ and $\C_3$ (class of a
3-cycle).
\cite[Princ.~5.15]{BFr02} shows
$\ni(A_5,\bfC_{\pm53})$ (absolute or inner) has one braid orbit with lifting invariant +1. By Riemann-Hurwitz, the genus $g$ of absolute 
 covers (degree 5 over $\prP^1_z$) in this Nielsen class is 1, from  $2(5+g-1)=10$. So Prop.~\ref{cnf.3} doesn't
 apply directly. Rather, \cite[\S5.5.2]{BFr02} shows how to compute beyond the genus 0 limitation. Now, take $p=2$. 

This Nielsen class clearly contains no g-$2'$
rep.  Further, similar examples work for any $r\ge 3$ conjugacy classes. For $r\ge 5$: 
juxtapose $\bg \in \ni(A_5,\bfC_{\pm 53})$ with $(g,g^{-1})$ or $(g,g,g)$ ($g\in \C_3$) appropriately. For $r=4$, replace
$\bfC_{\pm53}$ by $\bfC_{\pm 53^2}$. Call the shift (resp. conjugacy classes) of one of these reps.~$\bg'$ (resp.~$\bfC'$).

\begin{res} \label{aWeigCusp} For
$\bfC'=\bfC_{\pm 53^2}$, the natural map $\ni(G_1(A_5), \bfC') \to \ni(A_5,\bfC')$  is onto: no level 0 braid orbit is
obstructed.  The cusp represented by $$\bg''=((3\,4\,5),(5\, 4\, 3\, 2\, 1),(2\, 4\, 3\, 5\, 1),(3\,4\,5))$$ has an o-$p'$ cusp in
$\ni(G_1(A_5),
\bfC')$ over it. So, $\bg''$ is a Weigel cusp. 
\end{res} 

\begin{proof}[Comments] With $R\to A_5$  the $\Spin_5$ cover of $A_5$, $s_{R}(\bg'')=s_{R}((\bg'')\sh)=1$ as we explained above. The
only appearance of $\one_{A_5}$ in $M_0=\ker(G_1(A_5)\to A_5)$ is from $\ker(R\to A_5)$ (\cite[Cor.~5.7]{BFr02} or \cite[Part II]{Fr4})). So,
the hypotheses of Princ.~\ref{FP3}, 2nd paragraph, with $k=0$ apply; and the conclusion does also. 
\end{proof} 

If Conj.~\ref{noWeigB} holds for $\ni(G_1(A_5),\bfC'=\bfC_{\pm 53^2})$ in Res.~\ref{aWeigCusp}, then the conclusion to Prob.~\ref{weigCuspsDis}
is affirmative. 
\begin{prob} \label{weigCuspsDis} Are all o-$p'$ cusps gone at high
levels of the $\ni(G_1(A_5),\bfC_{\pm 53^2})$ \MT? Is it even possible this \MT\ is empty at high levels (agreeing with nonexistence of
infinite component branches having only $p$ cusp branches as in \S\ref{doWeigExist})? 
\end{prob} 

\begin{exmp}[$\ni(A_4,\bfC_{\pm 3^2})$ with $p=2$, see \S\ref{3rdMC}] 
There is an o-$p'$ cusp: 
$$\bg=((1\,2\,4),(1\,2\,3), (1\,3\,4),(1\,2\,4)).$$ Apply the proof of Res.~\ref{aWeigCusp} here. A direct application of Prop.~\ref{cnf.3} \wsp
since the genus 0 hypotheses holds \wsp shows $\beta_{2,3}(s({}_0\bg'))=-1$ while $\beta_{1,4}(s({}_0\bg''))=+1$ (in analogous notation). So,
the 2nd paragraph Prop.~\ref{FP3} conclusion is that the left side of \eqref{indepCusp} is -1, and $\bg$ is not in the image from
$\ni(\Spin_4,\bfC_{\pm3^2})$. \end{exmp} 

\section{Nub of the (weak) Main Conjecture} \label{wcNub}  Use notation, especially for genera, around
\eqref{genusseq}. Assume $B'=\{\sH'_k\}_{k=0}^\infty$ is an infinite branch of $\sT_{G,\bfC,p}$ defined over a number field $K$.  From
Prop.~\ref{nopcenter}, to consider the Main Conj.~we may assume  $G=G_0$ has the $p$-part of its center trivial. We make that
assumption throughout this section. This lets us use the 2nd part of Princ.~\ref{FP1}. 

We show the 
  Main Conj.~\ref{wconj} (for $r=4$) holds unless we are in one of three cases. These we stipulate by listing how $\bar \sH'_{k+1}/\bar
\sH'_k$ ramifies  when
$k>>0$:
\begin{itemize} \item either it doesn't ramify over cusps; \item it is equivalent to a degree $p$ polynomial; \item  or it is equivalent to a
degree
$p$ rational function branched only at two points. \end{itemize} 

\subsection{There should be no $\sT_{G,\bfC,p}$ genus 0 or 1 branches} \label{cuspGrowth}  We must consider two possibilities that would
 contradict the Main Conjecture:
\begin{edesc} \label{genusng} \item \label{genusng0} $g_{\bar \sH'_k}=0$ for all $0\le k < \infty$ ($B'$ has genus 0; $\Ge_{B'}$ consists of
0's); or 
\item \label{genusng1} For $k$ large, $g_{\bar \sH'_k}=1$ ($B'$ has genus 1; almost all of $\Ge_{B'}$ is 1's).\end{edesc} 

\subsubsection{Reduction of the Main Conj.~to explicit cases} \label{redMC} 
An elementary corollary of Riemann-Hurwitz says for $k>>0$,  
\eql{genusng}{genusng1} implies $\bar \sH'_{k+1} \to \bar \sH'_k$ doesn't ramify. From Princ.~\ref{FP1} this says:   
\begin{triv} \label{pcuspex} For no value of $k$ does $\bar \sH_k'$ have a $p$ cusp. \end{triv} 
Now assume, contrary to \eqref{pcuspex},  $\bp_k'\in \sH_k'$ is a $p$ cusp for some $k$. Denote the degree of $\sH'_{k+1}/\sH_k'$ by  
$\nu_k$ and the number of points $\bp_{k+1}'\in \sH_{k+1}'$ over $\bp_k'$ by  $u_k$. 
Thm.~\ref{wmcneedspcusps} says possibilities for \eql{genusng}{genusng0} that  \cite{FrWProof} must eliminate are these. For $k>>0$,
$\nu_k=p$, 
$u_k=1$ and $\bar \sH'_{k+1}/\bar \sH'_k$ is equivalent (as a cover over $K$) to either:  
\begin{triv} \label{pcusponea}  a degree
$p$ {\sl polynomial map}; or \end{triv} 
\begin{triv} \label{pcusponeb} a degree $p$ rational function ramified precisely over two $K$ conjugate
points. 
 \end{triv} 

\begin{thm} \label{wmcneedspcusps} If none of \eqref{pcuspex}, 
\eqref{pcusponea} or \eqref{pcusponeb} hold for the component branch $B'$, then
$B'$ satisfies the conclusion of Main Conj.~\ref{wconj}: High levels of $B'$ have no $K$ points. 

For $B'$ with \eql{fmst}{fmstb}  holding  (full elliptic ramification; including when $B'$ has fine reduced moduli \wsp 
\S\ref{fineredMod}) for $k>>0$, the Main Conj.~holds unless  \eqref{pcusponeb} holds.  
\end{thm} 

\begin{proof} Assume \eqref{pcuspex} doesn't hold and $g_k'=0$ for  large $k$. That is,       
\begin{equation} \label{ktok+1} 2(\deg(\bar
\sH_k'/\prP^1_j)-1)= 
\ind(\gamma_{0,k}')+\ind(\gamma_{1,k}')+\ind(\gamma_{\infty,k}'):\  \eql{genusng}{genusng0} \text{ holds}.\end{equation} 
Consider now what would allow  $g_{k+u}'$, $u\ge 0$ to also be 0. 

Denote the cardinality of the $p$ cusps on $\sH_k'$ by $t_k$. For each
 $p$ cusp, $\bp_k'\in \sH_k'$,  Princ.~\ref{FP1} says the following.  
\begin{triv} \label{FP1pca} Each $\bp_{k+1}'$ over $\bp_k'$   has index 
$p$ and $\bar \sH_{k+1}'/\bar \sH_k'$ has degree $\nu_k=p\cdot u_k$. \end{triv} 
\begin{triv}\label{FP1pcb}  Also,  $t_{k+1}\ge t_k\cdot u_k$. 
\end{triv}

Apply \eqref{FP1pca}, by replacing $k$ by $k+1$,  to any $\bp_{k+2}'\in \sH'_{k+2}$ over a $\bp_{k+1}'$. Conclude: 
\begin{itemize} \item   there is an index
contribution of $t_k\cdot u_k\cdot u_{k+1}\cdot (p-1)$ from all $\bp_{k+2}'\,$s to Riemann-Hurwitz from
$\bar \sH'_{k+2}$ to $\bar \sH_{k+1}'$; and 
\item   Riemann-Hurwitz applied to  $ \bar \sH_{k+2}'\to \bar \sH_{k+1}'$  contradicts \eqref{ktok+1} if $$t_k\cdot u_k\cdot
u_{k+1}\cdot (p-1)> 2(p\cdot u_{k+1} -1).$$ \end{itemize}  

Suppose $t_k\ge 2$. Then, we contradict \eqref{ktok+1} if  $(u_k-1)\cdot p \ge u_k$. This happens unless $u_k=1$ or $u_k=2=p$. In the latter case, with $t_k=2$, we would have $t_{k+1}=4$ from
\eqref{FP1pcb}. Then, putting $p=2$ you see a contradiction by shifting $k$ to $k+1$. So, the argument forces (with $t_k\ge 2$)  $u_k=1$, 
$t_k=2$, and no ramification outside these two cusps. Further, under these assumptions (and \eql{genusng}{genusng0}), \eql{fmst}{fmstb} must
hold for  $k>>0$.  

On the other hand, if $t_k=1$ for $k>>0$, then (with \eql{genusng}{genusng0}), \eqref{FP1pcb} forces $u_k=1$. That means
$\sH_{k+1}'/\sH_k$ is a cover of genus 0 curves of degree $p$ with one place totally ramified. This is equivalent to a cover represented by
a polynomial (see Prop.~\ref{ratPtsifg0or1}). 
\end{proof} 

\begin{res} \label{resRam} A branch $B'$ of $\sT_{G,\bfC,p}$ contradicts case \eql{genusng}{genusng0} if there is a $p$
cusp at level
$k$ and $\bar \sH_{k+u+1}'/\bar \sH_{k+u}'$ has degree  $\ge p+1$.  For $B'$ to contradict \eql{genusng}{genusng1}, we only  need
one $p$ cusp at a high level $k$: Princ.~\ref{FP1}  forces $\sH_{k+1}'/\sH_k'$ to  ramify. \end{res} 

\subsubsection{Why \eqref{pcuspex}, \eqref{pcusponea}, or \eqref{pcusponeb} would contradict Conj.~\ref{wconj}} Prop.~\ref{ratPtsifg0or1}
shows the exceptional cases in \S\ref{redMC} are serious.  
\begin{lem} \label{ellipinf}  For any projective genus 1 curve $X$ over a number field $K$, we can extend $K$ to assume $X(K)$ is an elliptic
curve with infinitely many points. \end{lem} 

\begin{proof} Extend $K$ to assume $X(K)\ne \emptyset$, and use one of those points as an origin to assume $X$ is an elliptic
curve. Now form $\mu_K: G_K\to \GL_2(\hat \bZ)$, the action of $G_K$ on all division points of $X$. Put $X$ in Weierstrass normal
form, so its affine version has the shape $$\{(x,y)\mid y^2=x^3-u_2x-u_3\}.$$ Next we show $X(L)$ cannot be finite for each number
field $L/K$. Suppose it is. 

First we show $\mu_K$ is an embedding of $G_K$. Suppose not. Let $\sigma\in G_K$ with $\mu_K(\sigma)=1$, but assuming $\sigma\ne 1$,
there is a finite extension $L/K$ on which $\sigma$ acts nontrivially. Take a primitive generator
$x_0$ for $L/K$ (that is, $L=K(x_0)$). Solve for $y_0$ so that $(x_0,y_0)\in X(L')$, with $L'=K(x_0,y_0)$. By assumption $(x_0,y_0)$ gives a
division point on $X$, and clearly $\sigma$ acts nontrivially on it. 

That gives that $\mu_K$ is an embedding. Yet, a
simple consequence of Hilbert's irreducibility Theorem is that there is a Galois $L/K$ with group $S_n$ for any large integer $n$. It
is an elementary group observation that $S_n$ for $n>5$ large cannot embed in $\GL_2(\bZ/N)$ for any value of $N$. This contradiction finishes
the proof. 
\end{proof}
  
\begin{prop}  \label{ratPtsifg0or1} A \MT\ for which either \eqref{pcuspex},  \eqref {pcusponea} or \eqref{pcusponeb} holds fails 
the conclusion of Conj.~\ref{wconj}.
\end{prop} 

\begin{proof}  Recall: 
We start with a component branch
$B'$ having definition field a number field $K'$. If $B'$ satisfies \eqref{pcuspex}, then Lem.~\ref{ellipinf} gives $k=k_0$, and $K$ with
$[K:K']<\infty$ and $|\bar
\sH_{k_0}'(K)|=\infty$. Now we have a tower of elliptic curves, all isogenous. Each therefore has infinitely many rational points. Only
finitely many of these can be cusps, and the rest will be rational points on $\sH_k'$, for each $k$. That shows, if \eqref{pcuspex} holds, we
do contradict Conj.~\ref{wconj}. 

Now consider \eqref{pcusponea}:  $\bar\sH_{k+1}'\to \bar\sH_k'$ is a degree $p$ cover (over $K$) of genus 0 curves with a
distinguished totally ramified point $\bp_k'\in \bar\sH_k'$. Then, both $\bp_k'$ and the unique point $\bp_{k+1}'$ over it are $K$ points. So,
again
$\sH_k'(K)$ is infinite and if \eqref{pcusponea} holds, then we contradict Conj.~\ref{wconj}.  

Finally, consider \eqref{pcusponeb}. Suppose $X\to Y$ is a $K$ map of genus 0 curves of degree $N$. Then, they both
define elements of order 2 in the Brauer-Severi group $H^2(G_K,\bar K^*)$. Denote these $[X]$ and $[Y]$. Then, $N\cdot [X]=[Y]$ (in additive
notation -- see the argument of \cite[Lem.~4.11]{BFr02} for example). In particular, if $N=p$ is odd, and $K$ is large enough that $X$ has a
rational point, then
$[X]=[Y]=0$ and both have infinitely many rational points. The case for $N=2$ is even easier for it is automatic that $2\cdot
[\bar\sH_{k+1}']=0$ ($=[\bar\sH_{k}']$). For this case we immediately have a tower of  degree 2 maps between $\prP^1\,$s. So,
finishing
\eqref{pcusponeb} reverts to the previous case.  
\end{proof} 

\subsection{What we need to complete the Conj.~\ref{genusState} proof} \label{cuspGrowthProbs}  The results of
\S\ref{cuspGrowth} show  the main point in finishing the Main Conjecture for $r=4$ is a $p$ cusp at some high level. Better yet, if the 
$\limsup$ of  $\deg(\bar
\sH_{k+u+1}'/\bar
\sH_{k+u}')$ is {\sl not\/} $p$, one such $p$
cusp \eqref{FP1pcb} guarantees the $p$ cusp count (at level $k$) is unbounded as $k\mapsto \infty$ . Prop.~\ref{unbpcusps} gives  examples
that show how to compute a (growing) lower bound to the $p$-cusp count with the levels. 

\subsubsection{Reducing to pure cusp branches} 
\S\ref{cuspResults} calls an infinite cusp branch
$B$ pure in cases \eql{branchtype}{branchtypea} and \eql{branchtype}{branchtypec} if these have no extraneous (finite) start strings  
of  g-$p'$ (possibly followed by a  string of o-$p'$) cusps. Continue that notation to define $B$ by a 
sequence of cusp sets $({}_k\bg)\text{Cu}_4\subset \ni(G_k,\bfC)^{\inn}$. We can assume $k$ is large. That allows starting at any desired level.
So we revert to where $B$ is one of the pure  infinite cusp branches
$B$ in
$\sC_{G,\bfC,p}$ with representatives $$\{{}_k\bg=(\row {{}_kg} 4)\in
\ni_k'\}_{k=0}^\infty.$$ Here $\ni_k'$ is the braid ($\bar M_4$ \wsp \S\ref{NielClDict}) orbit on $\ni(G,\bfC)^{\inn,\rd}$ of ${}_k\bg$. 
For all
$k\ge 0$, either:  
\begin{edesc} \label{cuspType} \item  \label{cuspTypea}  From Princ.~\ref{FP1}, $p|({}_k\bg)\mpr$; or   
\item  \label{cuspTypeb} From Princ.~\ref{FP2}, ${}_k\bg$ is a g-$p'$ rep.; or  
\item  \label{cuspTypec}  From  Princ.~\ref{FP3} (or Rem.~\ref{extendFP3}), ${}_k\bg$ is a Weigel cusp with 
$$s_{R_{H_{2,3}}}({}_k\bg)=1=s_{R_{H_{1,4}}}({}_k\bg).$$   
\end{edesc} 

\subsubsection{Using a g-$p'$ cusp branch to get $p$ cusps} \label{learngpcusp} \S\ref{BrFrattProp}
\cite[\S9]{BFr02} does many cases of \eql{cuspType}{cuspTypeb}, where $p=2$ and there is a g-$p'$ cusp that is the shift of an H-M rep. Here is
what we learned, by example,  about getting $p$ cusps from it. Our  example continues   
\S\ref{disappearop'}:   the $(A_5,\bfC_{3^4},p=2)$ \MT\ where level 0 had no 2 cusps. 

Prop.~\ref{cnf.3} applies with the $\Spin_5\to A_5$ cover to show both level 1
components have $p$ cusps (with $p=2$) \cite[Cor.~8.3]{BFr02}. The full analysis says the component,
$\sH_+(G_1(A_5),\bfC_{3^4})^{\inn,\rd}$, containing all the H-M cusps, has genus 12 and degree 16 over the unique component of
$\sH(A_5,\bfC_{3^4})^{\inn,\rd}$. It also has all the real (and so all the $\bQ$) points at level 1 \cite[\S 8.6]{BFr02}.  Further, all except
the shift of the H-M cusps are 2 cusps. The other component,
$\sH_-(G_1(A_5),\bfC_{3^4})^{\inn,\rd}$ is obstructed, so no full branch over it has $\tG 2(A_5)$ (the whole $2$-Frattini cover of $A_5$) as a
limit group. 

\begin{prop} \label{unbpcusps} The number of $p$ cusps at level $k$ in any H-M component branch over $\sH_+(A_5,\bfC_{3^4})^{\inn,\rd}$ is
unbounded in
$k$.
\end{prop} 

\begin{proof} The argument has this abstract idea. Let $B=\{\bp_k\}_{k=0}^\infty$ be a g-$p'$ cusp branch. Suppose for $k\ge k_0$ you can
braid $\bp_k$ to a $p$ cusp $\bp'_k$ with ramification index exactly divisible by $p$. Then, Princ.~\ref{FP1} allows, with $k=k_0+u$, 
inductively braiding 
$\bp_k$ to a sequence of cusps $\bp'_k(1),\dots, \bp_k'(u)$ with $\bp'_k(t)$ having ramification index exactly divisible by $p^t$, $u=1,\dots,
t$.  From their ramification indices over $j=\infty$, these give $u$ different $p$ cusps at level $k_0+u$. 

For $\ni(G_k(A_5),\bfC_{3^4})$ you can take $k_0=1$ and $\bp_k'$ is produced as the {\sl near\/} H-M rep.~associated to $\bp_k$
\cite[Prop.~6.8]{BFr02}. 
 \end{proof} 

\subsubsection{Limit groups and field of moduli examples} These examples show our progress in computing, and that the consequences are relevant
to the abstract results. 

\begin{prob} \label{H-compbr} What are the limit groups of full component branches (\S\ref{limNielsen}) over
$\sH_-(G_1(A_5),\bfC_{3^4})^{\inn,\rd}$?
\end{prob} 

\begin{exmp}[Continuing Prob.~\ref{H-compbr}] \label{G1A5} By contrast to examples in \S\ref{heisenberg} and \S\ref{Z3limgps}, we don't yet 
know the limit groups for $\sH_-(G_1(A_5),\bfC_{3^4})^{\inn,\rd}$. Example:  Each space $\sH(A_5,\bfC_{3^r})^{\inn,\rd}$, $r\ge
5$, has exactly two components $\sH_{\pm}(A_5,\bfC_{3^r})$ \cite[Thm.~1.3]{altGps}. 
Also,  $\sH_{+}(A_5,\bfC_{3^r})$ is a g-$2'$ component. So,  from
Princ.~\ref{FP2} it has $\tG 2(A_5)$ as a limit group. 

Further, $\sH_{-}(A_5,\bfC_{3^r})$ has a unique limit group, just $A_5$. This is because
the 1st Loewy layer (\S\ref{loewyLayers}) of $M_0(A_5)$ consists of just the Schur multiplier $\bZ/2$ of $A_5$ \cite[Cor.~5.7]{BFr02}. We know
the Schur multiplier of
$G_1(A_5)$ is just $\bZ/2$. Still, what if other $A_5$ irreducible modules appear in the first Loewy layer of the characteristic module $M_1$? 
Then, akin to Ex.~\ref{schurComps}, the braid orbit corresponding to $\sH_-(G_1(A_5),\bfC_{3^4})^{\inn,\rd}$ could have all limit groups
larger than  $G_1(A_5)$. \end{exmp}

\begin{prob} \cite[Thm.~1.3]{altGps} says $\sH(A_n,\bfC_{3^r})^{\inn,\rd}$, $r\ge n$, always has exactly two components, which we
can denote $\sH_{n,r,\pm}$. When $p=2$, $\sH_{n,r,+}$ always has ${}_2\tilde G(A_n)$ as one limit group. Further, 
the limit groups of $\sH_{n,r,-}$ never include ${}_2\tilde G(A_n)$. Still,  as in
Ex.~\ref{G1A5}, for which $(n,r)$ is $A_n$ a limit group? From  \cite[Obst.~Lem.~3.2]{FrK97} (as in Lem.~\ref{liftLem}), the result only
depends on
$n$: Whether there is another irreducible in the 1st Loewy layer of $M_0(A_n)$. \cite[Rem.~2.5]{FrK97} (based on \cite{Ben2}) shows there is a
Frattini cover of
$A_8$ that doesn't factor through $\Spin_8$. So, 
$A_8$ is never a  limit group of $\sH_{8,r,-}$. We know little about this for $n\not\in \{4,5,8,9\}$. 
\end{prob}  

Our next example shows how significant are the cusps
$\bp_k'$ in the  braid from
$\bp_k$ to $\bp_k'$ in the proof of Prop.~\ref{unbpcusps}. The topic shows how one \MT\ produces an infinite number of closely related
situations contrasting the field of moduli and the field of definition of covers corresponding to points on tower levels.

\begin{exmp}[Moduli field  versus definition field] \label{modvsdef} Recall the cusps $\bp_k'$ achieved from braiding from H-M
cusps  in the proof of Prop.~\ref{unbpcusps}. These and the H-M cusps are  are the only  real 
(coordinates in $\bR$) cusps  on the
$(A_5,\bfC_{3^4}, p=2)$
\MT\ at level $k> 0$.   Let $R_k\to G_k(A_5)$ be the representation cover antecedent (\S\ref{repComps}) to the Schur multiplier of $A_5$. 

Regard the branch as defined over $\bR$. Then, $\bR$ points over
any $1< j< \infty$ in the real component abutting to $\bp_k$ represent covers in $\ni(R_k,\bfC_{3^4})$ whose field of definition
is
$\bR$ equal to its field of moduli. By contrast, with similar words concluding \lq\lq real component abutting to
$\bp_k'$\rq\rq\ (not $\bp_k$) here the moduli field is  $\bR$, but it is not a definition field  \cite[Prop.~6.8]{BFr02}.
\end{exmp}

\subsection{Chances for a genera formula} \label{braid-to-p} Ques.~\ref{wildGuess} asks if  
a g-$p'$ cusp branch represented by $B=\{{}_k\bg \in \ni_k'\}_{k=0}^\infty$ (notation like that of Princ.~\ref{FP1}) can deliver an analytic
expression for genera  akin to that for a modular curve tower. Further, Prop.~\ref{unbpcusps} supports why we expect to be able to braid from a
g-$p'$ cusp at level 0, in numbers increasing with $k$, a collection of $p$ cusps resembling those on modular curve towers (as in \cite[Talk
1]{Ont05}). \S\ref{genfromgp'} lists the challenges for this. 
\S\ref{abelgen} suggests simplifying  to a, still valuable, abelianized version.  

\subsubsection{Challenging a genera formula} \label{genfromgp'} Our examples show Ques.~\ref{wildGuess} is difficult. 
  
\begin{edesc} \label{ch} \item \label{cha} Are there o-$p'$ cusps in the orbit of ${}_k\bg$?  
\item \label{chb} For $k>>0$ are there any $p$-cusps in the orbit of ${}_0\bg$. If so, given how many there are at level 0; how many will
there be at level $k$?  
\item \label{chc} Can we separate the braid orbit of ${}_k\bg$ from other braid orbits? \end{edesc} 

We comment on these challenges. Example of \eql{ch}{cha}:  Prop.~\ref{op'overgp'} gives a \MT\
with related pairs of g-$p'$ and o-$p'$ cusps, represented respectively by ${}_k\bg$ and ${}_k\bg'$,  at every level. Can
you braid between ${}_k\bg$ and ${}_k\bg'$? 

Here is an immediate case wherein we must distinguish between
\eql{ch}{cha} and \eql{ch}{chb}. If ${}_o\bg=(g_1,g_1^{-1},g_2,g_2^{-1})$ is an H-M rep., there are two
possibilities since $\lrang{g_1,g_2}=G_0$: Either this is a $p$ cusp or it is an o-$p'$ cusp. For the latter, we guess at 
high levels that either the only cusps above it are $p$ cusps. Princ.~\ref{FP3} presents this possibility (contrary to Conj.~\ref{noWeigB}): 
\begin{triv} \label{3-cycle} There is an infinite branch on the \MT, $(G_0,\bfC',p)$ with $\bfC'$ the  conjugacy classes of $g_1,g_2$ and 
$g_1g_2^{-1}$. \end{triv} Having such a branch  is equivalent to having the homomorphism $\psi': M_{\bg'}\to G_0$
defined by
$\bg'=(g_1,g_2,(g_1g_2)^{-1})$ extending to $\tilde \psi': M_{\bg'}\to \tG p$. \cite{FrWProof} notes a necessary condition   from
the genus of the 3 branch point cover $X\to
\prP^1_z$  representing $\psi'$. It must exceed the rank of $\ker(\tG p\to G_0)$. Apply \eqref{3-cycle} to the example of
\S\ref{nogp'}, with $\ni(A_5,\bfC_{\pm53})$.  The genus $g$ of the corresponding 
$X$ satisfies $$2(60+g-1)=2(60/5)\cdot 4+ (60/3)\cdot 2,$$ so $g=9$, while the rank of $\ker(\tG p\to G_0)$ is 4.  

Example of \eql{ch}{chc}: Thm.~\ref{highObst} gives examples with at least two components \wsp one H-M \wsp at each higher level
of a \MT. The cases we give {\sl replicate\/} (in the sense of antecedent Schur multipliers) a two (or more) component situation at level 1.
This regularity of behavior is what we {\sl expect\/} with g-$p'$ cusps. Yet, is it always like this? 

\subsubsection{Shimura-like levels and abelianized genera} \label{abelgen} A level $k$ \MT\ component,
$\sH_k'$,  has above it a tower one may compare with Shimura varieties. That 
goes
like this. 

Let $\ker_k=\ker(\tG p\to G_k)$ (\S\ref{gpStatements}). The sequence of spaces comes 
from
forming $\tG p/(\ker_k,\ker_k)=\tG {p,k}$. This gives a
$p$-Frattini extension of $G_k$ by the abelian group
$\ker_k/(\ker_k,\ker_k)=L_k$, as in the proof of Lem.~\ref{ppowermap}. The lift of $g\in G_k$ to $\tilde g\in \tG {p,k}$ 
gives an action of $g$ on $L_k$ by
the conjugation by $\tilde g$. 

Form the spaces
$\{\sH_{k,u}\}_{u\ge 0}^\infty$ corresponding to the Nielsen classes $\ni(\tG 
{p,k}/p^uL_k,\bfC)$, and denote by
$\{\sH_{k,u}'\}_{u\ge 0}^\infty$ those (abelianized) components over $\sH_k'=\sH_{k,0}'$. 

Let $R_k'\to G_k$ (resp.~$R_k\to G_k$) be  maximal among central, $p$-Frattini (resp.~exponent $p$ Frattini) 
extensions of 
$G_k$. Then, $\ker(R_k'\to G_k)$ (resp.~$\ker(R_k\to G_k)$) is
the maximal 
$p$ quotient (resp.~exponent $p$) of $G_k\,$s Schur  multiplier. Cor.~\ref{weigCor1} checks for an infinite branch above a given component 
by  inductively checking  Nielsen elements ${}_k\bg$ for $s_{R_k/G_k}({}_k\bg)=0$ at successive levels for  {\sl all\/} $k$.
\S\ref{exofWeigTest} has examples that require successive checks. 

\begin{thm} \label{wtestb} For $u>>0$, $\sH_{k,u}'$ is nonempty  if and only if $s_{R_k'/G_k}(O')=0$ (just one test).
\end{thm}   

Ques.~\ref{wildGuess} has this easier, yet very valuable,  
variant. 

\begin{prob}[Abelianized Tower Genera] \label{genuskk+1abel} Label the precise 
ingredients
needed to compute genera of the $\{\sH_{k,u}'\}_{u\ge 0}^\infty$ components. 
\end{prob}

\section{Strong Conjecture for $r=4$} \label{achSchurDomain} Our strong Main Conjecture \ref{MainConj} is an expectation that 
the best \MT s are akin to those of modular curves. \S\ref{initMC} shows how the \MT\ cusp language applies to modular curves. 
\S\ref{secondMC} strengthens that, noting cusp branches defined by g-$p'$ cusps and $p$ cusps generalize projective sequences of
modular curve cusps. Finally, \S\ref{3rdMC} starts a discussion (continued in the appendix) on a non-modular curve \MT\  whose low levels have
genus 0 and 1 components with worthy applications.

\subsection{Initial comparison of \MT s with modular curves} \label{initMC} Let $D_{p^{k+1}}$ be the dihedral group of order $2\cdot p^{k+1}$
with $p$ odd. 

\subsubsection{The strong Main Conjecture} 
\cite[Lect.~1]{Ont05} computes the genera of the modular curves $X_0(p^{k+1})$  and
$X_1(p^{k+1})$  as \MT\ levels. Example: $X_1(p^{k+1})$, defined by $\ni(D_{p^{k+1}},\bfC_{2^4})^{\inn,\rd}$ with $\C_2$ the involution
class, has these properties. 
\begin{edesc} \label{modCurve} \item  \label{modCurvea} There is one $\bar M_4$ orbit. 
\item  \label{modCurveb} We inductively compute all cusps at level $k$ using an H-M
rep.~ (width $p^{k+1}$), and the shift of H-M rep.~cusps are g-$p'$ cusps of width 1.  
\item  \label{modCurvec} $\gamma_0'$ or $\gamma_1'$ have no fixed points.
\item  \label{modCurved} $\sQ''$ (\S\ref{NielClDict}) acts trivially at all levels. 
\end{edesc}

\noindent \cite[Prop.~8.4]{FrS05} generalizes \eql{modCurve}{modCurvec} and \eql{modCurve}{modCurved}. This is the \MT\  version of Serre's
abelian variety lemma: (roughly) among automorphisms, only the identity fixes many torsion points. Use the notation of \S\ref{setUpSM} for a
\MT\ of rank $u\ge 0$. Again, assume  $r=4$ for these \MT s.

\begin{guess}[Strong Main Conjecture]   \label{MainConj}
$P_\bfC$ Version: Over all  $p\not\in
P_\bfC$, for only finitely many  
$V\in \sV_p(J)$, does $\sH(V\xs
J,\bfC)^{\inn,\rd}$ have genus 0 or 1 components. \end{guess} 

\noindent There is a $P_\bfC'$ version, though the weak
Conjecture and Conj.~\ref{MainConj} imply it.  

\begin{guess}[Mazur-Merel Version of the strong Main Conjecture] \label{MMMainConj} With hypotheses of Conj.~\ref{MainConj}, over
all  $p\not\in P_\bfC$, for only finitely many  
$V\in \sV_p(J)$, does $\sH(V\xs
J,\bfC)^{\inn,\rd}$ have a rational point. \end{guess} 

\subsubsection{Comparison with the Strong Torsion Conjecture} \label{compSTC} The following observation generalizing \cite[Thm.~6.1]{BFr02}
appears in
\cite[Prop.~4.10]{cadoret05b}. For the 
\MT\ of 
$(G,\bfC,p)$, the (weak) Main Conjecture {\sl for all values of $r$\/} follows from Conj.~\ref{STC}, called by \cite{Si92} and
\cite{Ka98} the Strong Torsion Conjecture. Let $K$ be a number field. 

\begin{guess}[STC] \label{STC} For $g,d\ge 1$,  there exists $n(g,d)\ge 1$ with this property. If $n \ge n(d, g)$, then there are no dimension
$g$ abelian varieties
$A$ defined over  $K$, with  $[K : Q] \le d$, and having a $K$ torsion point of order $n$.
 \end{guess} 

By contrast, the weak conjecture for the \MT\ given by $(D_{p}, \bfC_{2^r}, p)$ (necessarily for a nonempty \MT, $r=2g+2\ge 4$ is even) is
{\sl equivalent\/} to the following. 

\begin{guess} For $k$ large there is no cyclic group $C\cong \bZ/p^{k+1}$ of torsion on a hyperelliptic Jacobian of genus $g$
for which $G_K$ acts on $C$ through its cyclotomic action on $\lrang{e^{2\pi i/p^{k+1}}}$ 
\cite[\S5.2]{DFr94}. 
\end{guess} Further, the Strong Main Conjecture for a higher rank \MT\ doesn't follow from Conj.~\ref{STC} because the genus of the curves (and
so the dimension of the Jacobians) in question grows with primes $p$. 
  
\subsection{Modular curve comparison for Serre's OIT}   \label{secondMC} Principles \ref{FP2} and \ref{FP3} help toward describing all 
branches in $\sC_{G,\bfC,p}$. This guides the strong Conjecture in how it might effectively generalize
Serre's Open Image Theorem (OIT) \cite{SeAbell-adic}. 

\subsubsection{Frattini properties in the OIT} Here  are significant OIT ingredients. 

\begin{edesc} \label{serreProp} \item Acting by $G_{\bQ_p}$ on projective systems of points in neighborhoods of  H-M reps. on
$\{X_1(p^{k+1})\}_{k=0}^\infty$ gives a transvection in the projective sequence of monodromy inertia groups. 
\item \label{serrePropb} The geometric monodromy group, $\PSL_2(\bZ/p^{k+1})$, for  
$X_1(p^{k+1})\to \prP^1_j$  is a $p$-Frattini cover of the monodromy at level 0 if $p\ne 2$ or 3. 
\end{edesc} 

Here is how \eql{serreProp}{serrePropb} works ($p$ is odd). Let  $\{\bp_k\in
X_0(p^{k+1})\}_{k=0}^\infty$ be a projective sequence 
of points over $j'\in F$. Then $G_F$ acts on these to give a map   
$$G_F\mapright{\psi_{2,j'}}\lim_{\infty\leftarrow k}
\GL_2(\bZ/p^{k+1})/\{\pm I_2\}=GL_2(\bZ_p)/\{\pm I_2\}\mapright{\scriptscriptstyle{\text{ Det}}} \GL_1(\bZ_p).$$ 

The induced map $\psi_{1,j'}: G_F\to \GL_1(\bZ_p)$ is onto an open subgroup 
because (essentially) all the roots of 1 are present in the field generated by the division points on elliptic curves. This deduction
interprets  from the Weil pairing on elliptic curves. This is an alternating pairing on $p^{k+1}$ division points into $p^{k+1}$th roots
of one
\wsp interpreted as the cup product pairing from 1st ($\ell$-adic, but $\ell=p$) cohomology to the 2nd $\ell$-adic cohomology.
Rem.~\ref{SeOIT} states the \MT\ version of this. 

Let $G_F^0$ be the kernel of $\psi_{1,j'}$. Consider the restriction $\psi_{2,j'}^0: G_F^0 \to \PSL_2(\bZ_p)$, and composite by going  
$\mod p$  to get 
$\psi_{2,j'}^0 \!\!\mod p: G_F^0 \to \PSL_2(\bZ/p)$. 

\begin{res} \label{SeFrat} For  $p\ne 2$ or 3, if $\psi_{2,j'}^0 \!\!\mod p$ is  onto, then $\psi_{2,j'}^0: G_F^0\to \PSL_2(\bZ_p)$  is
onto.  If $p= 3$ (also for $p=2$), and $\psi_{2,j'}^0 \!\!\mod p^2$ is onto, then so is  $\psi_{2,j'}^0$. \end{res} 

\begin{proof}[Comments] First: The $\!\!\mod p$ map $\PSL_2(\bZ_p)\to \PSL_2(\bZ/p)$ is a Frattini cover  if
$p\ne 2$ or 3 \cite[IV-23 Lem.~2]{SeAbell-adic}. It isn't, however, the universal $p$-Frattini cover of $\PSL_2(\bZ/p)$, ever! For example, 
consider the case $p=5$:
$\PSL_2(\bZ/p)=A_5$. Then,
$M_0=\ker(G_{5,1}(A_5)\to A_5)$ (notation of \S\ref{gpStatements}) is a rank 6, $A_5$ module. It fits in a nonsplit short exact sequence
$0\to M'\to M_0\to M'\to 0$ with  $M'$ the adjoint representation of $\PSL_2(\bZ/5)$ (on $2\times 2$ trace 0 matrices 
\cite[Rem.~2.10]{Fr4}).

For $p=3$,  $\PSL_2(\bZ/3)$ is not simple. Yet,  $\PSL_2(\bZ_3)\to \PSL_2(\bZ/3^2)$ is Frattini. 

\cite[Thm.~C]{Weigel} computes the rank of $\ker(G_{p,1}(\PSL_2(\bF_q))\to \PSL_2(\bF_q))$ when $\bF_q$ is the finite field of order $q=p^u$.
The adjoint representation appears a lot.  This also computes this characteristic rank for the other primes dividing $|\PSL_2(\bF_q)|$, giving 
important empirical data for effective computation of Frattini ranks. 

Let $R_q$ be the Witt vectors for $\bF_q$. \cite[\S4]{V3} notes that $\GL_n(R_q)\to \GL_n(\bF_q)$ is a Frattini cover so long as $p>2$ does not
divide $n$, and if $p=3$, $n\ge 4$.  \cite[\S4]{Va03} uses this Frattini principle in the full context of Shimura varieties, continuing
the tradition of \cite{{SeAbell-adic}}. Those with  Shimura variety experience know that the semi-simple groups that arise, generalizing the 
$\PSL_2$ case (symplectic groups, for example), are from a moduli problem on abelian varieties. 
\end{proof} 

\begin{rem} It is elementary that $\psi_{2,j'}^0 \!\!\mod p$ (in Res.~\ref{SeFrat}) is onto  for a dense set $j'$ in any number
field. For $p\ne 2$ or 3, just apply Hilbert's Irreducibility Theorem to the irreducible cover $X_0(p)\to \prP^1_j$ (for $p=3$, to
$X_0(p^2)\to
\prP^1_j$).
\end{rem} 
 
\subsubsection{F(rattini)-quotients of \MT s} \label{F-q}  Consider a rank $u$ \MT\ from $F_u\xs J$ and 4 conjugacy 
classes in $J$ (\S\ref{setUpSM}). For $p\not\in P_{\bfC}$, assume $\tilde G^*=V^*\xs J\in \sG_{J,p}$ is a $\bfC$ $p$-Nielsen
limit. That means there are projective systems of $\{\bg_V\in \ni(V\xs J,\bfC)\}'$ with $'$ indicating running over finite 
$J$ quotients of $V^*$ covering $\bZ/p^u$. This projective system defines a cusp branch. 

By taking braid orbits,  these define a projective system of  \MT\ 
components on the full component graph
$\sT^f_{\bZ/p^u\xs J,\bfC,p}$. Use our previous notation $B$ for a cusp branch and $B'$ for the component branch $B$ defines. For a $J$ quotient
$V$ of
$V^*$ use $B_V$ and $B_V'$ for the corresponding cusp $\bg_V$ and its component. Let $F_\bfC$ be the definition field of all the inner
reduced Hurwitz spaces $\sH(G_k((\bZ/p)^u)\xs J,\bfC)^{\inn,\rd}$ as in \S\ref{cuspBranches}. To simplify, assume $F_\bfC=\bQ$. 
\newcommand{\bG}{\text{\bf G}}
\begin{defn} Suppose $V_0$ is a $J$ quotient of $(\bZ/p)^u$. We call the \MT\ for $(V_0\xs J,\bfC,p)$ an F-quotient of the \MT\ for
$((\bZ/p)^u\xs J,\bfC,p)$. Then, there is a natural map from
$\sT^f_{\bZ/p^u\xs J,\bfC,p}$ to
$\sT^f_{V_0\xs J,\bfC,p}$ (on cusps also) induced by the map $$\sH((\bZ/p)^u\xs J,\bfC)^{\inn,\rd}\to \sH(V_0\xs J,\bfC)^{\inn,\rd}\eqdef
\sH_{V_0}.$$  
\end{defn} We will refer to  $B_V'$ on branch $B'$ as if it is the corresponding Hurwitz space. Also, for $V_1\to V_2$
a homomorphism of $J$ groups, denote the corresponding Hurwitz space map as $B'_{V_1}\to B'_{V_2}$. Let $\bG_V$ be the geometric
monodromy group of
$B'_V\to
\prP^1_j$. 

In the best circumstances for the cusp branch $B$, as in \S\ref{appMTRank}, we expect this.   
\begin{edesc}  \label{serreAnal} \item \label{serreAnala} Computable $\bQ_p$ action: We can decipher the $G_{\bQ_p}$ orbit on $B$. 
\item \label{serreAnalb} Branch Frattini Property: Excluding finitely many  $V_2$ corresponding to $B_{V_2}'$ on the branch $B'$, all the maps
$\bG_{V_1}\to
\bG_{V_2}$ are
$p$-Frattini covers.  
\item  Smooth genera: The genera of $B'_V$
should have a modular curve-like  formula, coming from clear understanding of g-$p'$ and $p$-cusps on $B'$.   
\end{edesc}   

\S\ref{comp-tang} notes results on the $\bR$ and $\bQ_\ell$ nature of cusp branches, extending \eql{serreAnal}{serreAnala}. 

\subsubsection{More on Branch Frattini propery \eql{serreAnal}{serreAnalb}} \label{BrFrattProp} A weaker version of
\eql{serreAnal}{serreAnalb} would assert that $\bG_{V_1}\to
\bG_{V_2}$ is a $p$-group. In turn this implies all ramification groups are $p$-groups, and  Lem.~\ref{fineredBr} (condition
\eql{fmst}{fmstb}) implies exactly that. 
 
Propery \eql{serreAnal}{serreAnalb} is an analog of Serre's use of the $p$-Frattini property. We expect something like it for all reasonable
\MT s. For example, suppose we have a g-$p'$ (or even, shift of an H-M) cusp on a \MT. Then, we expect the geometric
monodromy groups
$\bG_k$ of $\bar
\sH(G_k(G),\bfC)\to \prP^1_j$ to satisfy \eql{serreAnal}{serreAnalb}. 

That is, for $k_0$ large and $k\ge k_0$, $\bG_k\to \bG_{k_0}$ should be 
a $p$-Frattini cover. For certain, however, we can't always take $k_0=0$.  For example, for the  
\MT\ for $(A_5,\bfC_{3^4},p=2)$ we have these facts. This continues   Ex.~\ref{highObstExs}, 
Ex.~\ref{level1A5}, \S\ref{disappearop'},\S\ref{learngpcusp}, Ex.~\ref{modvsdef} and  \S\ref{BrFrattProp}.

\begin{edesc} \item There is exactly one H-M component $B_1'$ at level 1. 
\item  the degree of $B_1'\to B_0'$ is 16, but 
\item  $|H_{1,0}|\!=|\ker(\bG_1\to \bG_0)|\!=3\cdot 2^6$ with an 
$S_3$ at the top \!\cite[App.~A]{BFr02}. \end{edesc} So, $H_{1,0}$ is not even a two group. We use proofs, not  {\bf GAP}
calculations, so we know why this is happening. Prob.~\ref{fratMT} starts with a fixed g-$p'$ branch (as in
\S\ref{F2Z3}). 

\begin{prob} \label{fratMT}  Show $H_{k+1,k}=\ker(\bG_{k+1}\to \bG_k)$ is a 2-group (resp.~$p$-group) for large $k$ for the
$(A_5,\bfC_{3^4},p=2)$ (resp.~$((\bZ/p)^2\xs
\bZ/3,\bfC_{\pm 3^2}, p\ne 3)$
\MT. \end{prob} 

My thinking \eql{serreAnal}{serreAnalb} might hold came from \cite{Iha86} (even though Ihara has $p$-groups, the opposite of $p$-perfect
groups). 

Of course, if we knew explicitly the  subgroups of $\PSL_2(\bZ)$ defining the \MT\ levels that would answer Prob.~\ref{fratMT}. Even one other
case than modular curves where we could test these problems would be reassuring. In fact, \cite{Berger} almost includes  the non-trivial
F-quotient of
$((\bZ/p)^2\xs
\bZ/3,\bfC_{\pm 3^2}, p\equiv 1 \mod 3)$. Only, he has taken for $\bfC$ the repetition 3 times of one conjugacy class,
and the other just once? He uses the Bureau representation of the braid group to effect his calculation. It promises 
answering such questions as Prob.~\ref{NewSeOIT} for at least this non-modular curve situation. 

\subsubsection{Complete fields and tangential base points} \label{comp-tang} Suppose $B$ is a cusp branch. Much work on the Inverse Galois
Problem is appropriate for service to this problem. 

\begin{prob} \label{ellSeq} What do we need to know to detect when $B$ is a projective sequence of $\bQ_\ell$ cusps, $\ell\ne p$ (including 
$\ell=\infty$)?
\end{prob} 

The  effective
computation for $\bR$ points on Hurwitz spaces in \cite{DFr90}  works to analyze higher \MT\ levels (as in \cite[\S6]{BFr02}, especially see the
use made in Ex.~\ref{modvsdef}). The model for $\bQ_\ell$ has followed this. It is necessary  for a positive answer to Prob.~\ref{ellSeq} that
the manifolds
$\bar \sH_k'$ have definition field $\bQ_\ell$. 

The basic proposition in that direction is \cite[Thm.~3.21]{Fr4}. It says: If all H-M reps.~in the Nielsen classes for
level $k$ lie in one braid orbit (so all the H-M cusps lie on $\bar \sH'_k$) then $\bar \sH_k'$ has definition field $\bQ$. Further, it gives
a criterion for this to happen at level 0 that implies it automatically at all other levels. Then, Harbater patching applies to 
produces a projective sequence of
$\bQ_\ell$ cusps on 
$\{\bar
\sH_k'\}_{k=0}^\infty$. 
\cite[Thm.~2.7]{pierre} has a precise statement from \cite{DDe04}. 

\cite{DEm04} redoes the author's result 
using a more classical compactification. One problem: When $r=4$, the criterion of  \cite[Thm.~3.21]{Fr4} never applies.
An example failure is the two H-M components  at Level 1 in
\S\ref{level1p2} (see Rem.~\ref{conjComp}). 

So, we require deeper methods to analyze the definition field of a component branch 
and its  cusps when $r=4$. Based on \cite{IM95} and
\cite{WeFM}, \cite[App.~D.3]{BFr02} describes a method that {\sl will\/} work with sufficient grasp of the group theory and use  
of an especially good cusp branch. 

Again, $B$ is a g-$p'$ cusp branch, defining a component branch $B'$ on a \MT. The desired archetype for a tangential base
point comes from  $X_0(p^{k+1})$. We identify this space with $\sH(\bZ/p^{k+1}\xs \,\bZ/2,\bfC_{2^4})^{\abs,\rd}$; the {\sl absolute\/} 
reduced Hurwitz space related to the nontrivial F-quotient in Serre's OIT. The unique cusp of width $p^{k+1}$ identifies with the unique H-M
cusp, and so it has  $\bQ$ as definition field. 

In the now classical picture, points on the space approaching this cusp preciously go to a
controlled $p$-catastrophe. A $p$-adic power series  representing $j$, parametrizes a Tate curve  ($p$-adic torus) degenerating with  $j\mapsto
\infty$ ($p$-adically). 

Generalizing such constructions to
g-$p'$ cusps cannot be trivial. Yet, the apparatus for exploiting them as Serre does in \cite[IV.29--IV.45]{SeAbell-adic} is already
in the Grothendieck-Teichm\"uller motivated formulas of Ihara-Matsumoto-Wewers (\cite{IM95}, \cite{WeFM}; \cite[App.~D]{BFr02} discusses this). 
Making it work, a' la 
\cite{Nakamura}, in our more general situation requires a dedicated project. Deciding the definition field
of the two genus 1 components  in  \eql{apps}{appsb} is a practical example of its value. 

The groups $H_{2,3}(\bg)$ and $H_{1,4}(\bg)$ give a {\sl type\/} to g-$p'$ cusps. \cite[Lect.~4]{Ont05} defines 
g-$p'$  rep.~types in  Nielsen classes for any $r$, making sense of Prob.~\ref{HMgp'} for all $r$. 

\begin{prob} \label{HMgp'} Show this analog of \cite[Thm.~3.21]{Fr4} for general g-$p'$ cusp branches of a given type holds. If there
are finitely many (resp.~one) braid orbit of this type, then $G_F$ has a finite orbit (resp.~is fixed) on their component branch(s).  
\end{prob}

\begin{rem} \label{conjComp} Examples show that outer automorphisms of $G_k$
can  conjugate distinct H-M components on $\sH(G_k,\bfC)$ (\eqref{apps} and \cite[\S9.1]{BFr02}). Is this  is a general phenomenon? Nor do we
know if there are {\sl always}, modulo braiding, just finitely many
$G_F$ orbits of (shifts of) H-M reps. This makes sense for all g-$p'$ cusps. 
\end{rem} 

\subsection{$F_2\xs \bZ/3$, $p=2$: Level 0, 1 components} \label{3rdMC}   Components on these levels bring up 
deeper aspects of complex multiplication and the inverse Galois problem. This example  shows how  such tools as 
the {\sl $\sh$-incidence matrix\/} can identify components at a \MT\ level. We now explain why at level 0 there are two components:
$$\sH(\tilde F_{2,2}/\Phi^1\xs J_3,\bfC_{\pm 3^2})^{\inn,\rd}=\sH_0^+\cup \sH_0^-.$$ Both have genus 0, and $\sH_0^+$ is an H-M component.  The
other has nontrivial lifting invariant;  there is nothing above it at level 1. Though both are families of genus 1 curves, and upper half plane
quotients, neither is a modular curve. 

\subsubsection{Setting up reduced Nielsen classes} This Nielsen class has
$G=A_4$ with 
$\bfC_{\pm 3^2}$ as two pairs of 3-cycles in each of the conjugacy classes with order 3. First look at the situation with $A_3$ replacing
$A_4$. 

The total Nielsen class $\ni(A_3,C_{\pm 3^2})^\inn$ contains six 
elements corresponding to the six possible arrangements of the conjugacy classes. Since $A_3$ is abelian,
the inner classes are the same. Also, the outer automorphism of $A_n$ ($n=3$ or 4) from conjugation by
$(1\,2)\in S_n$ restricts to $A_3$ to send a conjugacy class arrangement to its
complement. Here is a convenient list of the arrangements, and their complements: 
$$\begin{array}{rrr} [1] +-+-& [2] ++-- & [3] +--+ \\ {[4]} -+-+ & [5] --++ & [6] -++-.\end{array}$$
The group $\sQ''=\lrang{q_1q_3^{-1}, \sh^2}$ equates elements in this list with their complements. So, inner reduced classes and 
absolute (not reduced) classes are the same. Conclude: $\sH(A_3,\bfC_{\pm 3^2})^{\inn,\rd}\to \prP^1_j$
is a degree three cover with branch cycles
$$(\gamma_0^*,\gamma_1^*,\gamma_\infty^*)=((1\,3\,2),(2\,3),(1\,2)).$$  Check easily: If $(\row
g 4)$ maps to [1], and (with no loss) $g_1=(1\,2\,3)$, then either this is $\bg_{1,1}$ (in \eqref{g11}) or $g_1g_2$ has
order 2. Listing the four order 2 elements gives a total of five elements in the reduced
Nielsen class $\ni(A_4,\bfC_{\pm 3^2})^{\inn,\rd}$ lying over [1]. 

\subsubsection{Effect of $\gamma_\infty$ on $\ni(A_4,\bfC_{\pm 3^2})^{\inn,\rd}$} \label{gammaA4action} 
Start with an H-M rep over [1] in
$A_3$:
\begin{equation} \label{g11} \bg_{1,1}=((1\,2\,3), (1\,3\,2),(1\,3\,4), (1\,4\,3))\in\ni(A_4,\bfC_{\pm 3^2}).\end{equation}  The middle twist 
squared on this conjugates the middle two by $(1\,4)(2\,3)$ to give $$\bg_{1,2}=((1\,2\,3),
(4\,2\,3),(4\,2\,1), (1\,4\,3)).$$  The result is a $\gamma_\infty$ orbit of length 4. The middle twist 
squared on 
$$\bg_{1,3}=((1\,2\,3),(1\,2\,4),(1\,4\,2),(1\,3\,2))$$ leaves it fixed, giving a $\gamma_\infty$ orbit of
length 2. Similarly, the square of the middle twist on
$\bg_{1,4}=((1\,2\,3),(1\,2\,4),(1\,2\,3),(1\,2\,4))$ conjugates the middle pair by
$(1\,3)(2\,4)$ producing
$\bg_{1,5}=((1\,2\,3),(1\,2\,4),(2\,4\,3),(1\,4\,3))$. Again the middle twist gives an element of order
4 on reduced Nielsen classes. 

The H-M rep.~$\bg_{3,1}=((1\,2\,3), (1\,3\,2),(1\,4\,3), (1\,3\,4))\in\ni(A_4,\bfC_{\pm 3^2})$ maps to [3] in
$A_3$. Applying $\gamma_\infty$  gives 
$\bg_{3,2}=((1\,2\,3),(1\,2\,4 ) ,(1\,3\,2), (1\,3\,4))$, the same as  conjugating on the
middle two by $(2\,4\,3)$.  The result is a length 3 $\gamma_\infty$ orbit.

On Nielsen class representatives over [3], $\gamma_\infty$ has one
orbit of length 3 and two of length one. See this by listing the second and third positions (leaving
$(1\,2\,3)$ as the first). Label these as $$\begin{array}{rl}&1'=((1\,3\,2),(1\,4\, 3)), 
2'=((1\,2\,4),(1\,3\, 2)),3'=((1\,2\,4),(2\,3\, 4)),\\ &4'=((1\,2\,4),(1\,2\, 4)),  5'=((1\,2\,4),(1\,4\,
3)).\end{array}$$ 

\subsubsection{Using Wohlfahrt's Theorem} \label{WohlfahrtThm} 
For $ \Phi^\rd:  \sH^\rd\to U_\infty$, one of our reduced Hurwitz space covers, let $\Gamma  \le \SL_2(\bZ)$ define it as an upper half-plane 
quotient $\bH/\Gamma$ (\S\ref{jinv}). Now let  $N_\Gamma$ be the least common
multiple (lcm) of its cusp widths. Equivalently:  $N_\Gamma$ is the lcm of the ramification
orders of points of the compactification 
$\bar \sH^\rd$ over
$j=\infty$; or the lcm of the orders of $\gamma_\infty$ on reduced Nielsen classes.  

Wohlfahrt's Theorem
\cite{Wohlfahrt} says $\Gamma$ is congruence if and only if $\Gamma$ contains the congruence subgroup, $\Gamma(N_\Gamma)$, defined by
$N_\Gamma$.  We have a situation with a modular curve-like aspect, though we find these $j$-line covers aren't 
modular curves by seeing the cusps fail Wohlfahrt's condition. Here is our procedure. 

Compute $\gamma_\infty$ orbits on $\ni^\rd$. Then, check their distribution among $\bar
M_4=\lrang{\gamma_\infty,\sh}$ orbits ($\sH^\rd$ components). For each $\sH^\rd$ component
$\sH'$,  check  the lcm of $\gamma_\infty$ orbit lengths to compute $N'$, the modulus if it were a modular curve. Then, see whether a
permutation representation of $\Gamma(N')$ could produce 
$\Phi': \sH'\to \prP^1_j$, and the type of cusps now computed.  Denote $\Spin_4$ (\S\ref{pperfect}) by $\hat A_4$.  

Use notation  ending
 \S\ref{gammaA4action}. Note: Neither of $\sH_0^{\inn,\rd,\pm}$ have reduced fine moduli. The Nielsen braid orbit for
$\sH_0^{\inn,\rd,-}$ (resp.~$\sH_0^{\inn,\rd,+}$) fails \eql{failfm}{failfma} (resp.~and also \eql{failfm}{failfmb}): 
\begin{edesc} \label{failfm} \item  \label{failfma} $\sQ''$ has length 2 (not 4 as required in \eql{fmst}{fmsta}) orbits; and  
\item  \label{failfmb} $\gamma_1$ has a fixed point (Lem.~\ref{g0g1fix}; contrary to
\eql{fmst}{fmstb}). \end{edesc}

\begin{prop} \label{A4L0} Then, $\gamma_\infty$ fixes $4'$ and $5'$ and cycles $1'\to 2'\to
3'$. So there are two
$\bar M_4$ orbits on $\ni(A_4,\bfC_{\pm 3^2})^{\inn,\rd}$, $\ni^+_0$ and $\ni^-_0$, having respective
degrees 9 and 6 and respective lifting invariants to $\hat A_4$ of $+1$ and $-1$. The first, containing
all H-M reps., has orbit widths 2,4 and 3. The  second has orbit widths 1,1 and 4. Neither defines a
modular curve cover of
$\prP^1_j$. 

Denote the corresponding completed covers 
$\bar \psi^{\pm}_0: \bar \sH^{\inn,\rd,\pm}_0\to \prP^1_j$. Both $\bar \sH^{\inn,\rd,\pm}_0$ have genus 0. Both have
natural covers 
$\bar
\mu^{\pm}:
\bar \sH^{\inn,\pm}_0\to \prP^1_j$ by completing the map \begin{equation} \label{ja} \bp\in \sH^{\inn,\rd,\pm}_0\mapsto
\beta(\bp)\eqdef j(\Pic(X_\bp)^{(0)})\in \prP^1_j.\end{equation} Then, this case's identification of inner and absolute reduced classes gives 
\begin{equation} \label{jb} \bp\in
\sH^{\inn,\rd,\pm}_0\mapsto (j(\bp),j(\Pic(X_\bp)^{(0)})),\end{equation} a birational embedding of
$\bar \sH^{\inn,\rd, \pm}_0$ in
$\prP^1_j\times \prP^1_j$. 

If we denote the corresponding $H_4$ orbits on $\ni(A_4,\bfC_{\pm3^2})^\inn$  by $\ni^{\inn,\pm}$,
then $\sQ''$ orbits on both have length 2.\end{prop}

\subsection{Proof of Prop.~\ref{A4L0}} This proof takes up the next four subsections. 

\subsubsection{$\gamma_\infty$ orbits on $\ni(A_4,\bfC_{\pm 3^2})^{\inn,\rd}$} \label{A4lev0} First: 
$\gamma_\infty$ fixes
$4'$ and it maps
$5'$  to
$((1\,2\,3),(2\, 3\, 4),(1\,2\,4),(3\,1\,2)) \text{ (conjugate by $(1\,2\,3)$ to $5'$)}$. 
  
These computations establish the orbit lengths: 
$$\begin{array}{lr}  &(g_{1,1})\gamma_\infty=((1\,2\,3),(1\,4\,2),(1\,3\,2),(1\,4\,3))=(3')\sh, \\
&(g_{1,3})\gamma_\infty = ((1\,2\,3),(1\,4\,2), (1\,2\,4),(1\,3\,2))=(1')\sh. \end{array}$$ They put the
H-M rep. in the $\bar M_4$ orbit with $\gamma_\infty$ orbits of length 2,3 and 4 (in the orbit of
the
$1'\to 2'\to 3'$ cycle). Use $\ni^+_0$ for the Nielsen reps.~in this
$\bar M_4$ orbit. 

\subsubsection{Graphics and Computational Tools: $\sh$-incidence} \label{graph-comp}  
The $\sh$-incidence matrix of $\ni^{+}_0$ comes from the following data. Elements
$\bg_{1,1},\bg_{1,2},\bg_{1,3}$ over [1] are permuted as a set by $\sh$. They map by $\gamma_\infty$
respectively to
$\bg_{2,1},\bg_{2,2},\bg_{2,3}$ over [2]. Under $\gamma_\infty$ these map respectively to $\bg_{1,2},\bg_{1,1},\bg_{1,3}$,
while $\bg_{3,1}, \bg_{3,2},
\bg_{3,3}$ cycle among each other. So, there are three $\gamma_\infty$ orbits,
$O_{1,1}$,
$O_{1,3}$ and $O_{3,1}$ on $\ni^{+}_0$ named for the subscripts of a representing element. 

The data above shows 
$$|O_{1,1}\cap (O_{3,1})\sh|=2,\ |O_{1,3}\cap (O_{3,1})\sh|=1.$$ Compute:  $\sh$ applied to $\bg_{1,3}$ is
$\bg_{1,1}$ so 
$|O_{1,1}\cap (O_{1,3})\sh|=1$.  The rest has two sources: 
\begin{itemize} \item symmetry of the $\sh$-incidence matrix, and; 
\item  elements in a row (or
column) add up to ramification index of the cusp labeling that row (or column).  \end{itemize}  

\begin{table}[h] \caption{$\sh$-Incidence Matrix for $\ni^{+}_0$} \label{shincni0+}
\hspace{1.6in} \begin{tabular}{|c|ccc|} \hline Orbit & $O_{1,1}$\ \vrule  &
$O_{1,3}$\ \vrule & 
$O_{3,1}$\  \\ \hline $O_{1,1}$ 
&1&1&2\\ 
$O_{1,3}$ &1 &0&1\  \\ $O_{3,1}$ &2&1&0\  \\   \hline \end{tabular} \end{table}

Similarly, the $\sh$-incidence matrix of $\ni^{-}_0$ comes from the following data. Elements $\bg_{1,4},
\bg_{1,5}$ over [1] map by $\gamma_\infty$ respectively to
$\bg_{2,4},\bg_{2,5}$ over [2], and these map respectively to $\bg_{1,5},\bg_{1,4}$, while $\gamma_\infty$
fixes both $\bg_{3,4},\bg_{3,5}$. So, there are three $\gamma_\infty$ orbits,
$O_{1,4}$,
$O_{3,4}$ and $O_{3,5}$ on $\ni^{-}_0$. 

\begin{table}[h] \caption{$\sh$-Incidence Matrix for $\ni^{-}_0$} \label{shincni0-}
\hspace{1.6in} \begin{tabular}{|c|ccc|} \hline Orbit & $O_{1,4}$\ \vrule  &
$O_{3,4}$\ \vrule & 
$O_{3,5}$\  \\ \hline $O_{1,4}$ 
&2&1&1\\ 
$O_{3,4}$ &1 &0&0\  \\ $O_{3,5}$ &1&0&0\  \\   \hline \end{tabular} \end{table}

\begin{lem} \label{g0g1fix} In general, the
$\sh$-incidence matrix is the same as the matrix obtained by replacing
$\sh=\gamma_1$ by $\gamma_0$. Further, the only possible elements fixed by either lie in $\gamma_\infty$
orbits
$O$ with $|O\cap (O)\sh\ne 0|$.  

On 
$\ni^{+}_0$ (resp.~$\ni^{-}_0$), 
$\gamma_1$ fixes 1 (resp.~no) element(s), while $\gamma_0$ fixes 
none. 
\end{lem} 

\begin{proof} We explain the first paragraph. From 
$((\bg)\gamma_\infty^{-1})\gamma_0=(\bg)\gamma_1$ on reduced Nielsen classes, the range of
$\gamma_0$ and $\gamma_1$ are the same on any
$\gamma_\infty$ orbit. So, the
$\sh$-incidence matrix is the same as the matrix obtained by replacing
$\sh=\gamma_1$ by $\gamma_0$.

A fixed point of $\gamma_1=\sh$ in $O$,  a $\gamma_\infty$ orbit, would contribute to $O\cap
(O)\sh$. Since the $\sh$-incidence matrix is the same as that for replacing $\gamma_1$ by
$\gamma_0$, 0's along the diagonal also imply there is no $\gamma_0$ fixed point.  

We now show the statement about fixed points of $\gamma_1=\sh$. Any fixed points must come
from a nonzero entry along the diagonal of the $\sh$-incidence matrix. For
$\ni^{+}_0$,  there is precisely one reduced Nielsen class $\bg$ in $O_{1,1}\cap (O_{1,1})\sh$. Write
$\bg=(\bg')\sh$. Apply $\sh$ to both sides, and conclude $(\bg)\sh=\bg'$. Therefore, as there is only one
element with this property, $\bg=\bg'$.    Now return to the example details. 

Apply the above to $\ni^{-}_0$. Since  $|O_{1,4}\cap
(O_{1,4})\sh=2|$, there are either two fixed points, or none. Since $\sh$ preserves the fiber over [1], we need only  check if
$(\bg_{1,4})\sh$ is reduced equivalent to $\bg_{1,4}$. Apply $q_1^{-1}q_3$ to $(\bg_{1,4})\sh$: 
the result is $((1\,2\,3),(3\,4\,2), (1\,3\,4), (1\,2\,4))$. Conjugate this by $(1\,2\,3)^{-1}$ to get
$\bg_{1,5}$. So, $\gamma_1$ has no fixed points on $\ni^-_0$. 
Since $\gamma_0$ moves the fibers over $[1],[2],[3]$ in a cycle, it fixes no Nielsen class elements. 
\end{proof}

We know the degrees of $\bar \psi_o^{\pm}$ are
respectively 9 and 6.  Lem.~\ref{g0g1fix} gives the genus $g_0^{\pm}$ of $\bar
\sH^{\inn,\pm}_0$ from Riemann-Hurwitz: \begin{equation} \label{compA43gen} \begin{array}{rl}
2(9+g_0^+-1)&=3\cdot 2+ (9-1)/2 + (1+2+3)=16,
\text{ or } g_0^+=0;\\
2(6+g_0^--1)&=2\cdot 2+ 6/2 + 3=10,
\text{ or } g_0^-=0.\end{array} \end{equation}

\begin{rem} In the
$\bar M_4$ orbit on $\ni^{\inn,-}_0$ there is a nonzero diagonal entry, though neither
$\gamma_0$ nor
$\gamma_1$ has a fixed point in the corresponding $\gamma_\infty$ orbit.  \end{rem}

\subsubsection{Checking $s_{\hat A_4/A_4}$ of \S\ref{smallInv} on two $\bar M_4$ orbits} 
Apply $\sh$ to $4'$. This shows $g_{1,4},g_{1,5},4',5'$ all lie in one $\bar M_4$ orbit. Any H-M rep.~has
lifting invariant +1, and since it is a $\bar M_4$ invariant, all elements in $\ni^+_0$ have lifting
invariant +1. For the other orbit, we have only to check the lifting invariant on $4'$, written in 
full as $$\bg_{1,4}=((1\,2\,3),(1\,2\,4),(1\,2\, 4),(4\,3\,2))=(\row g 4).$$ 
Compute the lifting invariant as $\hat g_1\hat g_2 \hat g_3\hat g_4$. Since $g_2=g_3$ (and their lifts
are the same), the invariant is $\hat g_1 \hat g_2^2 \hat g_4$. Apply Prop.~\ref{cnf.3} (not necessary,
though illuminating).  The genus zero hypothesis for a degree 4 cover holds for $((1\,2\,3),(1\,4\,2),
(4\,3\,2))$: $$s_{\hat A_4/A_4}(\bg_{1,4})=(-1)^{3\cdot (3^2-1)/8}=-1.$$  

\subsubsection{Why $\sH_0^{\pm}$ aren't modular curves} 
From \S\ref{WohlfahrtThm}, if the degree nine cover is modular, the monodromy group of the cover is a
quotient of
$\PSL_2(12)$.   If the degree 6 orbit is modular, the monodromy group is a quotient of
$\PSL_2(4)$. Since $\PSL_2(\bZ/4)$ modular curve has the $\lambda$-line as a quotient, with 2,2,2 as the cusp lengths,
these cusp lengths are wrong for the second cover to correspond to the $\lambda$-line. Similarly, for the
degree nine cover, as
$\PSL_2(\bZ/12)$ has both $\PSL_2(\bZ/4)$ and $\PSL_2(\bZ/3)$ as a quotient, the cusp lengths are
wrong. 

We can check the length of a $\sQ''$ orbit on $\ni^{\inn,+}_0$ and
$\ni^{\inn,-}_0$ by checking the length of the orbit of any particular element. If an orbit has an
H-M rep.~like $\bg_{1,1}$ it is always convenient to check elements of $\sQ''$ on it:
\begin{equation} \label{Q''check} \begin{array}{rl} (\bg_{1,1})\sh^2=&(1\,3)(2\,4)\bg_{1,1}(1\,3)(2\,4)\text{ and;} \\
(\bg_{1,1})q_1q_3^{-1}=&(1\,3)\bg_{1,1}(1\,3).\end{array} \end{equation}  

The top line of \eqref{Q''check} says $\sh^2$ fixes $\bg_{1,1}$. The bottom line, however, says 
$(\bg_{1,1})q_1q_3^{-1}$ is absolute, but not inner equivalent to $\bg_{1,1}$. For $\ni^{\inn,-}_0$, $\bg_{1,4}$ is transparently
fixed by
$\sh^2$, and  
$(\bg_{1,4})q_1q_3^{-1}=(3\,4)\bg_{1,4}(3\,4)$. Conclude the orbit length of $\sQ''$ on both
$\ni^{\inn,+}_0$ and
$\ni^{\inn,-}_0$ is 2. 

We finish Prop.~\ref{A4L0} by producing the map $\beta$ in  \eqref{jb}, and thereby concluding Prop.~\ref{appAndre}. Each $\bp\in \sH(\tilde
F_{2,2}/\Phi^1\xs J_3,\bfC_{\pm 3^2})^{\abs,\rd}$ gives a degree 4 cover $\phi: X_\bp\to \prP^1_z$ with four 3-cycle branch points. From
R-H, the genus $g$ of $X_\bp$  satisfies $2(4+g-1)=8$, or $g=1$. It may not, however, be an elliptic curve, though its degree 0 Picard variety
$\Pic(X_\bp)^{(0)}$ is. Define $\beta$ by taking its $j$-invariant.   

\begin{prop} \label{appAndre} The absolute space $\sH(\tilde F_{2,2}/\Phi^1\xs J_3,\bfC_{\pm
3^2})^{\abs,\rd}$ at level 0 embeds
in $\prP^1_j\times \prP^1_j$, but is not a Modular curve. So,  Andr\'e's Thm.~\cite{An98} says it contains at most finitely many Shimura-special
points (unlike  the $J_2$ case). \end{prop}

\begin{guess} The conclusion of Prop.~\ref{appAndre} is true for all other $p\ne 3$. \end{guess} Yet, we have a problem: What does 
{\sl Shimura special\/} mean when $p\ne 2$ or 3?

\subsubsection{Level 1 of $(A_4,\bfC_{\pm 3^2},p=2)$} \label{level1p2}  Level 1 of the \MT\ covers $\sH_0^+$:  $$\sH(\tilde F_{2,2}/\Phi^2\xs
J_3,\bfC_{\pm 3^2})^{\inn,\rd}\to \sH_0^+.$$ We know level 1 has two genus 0 components, 
$\sH_1^{-,c},\sH_1^{-,c'}$, complex conjugate and {\sl spin\/} obstructed; two 
genus 3 components, $\sH_1^{+,3}, \sH_1^{-,3}$, one spin obstructed, the other obstructed
by another Schur multiplier; and two genus 1 components, $\sH_1^{+,\beta}$, $\sH_1^{+,\beta^{-1}}$
both H-M comps \cite{FrS05}. 

Significance of $\sH_1^{+,\beta}$, $\sH_1^{+,\beta^{-1}}$:  
\begin{edesc} \label{apps} \item \label{appsa} $\text{Out}(\tilde F_{2,2}/\Phi^2\xs J_3)$ conjugates
$\sH_1^{+,\beta}$ to
\!$\sH_1^{+,\beta^{-1}}$\!\!.  
\item  \label{appsb} The following are equivalent for $K\le \bR$ a number field \cite[Ex.~9.2]{BFr02}. 
\begin{itemize} \item There are $\infty$-ly 
 many (reduced inequivalent -- \S\ref{jinv}) 4 branch point, $K$ regular realizations of the  
2-Frattini extension $G_1(A_5)$ of $A_5$.  
 \item $\sH_1^{+,\beta}$ has $\infty$-ly many $K$ points. 
\end{itemize} 
\end{edesc} 


\begin{appendix} 
\renewcommand{\labelenumi}{{{\rm (\teql \alph{enumi})}}} 

\section{Nielsen classes for $F_2\xs \bZ/2$}  \label{F2Z2} 
\S\ref{MTmc} does the Nielsen class version of all modular curves, by considering them coming from a rank 2 MT. Prop.~\ref{H2NC} shows 
there is a unique limit group $(\bZ_p)^2\xs \bZ/2$ \wsp not the whole universal $p$-Frattini cover \wsp for each $p\ne 2$. Then,
\S\ref{heisenberg} shows the Heisenberg group kernel acts here as a universal obstruction, running over all odd $p$.  

\subsection{Limit groups for the rank 2 MT of modular curves} \label{MTmc} Following \S\ref{F-q}, we consider the nonempty Nielsen classes of
the form $\ni(V\xs \bZ/2,\bfC_{2^4})$, $V\in \sV_p'$ (a nontrivial  $\tilde F_{2,p}$ quotient on which $\bZ/2$ acts, as in
\S\ref{setUpSM}).  The following formalizes an argument of
\cite[p.~114]{Fr4}. Form the projective completion of
$$K_4=\lrang{\psigma=\row
\sigma 4\mod \sigma_1\sigma_2\sigma_3\sigma_4=1 \text {
(product-one)}}.$$  Denote the result by $\hat  K_\psigma$. Use the
notation of
\S\ref{gpStatements}. 

\begin{prop} \label{H2NC}  Let $\hat D_\psigma$ (compatible with 
Cor.~\ref{weigCor1}) be the  quotient of $\hat  K_\psigma$ by   
$$\sigma_i^2=1,\   i=1,2,3,4\ (\text{so\
}\sigma_1\sigma_2=\sigma_4\sigma_3).$$ Then, $\prod_{p\ne 2} \bZ^2_p\xs 
J_2\equiv 
\hat D_\psigma$ and $\bZ^2_p\xs J_2$ is the unique $\bfC_{2^4}$
$p$-Nielsen class limit. 

The component graph of $\sC^f_{(\bZ/p)^2\xs\bZ/2,\bfC_{\pm 3^2},p}(\bZ^2_p\xs 
J_2)$ (as in \S\ref{extGraphs}) is a principle homogeneous space for $G(\bQ^\cyc/\bQ)$. 
\end{prop} 

\begin{proof} We show $\hat D_\psigma$ is  
$\tilde \bZ^2\xs J_2$; $\sigma_1\sigma_2$ and 
$\sigma_1\sigma_3$ are generators of $\tilde\bZ^2$; and  then that 
$\sigma_1$  acts on $\tilde\bZ^2$
by multiplication by $-1$. 

First: 
$\sigma_1(\sigma_1\sigma_2)\sigma_1=\sigma_2\sigma_1$ shows
$\sigma_1$ conjugates $\sigma_1\sigma_2$ to its inverse. Also, 
$$(\sigma_1\sigma_2)(\sigma_1\sigma_3)\!=\! 
(\sigma_1\sigma_3)\sigma_3(\sigma_2\sigma_1)\sigma_3\!=\!
(\sigma_1\sigma_3)(\sigma_1\sigma_2)$$ 
shows the said generators commute. The maximal possible quotient is 
$\bZ_p^2\xs\{\pm1\}$. 

Now we show for $G=V\xs
J_2$, $V$ a nontrivial quotient of $\bZ_p^2$, that $\ni(G,\bfC_{2^4})$ is nonempty.   Use 
a cofinal family of $V\,$s, $(\bZ/p^{k+1})^2$, $p\ne 2$.
Two proofs, one pure Nielsen class, the other with elliptic curves, appear in  \cite[\S6.1.3]{Fr05}. That shows
$\bZ_p^2\xs\{\pm1\}$ is a limit group. 

Uniqueness of the limit group does follow if we know there is just one braid orbit on the respective inner Nielsen classes. Alas, that
isn't  so. 

To finish  we use absolute Nielsen classes as an aid. Apply the elementary divisor theorem
to $(\bZ_p)^2$: Up to change of basis we may assume $V=\bZ_p/p^{u_1}\times \bZ_p/p^{u_2}$ with $u_1\le u_2$. 
If $u_1=0$, \cite[p.~156]{FrGGCM} shows there is just one braid orbit: in agreement with 
 identifying $\sH(D_{p^{u_2+1}},\bfC_{\pm 3^2})^{\inn,\rd}$ with the irreducible modular curve $Y_1(p^{k+1})$. 

This argument also applies to the general
case to reduce to when   $u_1=u_2$. That case is  the first two paragraphs of the proof of  \cite[Prop.~6.3]{Fr05}. Its essential gist, where
$\abs$ refers to modding out by $\GL_2(\bZ/p^{u+1})$ on Nielsen classes: 
\begin{edesc} \label{elmod} \item \label{elmoda} There
is just one element in $\ni((\bZ/p^{u+1})^2\times\bZ/2,\bfC_{\pm 3^2})^{\abs,\rd}$;
\item \label{elmodb} each of the $\phi(p^{u+1})/2$ inner classes 
defines a unique component of $\sH((\bZ/p^{u+1})^2\times\bZ/2,\bfC_{\pm 3^2})^{\inn,\rd}$; and 
\item \label{elmodc} the classes of \eql{elmod}{elmodb} are conjugate under the action of $G(\bQ(e^{2\pi i/p^{u+1}})/\bQ)$. \end{edesc} With 
$u$ varying this gives the last statement of the result.   
\end{proof}

\begin{rem}[Comments on \eql{elmod}{elmodb} and \eql{elmod}{elmodc}]  \label{SeOIT} Use the notation above. Excluding multiplication by -1, the
outer automorphisms 
$(\bZ/p^{k+1})^2\xs (\bZ/p^{k+1})^*$ of $(\bZ/p^{k+1})^2\xs \{\pm 1\}$ act through $\GL_2/\SL_2$ on $(\bZ/p^{k+1})^2$. By contrast the $H_4$
action is through $\SL_2(\bZ/p^{k+1})$ (explicitly in the proof). That is why you can't braid between The components of
$\sH((\bZ/p^{k+1})^2\xs\bZ/2,\bfC_{2^4})^{\inn,\rd}$. Yet, they form a single orbit under 
$G(\bQ(\cos(2\pi/p^{k+1}))/\bQ)$. This is the Hurwitz space interpretation of the {\sl Weil pairing}.  

The group $(\bZ/p)^2\xs J_2$ has quotients of the form $\bZ/p\xs J_2=G^*$. Corresponding to
that
$\bZ^2_p\xs J_2$ has the universal $p$-Frattini cover $\bZ/p\xs J_2$ of $G^*$ as a quotient. This is the source of the complex
multiplication situation in Serre's OIT (\S\ref{secondMC}). 
\end{rem} 

\subsection{Heisenberg analysis of modular curve Nielsen classes} \label{heisenberg} We briefly remind the reader of Loewy layers
and apply Jenning's Thm.~in \S\ref{loewyLayers}. Then, \S\ref{heisObst} applies this to explain a universal obstruction from a Heisenberg
group.  

\subsubsection{A Loewy layer example} \label{loewyLayers} \cite[p.~3]{Ben1} explains Loewy layers of a $\bZ/p[G]$ module $M$, though with no
examples. Most readers won't realize they are almost always hard to compute (if $p||G|$). 

Let
$J_{G,p}=J$ be the intersection of the maximal left (or right) ideals of $\bZ/p[G]$: The Jacobson radical of $\bZ/p[G]$. The basic lemma is that
$M/J_{G,p}M$, the {\sl first\/} Loewy layer of $M$,  is the maximal semi-simple quotient of $M$ for the action of $G$. Then, to continue the
series inductively apply this with $J_{G,p}M$ replacing $M$. 

Usually, however, this is far less information than you want. 
\cite[Part II]{Fr4} is where I needed modular representations for the first time. This explains the following point: Knowing $M$ from its Loewy
layers requires adding info on the nonsplit subquotients $M'$ of $M$ of the  form $0\to S_1\to M' \to S_2\to 0$ with $S_1$
(resp.~$S_2$) irreducibles in the
$\ell+1$st) (resp.~$\ell$th)  Loewy layer. An arrow from the
$\ell+1$st at $S_1$ to a copy of $S_2$ in the $\ell$th Loewy layer represents $M'$. These arrows give (anti-)directed paths from layer 1 to any
other layer $\ell$. 

For $G$ a $p$-group, and $M=\bZ/p[G]$, $J$ is the augmentation ideal: $$\ker(\sum_{g\in G} a_g g\mapsto 
\sum_{g\in G} a_g).$$ Jenning's Thm.~\cite[Thm.~3.14.6]{Ben1} (based on
\cite{quillen}) gives Loewy layer dimensions  with a Hilbert polynomial $H_{G}(t)$ (variable $t$). The
only
$p$-group irreducible is
$\one_G$. So, add the Loewy arrows from levels $\ell+1$ to $\ell$ and we know everything. 

\newcommand{\jF}{{F^\dagger}} 
Let $\jF_u(G)=\{g\in G\mid g-1\in J^u\}$. So, $\jF_1(G)=G$. Then, the input for $H_{G}(t)$ consists of the dimensions $n_1,n_2,\dots, n_u,\dots$
of the graded pieces of a Lie algebra due to Jenning's. The $u$th graded piece is $\jF_u/\jF_{u+1}$. Part of the proof shows $\jF_u$ is
generated by commutators and $p$th powers from $F\,$s with lower subscripts. In particular, if $G=(\bZ/p)^n$, then $n_1=n$ and $\jF_u/\jF_{u+1}$
is trivial for $u\ge 2$. So, the general expression $\prod_{u\ge 1} (\frac{1-t^{pu}}{1-t^u})^{n_u}$ becomes just
$H_{(\bZ/p)^n}(t)=(\frac{1-t^{p}}{1-t})^{n}$. 

\begin{lem} \label{loewyzp2} Then, $H_{(\bZ/p)^2}(t)=(1 + t + \ldots + t^{p-1})^2$ and the respective Loewy layers of $\bZ/p[(\bZ/p)^2]$ have
the dimensions
$1,2,\dots,p,p-1,\dots, 1$.  Given generators $x_1,x_2$ of the $\bZ/p$ module $(\bZ/p)^2$,  the symbols
$x_1^{\alpha}x_2^{\ell-\alpha}$,
$0\le \alpha,\ell-\alpha <p$ represent generators of copies of 
$\one$ at Loewy layer
$\ell$.  Arrows from $\one$ associated to  $x_1^{\alpha}x_2^{\ell-\alpha}$  go to copies of $\one$
associated to $x_1^{\alpha}x_2^{\ell-1-\alpha}$ and to $x_1^{\alpha-1}x_2^{\ell-\alpha}$ under the above constraints.     
\end{lem} 

\begin{proof} Calculate the coefficients of $(1 + t + \ldots + t^{p-1})^2$ to see the numerical series correctly expresses
the dimensions. The Loewy arrow statements come from identifying those subquotients of  $R=\bZ/p[G]$ that are module extensions of $\one$ by
$\one$. For this use  the Poincar\'e-Birkoff-Witt basis for the universal enveloping algebra of  $R$ \cite[p.~88]{Ben1}. \end{proof}

\subsubsection{A Heisenberg obstruction} \label{heisObst} 
The situation of Prop.~\ref{H2NC} is an example of
Cor.~\ref{weigCor1}.  First, $(\bZ \times \bZ) \xs \bZ/2$ is an oriented $p$-Poincar\'e 
duality group if $p$ is odd:  the finite-index subgroup $\bZ \times
\bZ$ is a surface group (the fundamental group of the torus).  Denote the matrix $\begin{pmatrix} 1 &x &z \\ 0& 1 &y \\ 0& 0& 1\end{pmatrix}$
by
$M(x,y,z)$ and consider 
$$\bH_{R,3}=\{M(x,y,z)\}_{x,y,z\in
R},$$ the {\sl  Heisenberg group\/} with entries in the commutative ring $R$. Let $H\le S_n$. Then, there is a 1-dimensional $\bZ/p[S_n]$ (so
also a $\bZ/p[H]$) module whose action is $m\mapsto (m)g=(-1)^{\Det(g)}m$. Denote $M$ by $\one^-$. This extends to a $\bZ_p[H]$ action on
$\bZ_p$. Denote this module as
$\bZ_p^-$.  

In our usual notation, let $G_0=(\bZ/2)^2\xs \bZ/2$ and  denote the 1st characteristic $p$-Frattini
cover of $G_0$ by $G_1$.  Prop.~\ref{univHeis} uses a universal Frattini extension. It specializes for all odd primes
$p$ to the
$\bZ/p$ quotient obstructing (as in Def.~\ref{obstMT}) the unique braid orbit in  $\ni(G_0,\bfC_{2^4})$ from lifting to $\ni(G_1,\bfC_{2^4})$,
as in Cor.~\ref{weigCor1}. In fact, by pullback we see it as the limit group obstruction in Cor.~\ref{weigCor2}. 

\begin{prop}\label{univHeis} The map $\bH_{\bZ/p,3}\to (\bZ/p)^2$ by $M(x,y,z)\mapsto (x,y)$ is a Frattini extension. The $p$-Frattini module
$M_0(G_0)$ of $G_0$ has $\one_{G_0}\oplus \one_{G_0}^-\oplus \one_{G_0}^{-}$ at its head. The extension defined by
$\one_{G_0}$  gives the Heisenberg group, obstructing the \MT\ at level 1. Still, it gives an infinite limit group 
$(\bZ_p )^2
\xs
\bZ/2$ by regarding $\bZ_p \times \bZ_p$ as $\bZ_p^- \times \bZ_p^-$.\end{prop}

\begin{proof} The characteristic Frattini  cover $\psi_{1,0}: G_1((\bZ/p)^2)\to (\bZ/p)^2$  factors through 
$\psi_\ab=(\bZ/p^2)^2\to (\bZ/p)^2$ (modding out by $p$). The nontrivial element of $\bZ/2$ acts by multiplication by $-1$ on $(\bZ/p^2)^2$.
In fact, $\psi_\ab$ is the maximal abelian extension through which $\psi_{1,0}$ factors. 

Loewy layers of any $(\bZ/p)^2\xs \bZ/2$
module are copies of $\one$ and $\one^-$. So, any proper extension of $\psi_\ab$ through which $\psi_{1,0}$ factors, also  factors 
through
$\psi': H\to (\bZ/p)^2$ with $\ker(\psi')$ of dimension 3 and $H$ not abelian.  

We choose the Heller construction (in \cite[Part II]{Fr4}, for example) to  describe the characteristic module
$$M_0((\bZ/p)^2\xs
\bZ/2)=\ker(G_1((\bZ/p)^2\xs \bZ/2)\to (\bZ/p)^2\xs \bZ/2) (p \text{ odd}).$$ Here is the rubric for this simple, though still nontrivial
case. Suppose $G_0$ is $p$-split: $G_0=P^*\xs H$ with $(|H|,p)=1$ and $P^*$ the $p$-Sylow, as in our case. Use the
Poincar\'e-Birkhoff-Witt basis of the universal enveloping algebra (from the proof of Lem.~\ref{loewyzp2}) to deduce the action of $H$ from its
conjugation action on
$P^*$. In our case, the $\ell$th Loewy layer of $\bZ/p[P^*]\eqdef P_\one$, with $P^*=(\bZ/p)^2$ consists of sums of
$\one$ (resp.~$\one^-$) if $\ell$ is even (resp.~odd) from 0 to $2p-2$ (resp.~1 to $2p-1$).  That is the projective indecomposable module
for $\one$. 

Now list the Loewy display for the projective indecomposable modules for
$G_0$ by tensoring the Loewy layers of the projective indecomposables for $\one$ with the semi-simple modules for $H$ 
\cite[p.~737]{darren2}. In our case, the semi-simples for $\bZ/2$ are just $\one$ and $\one^-$ giving $P_\one$ and $P_{\one^-}$ as the
projective indecomposables, the latter having the same look as the former except you switch the levels with $\one$ with those with $\one^-$.
Finally, 
$M_0$ is $\Omega_2\eqdef \ker(\psi_2: P_{\one^-}\oplus P_{\one^-}\to \ker (P_\one\to \one))$ with this understanding: $\ker (P_\one\to \one)$
has at its head $\one^-\oplus \one^-$ and $\psi_2$ is the map from the minimal projective ($P_{\one^-}\oplus P_{\one^-}$) that maps onto $\ker
(P_\one\to
\one)$. 

Using the arrows between Loewy layers that appear in Lem.~\ref{loewyzp2}, we can be explicit about constructing $\psi_2$ (knowing the
result is independent of our choices). For example, map the first copy of $P_{\one^-}$ in $P_{\one^-}\oplus P_{\one^-}$ so the image
$P'$ has
$\one^-$ at its head coming from the 3rd layer of $P_{\one^-}$. 

Then, map the second copy of $P_{\one^-}$ in $P_{\one^-}\oplus
P_{\one^-}$ to see the head of the image in  $\ker (P_\one\to \one)/P'$ is $\one\oplus \one^{-}$. These summands come from the
respective 2nd and 3rd Loewy layers of the copy of $P_{\one^-}$. That concludes the head of $M_0$.
The rest follows by identifying $\bH_{\bZ/p,3}\xs \bZ/2$ with the quotient of $G_1$ that extends $G_0$ by $\one_{G_0}$. 
\end{proof} 

\section{Nielsen classes for $F_2\xs \bZ/3$} \label{F2Z3} \S\ref{3rdMC} used the $p=2$ case of the \MT\ with $\bZ/3$ acting on
$F_2$. \S\ref{Z3limgps} gives our present knowledge of limit groups here. Finally, Ex.~\ref{schurComps} shows the effect of Schur multiplier statements
from  
\S\ref{schurDisc}: They account for much, but not all, of the six level 1 components for the case $p=2$. \S\ref{Z3compMult} gives 
a meaning to complex multiplication by considering the F-quotient from \S\ref{F-q} when $p\equiv 1 \mod 3$.

\subsection{Limit groups for another rank 2 \MT} \label{Z3limgps}

The next result works by proving the existence of H-M reps.~(whose shift gives example g-$p'$ cusps as in Ex.~\ref{HMbranch}). So, this 
produces $\tilde F_{2,p}\xs J_3$ as a limit group for each $p\ne 3$ from  Princ.~\ref{FP2}. 

Recall the action of $\alpha$ from \eql{F2Zn}{F2Zn3}. It induces the
matrix $\smatrix{0}{-1}1{-1}$, with characteristic polynomial $x^2+x+1$, on the $(\bZ/p)^2$ quotient of $F_2$. Denote $F_2/(F_2,F_2)$ by $L_2$
and its completion at $p$ by $L_{2,p}$.  

\begin{prop} \label{H3NC} The $(A_4,\bfC_{\pm 2})$ \MT\ for $p\ne 3$ has $\tilde F_{2,p}\xs J_3$   as a 
 limit group. For $p=2$,  the $(A_4,\bfC_{\pm 2})$ \MT\ also has $L_{2,p}\xs J_3$ as a  limit group. 
\end{prop} 

\begin{proof} Let $G=G_p=(\bZ/p)^2\xs J_3$:
$\lrang{\!\alpha}=J_3$.  We find  $$g_1=(\alpha,\bv_1) \text{ and }
g_2=(\alpha,\bv_2)$$ so that $\lrang{g_1,g_2}=G$.  The H-M rep.~$(g_1,g_1^{-1},g_2,g_2^{-1})$  is in $\ni(G,\bfC_{\pm 3^2})^{\inn}$.
Conjugate in $G$, to take a representative in the inner class with $\bv_1={\pmb 0}$. Consider 
$$g_1g_2^{-1}=(1,-\bv_2)\text{ and }
 g_1^2g_2=(1,\alpha^{-1}(v_2)).$$ So, $g_1,g_2$ generate precisely
when $\lrang{-\bv_2, \alpha^{-1}(v_2)}=(\bZ/p)^2$. Such a $\bv_2$ exists because 
the
eigenvalues of
$\alpha$ are distinct. So $(\bZ/p)^2$ is a cyclic $\lrang{\alpha}$
module. 

Now consider $\ni(G,\bfC_{2^4})^\inn$ with  $G=U\xs J_3$ and $U$ (a $\bZ/3$ module) having $(\bZ/p)^2$
as a quotient. There is a surjective map
$\psi: G\to (\bZ/p)^2\xs J_3$: a Frattini cover. So, if $g_1',g_2'$ 
generate $(\bZ/p)^2\xs J_3$, then respective order 3
 lifts of $g_1',g_2'$ to
$g_1,g_2\in G$ automatically generate $G$. Princ.~\ref{FP2} now applies: For $p\ne 3$, an H-M cusp branch gives $\tilde F_{2,p}\xs J_3$  
as a 
 limit group. 

Now we turn to the case $p=2$, and consider the other, not H-M rep., braid orbit on $\ni(A_4,\bfC_{\pm 3})$ given in Prop.~\ref{A4L0}.
\cite[Cor.~5.7]{BFr02} gives this Loewy display for $M_0=\ker(G_1(A_4)\to A_4)$: $0\to U_3\to U_3\oplus \one$ with $U_3$ the 2-dimensional
irreducible for $\bZ/2[A_4]$. In the augmented Loewy display, there is an arrow from the leftmost $U_3$ to each summand of $U_3\oplus
\one$. 

Let
$\bg$ be a representative of the orbit $\ni^-_0$ obstructed by
$\hat A_4\to A_4$. The completion at $p=2$ of the quotient $F_2/(F_2,F_2)\xs \bZ/3$ is $L_{2,2}\xs \bZ/3$, a 2-Frattini cover of $A_4$. Notice
that $\one_{A_4}$ is not a subquotient in this group. Therefore, Cor.~\ref{weigCor2} implies the map $M_\bg\to A_4$ extends to $M_\bg\to
L_{2,2}\xs \bZ/3$. Indeed, it is a Nielsen limit group through the braid orbit of $\bg$. 
\end{proof}

\begin{exmp}[The $(A_5,\bfC_{3^4},p=2)$ \MT] \label{A5limitgps} We continue Ex.~\ref{level1A5}. Let $O_2$ be the non-H-M braid
orbit of $\ni(G_1(A_5),\bfC_{3^4})$. \cite[Prop.~9.14]{BFr02} shows 
$G_1(A_5)$  embeds in $A_N$ for several values of $N$ (40, 60, 80, 120) with an additional property: With $$\Spin_N\times_{A_N} G_1(A_5)\eqdef
\Spin_N'\to G_1(A_5),$$ we have 
$s_{\Spin_N'}(O_2)=-1$.  

Let   $R_k'\to G_k$ be the
$k\nm 1$st  antecedent to $\Spin_N'\to G_1(A_5)$ (\S\ref{repComps}). As noted in Ex.~\ref{highObstExs}, the hypotheses of Thm.~\ref{highObst}
hold for this example and each level $k\ge 1$ of the \MT\ has an H-M component with at least two distinct
limit groups. 
\end{exmp}

\begin{exmp}[$\ni(G_1(A_4), \bfC_{\pm 3^2})$ braid orbits]  \label{schurComps} Again $p=2$. Similar to Ex.~\ref{A5limitgps},
and again using Ex.~\ref{highObstExs},  each level $k\ge 2$ has two H-M components, and each such component has at least four  distinct
limit groups.    
\end{exmp} 

\begin{prob} Let $\sH_k'$ be one of the H-M components in Ex.~\ref{schurComps}. Is the number of limit groups through $\sH_k'$ bounded with
$k$? \end{prob} 
 
\subsection{Complex multiplication for the $\bZ/3$ case} \label{Z3compMult}  Use the notation above.  If $p\ne 2,3$, 
$\alpha$ on
$(\bZ/p)^2$ has eigenvalues defined over $\bZ/p$ precisely when $-3$ is a square $\!\!\mod p$. From quadratic reciprocity, these are the
$p\equiv 1\mod 3$. 
 Exactly then,  
$\tilde F_{2,p}\xs J_3$ has quotients of the form
$\bZ/p\xs J_3=G^*$. Corresponding to that, the universal $p$-Frattini cover $\bZ_p\xs J_3$ of $G^*$ is a quotient of 
$\bZ^2_p\xs J_3$.  

\begin{prob} \label{NewSeOIT} When $p\equiv 1 \mod 3$, does a $\bZ_p\xs J_3$ quotient of $\tilde F_{2,p}\xs J_3$ correspond
to \lq\lq complex multiplication case\rq\rq\ for special values $j'\in \prP^1_j$ (as in the $J_2$ case in Rem.~\ref{SeOIT})? For all $j'\in
\prP^1_j$ over a number field, does this give a full analog of Serre's OIT in the $J_3$ case? \end{prob} 

The nontrivial F-quotient when $p\equiv 1\mod 3$ is like that for modular curves, a \MT\ case where $M_k'=\ker(G_{k+1}(\bZ/p\xs J_3)\to
G_k(\bZ/p\xs J_3))$ has rank 1 (as in Prop.~\ref{densityOnes}). What we know of $M_k=\ker(G_{k+1}((\bZ/p)^2\xs J_3)\to
G_k((\bZ/p)^2\xs J_3))$ (as a $G_k$ module, to which the conclusion of Prop.~\ref{densityOnes} applies) is from Semmen's thesis \cite{darren2}.
Such  information is significant in analyzing \eql{serreAnal}{serreAnalb}.

\section{Related Luminy talks and typos from \cite{BFr02}} Other Luminy talks contain material whose perusal simplifies our explaining the use
of the Hurwitz monodromy group and the background for this paper \S\ref{Lumtalks}. Our approach to explaining progress on \MT s is to use
\cite{BFr02} as a reference book in translating between geometric and arithmetic statements until the completion of
\cite{FrBook}. Our web site version of the former has typos corrected as they appear. 

\subsection{Conference talks that explain significant background points} \label{Lumtalks} Expositional elements
of the following papers support their use in \MT s.  

\begin{itemize} \item Matthieu Romagny and Stefan Wewers introduced 
Nielsen classes and material on Hurwitz spaces.
\item Kay Magaard introduced braids acting (through Hurwitz monodromy  $H_r$; \S\ref{NielClDict}) on Nielsen classes, necessary for
computations.
\item Pierre D\`ebes defined a (rank 0) Modular Tower (\MT),
comparing that with modular curves.
\item The (weak; rank 0) Main Conjecture is that there are no rational
points at suitably high tower levels. Pierre's talk reduced this conjecture, for
four branch point towers, to showing the genus rises with the levels.
\item Darren Semmen presented the profinite Frattini category. This 
showed how Schur multipliers control properties of the
Modular Tower levels. \end{itemize}

\subsection{Typos from the printed version of \cite{BFr02}} \label{BFrtypos} 
\begin{itemize} 
\item  p. 55, line 4 of 2nd paragraph: to the near H-M and H-M [not H -M] 
 p. 87, line 4. It also
explains  H-M [not H -M] and near H-M p. 87, line 8.  {\sl complements\/} of H-M and near H-M
[not H -M] p. 89, after (8.6): H-M or near H-M  [not H -M] rep.~is
p. 180, 3rd line of 2nd par.:  [not H -M] 

\item p. 92: It said: \lq\lq The cusp pairing for $r=4$ 
should extend to the case $r\ge 5$, though we don't yet know how.\rq\rq\  

We knew how to do that by the time the paper was complete, though we forgot to delete this line.
It now says: \lq\lq The cusp pairing for $r=4$  extends to the case $r\ge 5$ (\S2.10.2).\rq\rq  

\item p. 93: 1st par. \S1.4.7 (end): Change Merel-Mazur to Mazur-Merel. 
\item p. 94: (and image of $g_1^{-1}g_2$ in $A_5$ of order 5)

\item p. 103--104. Use of $\sQ''$ in Def.~2.12 on p.~103 precedes its definition on p.~104. 

\item Bottom of p. 107: $|\ni_k^\inn|= (p^{k+1}+p^k)\phi(p^{k})/2$ should be 
$$|\ni_k^\inn|= (p^{k+1}+p^k)\phi(p^{k+1})/2.$$

\item Statement of Prop.~2.17. [States the condition $o(g_1,g_2)$ is odd, after it uses that condition.] It should say this. 
Let $g_1g_2=g_3$, and $g_2g_1=g_3'$. Let $o(g_1,g_2)=o$ (resp.~$o'(g_1,g_2)=o'$) be the length of the orbit of $\gamma^2$ (resp.~$\gamma$) on
$(g_1,g_2)$. If $g_1=g_2$, then $o=o'=1$. 

\noindent Proposition 2.17 {\sl Assume $g_1\ne g_2$.  The orbit
of 
$\gamma^2$ containing $(g_1,g_2)$ is $(g_3^j g_1 g_3^{-j},g_3^jg_{2} g_3^{-j})$, 
$j=0,\dots, \ord(g_3)\nm1$. So,   
$$o=\ord(g_3)/|\lrang{g_3}\cap Z(g_1,g_2)|\eqdef  o(g_1,g_2).$$  Then, $o'=2\cdot o$, unless  $o$  is odd, and with
$x=(g_3)^{(o-1)/2}$ and $y=(g_3')^{(o-1)/2}$ (so $g_1y=xg_1$ and  $yg_2=g_2x$),  $yg_2$ has order 2.  Then,  
$o'=o$.}

\item p. 129: Title of Section 4 should be: [Moduli] and reduced Modular Towers (change \lq\lq Modular\rq\rq\ to
\lq\lq Moduli\rq\rq). 

\item p. 140: Reference to [Fr01] changed to [Fr02]: and a more precise quote:
[Fr02, Prop. 2.8]:  M. Fried, Moduli of relatively nilpotent extensions, Institute of Mathematical
Science Analysis 1267, June 2002, Communications in Arithmetic Fundamental Groups,
70--94. 

\item p. 160, line 22: as $(u,v)=(\mp(\bg),\wid(bg))$ should be,  as $$(u,v)=(\mp(\bg),\wid(\bg)).$$ 
\item p. 172: 1st par.~of Prop.~8.12, change \lq\lq (and $g_1^{-1}g_2$ of order 5)\rq\rq\ to \lq\lq (and image of $g_1^{-1}g_2$ in $A_5$ of
order 5).\rq\rq 

\item p. 180: 1st line of 2nd paragraph of \S9:  Orbits of $\lrang{\gamma_1,q_2}$ should be [Orbits
of 
$\lrang{\gamma_1,\gamma_\infty}$], to emphasize here we view $q_2$ as in $\bar M_4$. 

\item Bottom p. 184: $G_{k+1}$ [acts] trivially on \dots

\item p. 188: Def. 9.11: ${}^T\hat H$ should be ${}^T \hat G$. 

\item Ex.~9.19: The 3rd sentence should be:  
For this case, $\tr(T_{H'}(m))=4=\tr{T_{H'}(m')}$ and $\tr(T_{H'}(mm'))=8$: 
$m,m' \in C_{18}$ and $mm'\in C_{16}$.
\end{itemize} 

\end{appendix}


\providecommand{\bysame}{\leavevmode\hbox to3em{\hrulefill}\thinspace}
\providecommand{\MR}{\relax\ifhmode\unskip\space\fi MR }
\providecommand{\MRhref}[2]{%
  \href{http://www.ams.org/mathscinet-getitem?mr=#1}{#2}
}
\providecommand{\href}[2]{#2}

\end{document}